\documentclass[10pt]{amsart}

\usepackage{xypic,amscd,amssymb,latexsym}
\usepackage{amsmath}	
\usepackage{graphicx}
\usepackage{pdfsync}
\usepackage[breaklinks=true]{hyperref}

\usepackage{cancel}

\usepackage{mathrsfs}

\usepackage{url}

\usepackage{changepage}

\usepackage{pstricks}
\usepackage{pst-poly}

\usepackage{color}

\definecolor{pinky}{rgb}{1.0, 0, 1.0}

\begin{document}
\pagenumbering{arabic}

\newtheorem{theorem}{Theorem}[section]
\newtheorem{proposition}[theorem]{Proposition}
\newtheorem{lemma}[theorem]{Lemma}
\newtheorem{corollary}[theorem]{Corollary}
\newtheorem{remark}[theorem]{Remark}
\newtheorem{definition}[theorem]{Definition}
\newtheorem{question}[theorem]{Question}
\newtheorem{claim}[theorem]{Claim}
\newtheorem{conjecture}[theorem]{Conjecture}
\newtheorem{defprop}[theorem]{Definition and Proposition}
\newtheorem{example}[theorem]{Example}
\newtheorem{deflem}[theorem]{Definition and Lemma}

\def\qed{{\quad \vrule height 8pt width 8pt depth 0pt}}

\newcommand{\cplx}[0]{\mathbb{C}}

\newcommand{\fr}[1]{\mathfrak{#1}}

\newcommand{\vs}[0]{\vspace{2mm}}

\newcommand{\til}[1]{\widetilde{#1}}

\newcommand{\mcal}[1]{\mathcal{#1}}

\newcommand{\ul}[1]{\underline{#1}}

\newcommand{\ol}[1]{\overline{#1}}

\newcommand{\wh}[1]{\widehat{#1}}

\newcommand{\smallvectwo}[2]{
\left(
\begin{smallmatrix}
#1 \\
#2 \\
\end{smallmatrix}
\right)
}

\newcommand{\smallmattwo}[4]{
\left(
\begin{smallmatrix}
#1 & #2 \\
#3 & #4 \\
\end{smallmatrix}
\right)
}

\title[Three-dimensional quantum gravity from the quantum pseudo-K\"ahler plane]{Three-dimensional quantum gravity \\ from the quantum pseudo-K\"ahler plane}

\author{Hyun Kyu Kim}

\address{School of Mathematics, Korea Institute for Advanced Study (KIAS), 85 Hoegiro Dongdaemun-gu, Seoul 02455, Republic of Korea}

\email[H.~Kim]{hkim@kias.re.kr}

\thanks{This paper is to appear in Commun. Math. Phys. \quad DOI: 10.1007/s00220-023-04783-w \hspace{-50,65mm}}

\begin{abstract}
A new canonical Hopf algebra called the quantum pseudo-K\"ahler plane is introduced. This quantum group can be viewed as a deformation quantization of the complex two-dimensional plane $\mathbb{C}^2$ with a pseudo-K\"ahler metric, or as a complexified version of the well-known quantum plane Hopf algebra. A natural class of nicely-behaved representations of the quantum pseudo-K\"ahler plane algebra is defined and studied, in the spirit of the previous joint work of the author and I. B. Frenkel. The tensor square of a unique irreducible representation decomposes into the direct integral of the irreducibles, and the unitary decomposition map is expressed by a special function called the modular double compact quantum dilogarithm, used in the recent joint work of the author and C. Scarinci on the quantization of 3d gravity for positive cosmological constant case. Then, from the associativity of the tensor cube, and from the maps between the left and the right duals, we construct unitary operators forming a new representation of Kashaev's group of transformations of dotted ideal triangulations of punctured surfaces, as an analog of Kashaev's quantum Teichm\"uller theory. The present work thus inspires one to look for a Kashaev-type quantization of 3d gravity for positive cosmological constant.
\end{abstract}

\maketitle

\tableofcontents

\section{Introduction}
\label{sec:introduction}

Quantization of Teichm\"uller spaces of Riemann surfaces, first established independently by Kashaev \cite{Ka1} and by Chekhov and Fock \cite{F} \cite{CF}, has several interesting mathematical consequences like the construction of new families of unitary projective representations of mapping class groups of surfaces on Hilbert spaces, or 3-manifold invariants and even 3d Topological Quantum Field Theories. The original interest in the subject seems to be due to the possibility of its being an approach to the problem of 3 or 2+1 dimensional quantum gravity. Nevertheless, a direct relationship between the quantum Teichm\"uller theory and the 3d quantum gravity had not been established until a recent joint work of the author with Carlos Scarinci \cite{KS}. The natural theory that we quantized in that paper depends on the real parameter $\Lambda$ called the \emph{cosmological constant}; in fact, what matters seriously is only the sign of this constant, so there are three qualitatively different cases to consider. We observe that, in the case of negative cosmological constant, when the relevant 3d spacetime manifold is topologically $\mathbb{R} \times S$ for a punctured surface $S$, the phase space of the theory can be realized essentially as Fock and Goncharov's `symplectic double' \cite{FG09} of the Teichm\"uller space of the surface $S$, having a structure of a cluster variety. In particular, the Fock-Goncharov quantization \cite{FG09}, which is a vast generalization of the Chekhov-Fock quantization of Teichm\"uller spaces \cite{CF}, can be directly applied to yield a quantization of this theory in this case. Hence one might now say that the 3d quantum gravity for the negative cosmological constant is equivalent in some sense to a certain `double' version of the quantum Teichm\"uller theory. 

\vs

Meanwhile, in \cite{KS} we also established a satisfactory equivariant quantization for each of the remaining two cases for the different signs of the cosmological constant. The formulation of the quantization is similar to that in the negative cosmological constant case, but we use different special functions for the unitary operators associated to transformations of ideal triangulations of the punctured surface $S$.  As shall be explained in more detail in \S\ref{sec:Kashaev_quantization}, these unitary operators constitute the heart of all the above mentioned quantization stories; in particular, these unitary operators have yielded two new families of unitary representations of mapping class groups on Hilbert spaces, for the cases of non-negative cosmological constant. 

\vs

The special function used in the quantization of Teichm\"uller spaces and cluster $\mathscr{X}$-varieties \cite{FG09} \cite{Ka1} \cite{F} \cite{CF} as well as in the quantization of a moduli space of 3d gravity for the negative cosmological constant case \cite{KS}  is the \emph{non-compact quantum dilogarithm} function $\Phi^\hbar(z)$ of Faddeev and Kashaev \cite{FK} \cite{Faddeev} which is defined for a real parameter $\hbar>0$ as a meromorphic function on the complex plane given on the strip $|\mathrm{Im}(z)|<\pi(1+\hbar)$ as the formula
\begin{align}
\label{eq:Phi_hbar_first_time}
\Phi^\hbar(z) = \exp\left(-\frac{1}{4} \int_\Omega \frac{e^{-{\rm i} pz}}{\sinh(\pi p) \sinh(\pi \hbar p)} \frac{dp}{p} \right),
\end{align}
known already to Barnes \cite{B01} a hundred years ago, where $\Omega$ is the real line contour avoiding the origin by a small half circular detour above the origin. This function can be regarded as the \emph{limit} of a certain ratio of the \emph{compact quantum dilogarithm} \cite{FK} functions
\begin{align}
\label{eq:psi_q_first_time}
\psi^q(z) = \prod_{n=0}^\infty (1+q^{2n+1} z)^{-1},
\end{align}
defined for each complex number $q$ with $|q|<1$. More precisely, consider any complex number $h$ satisfying ${\rm Im}(h)\ge 0$ and ${\rm Re}(h)>0$. The function $\Phi^h(z)$ can be defined using the same formula as in the right hand side of eq.\eqref{eq:Phi_hbar_first_time} with $\hbar$ replaced by $h$, which gives a function on the strip $|{\rm Im}(z)| < \pi(1 + {\rm Re}(h))$, which in turn can be analytically continued to a meromorphic function. It is known that in case ${\rm Im}(h)>0$, $\Phi^h(z)$ can be expressed as a ratio of two compact quantum dilogarithm functions:
$$
\Phi^h(z) = \psi^{\exp(\pi{\rm i} h)}(e^z)/\psi^{\exp(-\pi{\rm i}/h)}(e^{z/h}).
$$
The ratio on the right hand side does not make sense when ${\rm Im}(h)=0$ because then $|e^{\pi {\rm i} h}|=1=|e^{-\pi {\rm i}/h}|$. However, the limit of this ratio as ${\rm Im}(h)\to 0$ makes sense, that is, the non-compact quantum dilogarithm $\Phi^\hbar$, for $\hbar \in \mathbb{R}_{>0}$, can be viewed as the limit
$$\Phi^\hbar(z) = \lim_{h\to \hbar} \psi^{\exp(\pi{\rm i} h)}(e^z)/\psi^{\exp(-\pi{\rm i}/h)}(e^{z/h}).$$

\vs

In the meantime, the special functions used in \cite{KS} for the positive cosmological constant case are the functions $\Phi^{\pm {\rm i}\hbar}$, which can be obtained by replacing $\hbar$ by $\pm {\rm i} \hbar \in {\rm i}\mathbb{R}$ in the contour integral formula in eq.\eqref{eq:Phi_hbar_first_time}; that is, we consider $\Phi^h$ when ${\rm Re}(h)=0$. However, the previous contour $\Omega$ does not work, as some poles of the integrand lie on this contour. Still, the same-looking contour integral formula
\begin{align}
\label{eq:Psi_hbar_first_time}
\Phi^{\pm {\rm i}\hbar}(z) = \exp\left(-\frac{1}{4} \int_{\Omega'} \frac{e^{-{\rm i} pz}}{\sinh(\pi p) \sinh(\pm \pi {\rm i}  \hbar p)} \frac{dp}{p} \right),
\end{align}
yields a well-defined function, where $\Omega'$ is the contour obtained by rotating $\Omega$ by a nonzero angle about the origin, either clockwise or counterclockwise depending on the sign of $\pm$; see \cite{KS} for more detail. In a sense, one can view the situation as the $\hbar$-parameter of Faddeev-Kashaev's $\Phi^\hbar$ being analytically continued to a complex parameter $h$, and then being specialized to the imaginary axis to yield the functions $\Phi^{\pm {\rm i} \hbar}$. However, it seems that while the function $\Phi^h$ for the case ${\rm Re}(h)>0$ was considered and studied by Kashaev and others, the limiting case ${\rm Re}(h)\to 0$ had not been much considered in the literature until \cite{KS}. In particular, as mentioned, if one wants to express the result not in terms of a limit of $\Phi^h$ for ${\rm Re}(h)>0$ but directly in terms of a contour integral, then the contour of integration has to be modified from $\Omega$.

\vs

Meanwhile, the function used in \cite{KS} for the zero cosmological constant case is the function $F_0(x,y)$ of two real variables $x,y$ defined by
\begin{align}
\label{eq:F_hbar_0_first_time}
F_0(x,y) = \exp\left( -\frac{y}{2\pi{\rm i}} \int_\Omega \frac{e^{-{\rm i} px}}{\sinh(\pi p)} \frac{dp}{p} \right),
\end{align}
which in particular does not involve the quantum parameter $\hbar$ at all. 

\vs

All three functions in the equations \eqref{eq:Phi_hbar_first_time}, \eqref{eq:Psi_hbar_first_time}, \eqref{eq:F_hbar_0_first_time} can in fact be described in a uniform manner, using the ring of generalized complex numbers; see \cite{KS}. In the meantime, one way of writing the latter two functions in a more familiar form is
\begin{align}
\label{eq:Phi_ihbar_as_ratio}
\Phi^{{\rm i}\hbar}(z) = \psi^{\exp(-\pi\hbar)}(e^z)/\psi^{\exp(-\pi/\hbar)}(e^{z/({\rm i}\hbar)}), 
\qquad 
\Phi^{-{\rm i}\hbar}(z)  = \ol{ \Phi^{{\rm i}\hbar}(\ol{z}) }^{-1},
\end{align}
\begin{align}
\nonumber
F_0(x,y) = (1+e^{x})^{y/(\pi{\rm i})},
\end{align}
where $\Phi^{\pm {\rm i}\hbar}(z)$ is now genuinely a ratio of two compact quantum dilogarithms, not just a limit of such.  In \cite{KS} we suggest to call these new functions $\Phi^{\pm {\rm i} \hbar}$ in eq.\eqref{eq:Psi_hbar_first_time} and eq.\eqref{eq:Phi_ihbar_as_ratio} the \emph{modular double compact quantum dilogarithms}. In the quantization of a moduli space of 3d gravity for the positive cosmological constant case \cite{KS}, these functions $\Phi^{\pm {\rm i}\hbar}$, which are the ratios of two honest compact quantum dilogarithm functions $\psi^{\exp(-\pi\hbar)}$ and $\psi^{\exp(-\pi/\hbar)}$ (when $\hbar>0$), are crucially used in the construction of the mutation intertwiners, yielding a new family of unitary representations of mapping class groups on Hilbert spaces.

\vs

It is a natural question to ask whether various established properties about the previously known unitary representations of mapping class groups coming from quantum Teichm\"uller theory involving the non-comapct quantum dilogarithm $\Phi^\hbar$ would also hold for the new representations coming from 3d quantum gravity involving new functions $\Phi^{\pm {\rm i}\hbar}$ and $F_0$. The present paper gives a positive answer to this question for one particular property of quantum Teichm\"uller theory discovered by I. B. Frenkel and the author \cite{FKi}, about an interesting connection between Kashaev's quantum Teichm\"uller theory \cite{Ka1} for genus zero surfaces and the representation theory of one of the most basic quantum groups called the \emph{quantum plane}. The quantum plane, denoted by $\mathcal{B}_q$, which a priori seems to have nothing to do with the Teichm\"uller theory of surfaces, is a Hopf $*$-algebra defined by
$$
\langle \, X^{\pm 1}, Y \, | \, XY = q^2 \, YX \, \rangle
$$
as an algebra over $\mathbb{C}$, for a complex paramater $q$ with $|q|=1$ (hence $q^*=q^{-1}$), with the coproduct
$$
\Delta X = X\otimes X, \qquad \Delta Y = Y \otimes X + 1 \otimes Y,
$$
with corresponding suitable antipode and counit, and the $*$-structure
$$
X^* = X, \quad Y^* = Y.
$$
In particular, it can be viewed as a Borel subalgebra of the well-known quantum group $\mathcal{U}_q(\frak{sl}(2,\mathbb{R}))$. We then studied a certain natural class of nicely behaved representations of the quantum plane algebra which was investigated by Schm\"udgen \cite{Sch92} and extended by Ponsot and Teschner \cite{PT}, where $X$ and $Y$ are represented as unbounded positive self-adjoint operators on a Hilbert space. Via natural constructions in the representation theory of Hopf algebras which we elaborate in \S\ref{sec:algebraic_story} and throughout the main text of the present paper, we showed that unitary intertwiners between these representations of the quantum plane recover the unitary operators for Kashaev's quantum Teichm\"uller theory. 

\vs

Here we briefly summarize this construction of \cite{FKi}. First, one natural irreducible representation $\pi$ of the quantum plane $\mathcal{B}_q$ on the Hilbert space $\mathscr{H} = L^2(\mathbb{R},dx)$ is studied, given by $\pi(X) = e^{-2\pi b p}$ and $\pi(Y) = e^{2\pi b x}$, with $p = \frac{1}{2\pi {\rm i}} \frac{d}{dx}$, and $\hbar = b^2$. It is observed that this is a unique irreducible representation in some sense. The tensor product representation $\mathscr{H} \otimes \mathscr{H}$ is studied, on which the quantum plane Hopf algebra $\mathcal{B}_q$ acts via the coproduct. This decomposes into a direct integral of the irreducibles $\mathscr{H}$, and this decomposition is realized as a map
\begin{align}
\label{eq:intro_F}
{\bf F} : \mathscr{H} \otimes \mathscr{H} \to M \otimes \mathscr{H},
\end{align}
where $M = L^2(\mathbb{R})$ is the multiplicity module, which is a trivial module of $\mathcal{B}_q$. For the triple tensor product, there are two decompositions $(\mathscr{H} \otimes \mathscr{H})\otimes \mathscr{H}$ and $\mathscr{H} \otimes (\mathscr{H} \otimes \mathscr{H})$, and the identity map between these spaces is encoded as a map
$$
{\bf T} : M \otimes M \to M \otimes M
$$
between the multiplicity modules. For the quadruple tensor product, there are five decompositions, and the identity maps between them can be expressed using ${\bf T}$ among the multiplicity modules; the (co-)associativity automatically yields the pentagon identity
\begin{align}
\label{eq:intro_pentagon}
{\bf T}_{23} {\bf T}_{12} = {\bf T}_{12} {\bf T}_{13} {\bf T}_{23}.
\end{align}
Meanwhile, the multiplicity module $M$ appearing in the decomposition in eq.\eqref{eq:intro_F} can be interpreted as the space of intertwiners
\begin{align}
\label{eq:intro_M}
M = {\rm Hom}_{\mathcal{B}_q}(\mathscr{H} , \mathscr{H}\otimes \mathscr{H}).
\end{align}
The left dual $\mathscr{H}'$ and the right dual ${}' \mathscr{H}$ representations of $\mathscr{H}$ are studied, by using the antipode of $\mathcal{B}_q$, as is in the standard theory of representations of Hopf algebras. The ${\rm Hom}$ spaces are studied using these dual representations, and the invariant subspaces ${\rm Hom}_{\mathcal{B}_q}$ are studied. Using only the natural constructions in representation theory, a way of cyclically permuting the roles of three spaces $\mathscr{H}$ in eq.\eqref{eq:intro_M} is developed, which yields a map
$$
{\bf A} : M \to M
$$
of order three:
\begin{align}
\label{eq:intro_A_order_three}
{\bf A}^3 = {\rm id}.
\end{align}
Again by considering representation theoretic statements and diagram chasing, the consistency equations involving ${\bf A}$ and ${\bf T}$ are proved:
\begin{align}
\label{eq:intro_two_relations}
{\bf A}_1 {\bf T}_{12} {\bf A}_2 = {\bf A}_2 {\bf T}_{21} {\bf A}_1, \quad
{\bf T}_{12} {\bf A}_1 {\bf T}_{21} = {\bf A}_1 {\bf A}_2 {\bf P}_{(12)},
\end{align}
where the subscript indices indicate the location of the tensor factors on which the operators are being applied, and ${\bf P}_{(12)}$ is the permutation map between the first and the second tensor factors. The space of intertwiners ${\rm Hom}_{\mathcal{B}_q}(\mathscr{H}^{\otimes n}, \mathscr{H}^{\otimes m})$ can be expressed as a tensor product of several copies of $M$, where each $M$ stands for the situation as in eq.\eqref{eq:intro_M}. This can be depicted pictorially as a polygon, each side representing a copy of $\mathscr{H}$, while a triangulation of this polygon yields a decomposition into a tensor product of $M$'s, where each $M$ is associated to each triangle. To be more precise, for each triangle a distinguished corner needs to be chosen to pin down the roles of the three $\mathscr{H}$'s being involved. A triangulation with distinguished corners is called a dotted triangulation. To a change of dotted triangulation given by changing one distinguished corner is associated the operator ${\bf A}$, and to a change of dotted triangulation that flips one diagonal edge of a quadrilateral to the other diagonal is associated the operator ${\bf T}$. So these operators ${\bf A}_j$ and ${\bf T}_{jk}$ represent the groupoid of changes of dotted triangulations. 

\vs

On the other hand, in a completely different story of quantum Teichm\"uller theory developed by Kashaev \cite{Ka1}, in order to quantize a certain version of the Teichm\"uller space of a punctured surface, a dotted triangulation of the surface is required, and per each change of dotted triangulations is associated a unitary operator. The main result of \cite{FKi} states that the operators ${\bf A}_j$ and ${\bf T}_{jk}$ coming from the representation theory of the quantum plane Hopf algebra $\mathcal{B}_q$ (or more precisely, the `modular double' version of $\mathcal{B}_q$) exactly coincide with Kashaev's operators $\til{\bf A}_j$ and $\til{\bf T}_{jk}$ developed in the quantum Teichm\"uller theory of punctured surfaces. We note that ${\bf T}_{jk}$ and $\til{\bf T}_{jk}$ involve the non-compact quantum dilogarithm function $\Phi^\hbar$, while ${\bf A}_j$ and $\til{\bf A}_j$ are some analogs of the Fourier transformation, or the Weil intertwiner for metaplectic group representations. We also note that, results similar to that of \cite{FKi} were obtained by Kashaev \cite{Ka0}, Bai \cite{Bai} and the author \cite{Kim_JPAA} in the case of finite dimensional representations, when $q$ is a root of unity.

\vs

In the present paper, we introduce a rather canonical Hopf $*$-algebra $\mathcal{C}_\mathbf{q}$ which we call the \emph{quantum pseudo-K\"ahler plane}, which can be viewed either as a deformation quantization of the complex two-dimensional plane $\mathbb{C}^2 = \{(z_1,z_2) : z_1,z_2\in\mathbb{C}\}$ with the pseudo-K\"ahler metric $dz_1\otimes d\ol{z}_2 + dz_2\otimes d\ol{z}_1$, or as a `complexified' version of the quantum plane; as an algebra it is
$$
\langle ~ X^{\pm 1}, Y, (X^*)^{\pm 1}, Y^* ~ | ~ XY = \mathbf{q}^2\, YX, ~~ [X,X^*]=[Y,Y^*]=[X,Y^*]=0 ~ \rangle,
$$
for a \emph{real} quantum parameter $\mathbf{q}$ with $0<{\bf q}<1$ (hence $\mathbf{q}^*=\mathbf{q}$), having the coproduct
$$
\Delta X = X\otimes X, \qquad \Delta Y = Y\otimes X + 1\otimes Y,
$$
with the $*$-structure sending $X\leftrightarrow X^*$ and $Y \leftrightarrow Y^*$ as can be expected from the notations. Then we define and study certain natural class of nicely behaved representations of the quantum pseudo-K\"ahler plane algebra as analogs of the ones for the quantum plane, where $X$ and $Y$ are represented this time as normal operators on a Hilbert space. We study the representation theory of this Hopf algebra $\mathcal{C}_{\bf q}$ in the style of \cite{FKi} as summarized above, and obtain the following main result. We omit much detail in the introduction for simplicity, such as the modular double phenomenon.
\begin{theorem}[the main result; see Thm.\ref{thm:main} for more detail]
A certain natural representation $\mathscr{H}$ of the quantum pseudo-K\"ahler plane algebra $\mathcal{C}_{\bf q}$ exists, so that the decomposition map ${\bf F} : \mathscr{H} \otimes \mathscr{H} \to M \otimes \mathscr{H}$ of its tensor square into itself, where $M = {\rm Hom}_{\mathcal{C}_{\bf q}}(\mathscr{H},\mathscr{H}\otimes \mathscr{H})$ is the trivial multiplicity module, as well as a certain natural `permutation' map on ${\rm Hom}_{\mathcal{C}_{\bf q}}(\mathscr{H},\mathscr{H}\otimes \mathscr{H})$, yield operators ${\bf T} : M \otimes M \to M \otimes M$ and ${\bf A} : M \to M$ that satisfy the operator identities in the equations \eqref{eq:intro_pentagon}, \eqref{eq:intro_A_order_three} and \eqref{eq:intro_two_relations}, up to multiplicative constants. The operator ${\bf T}$ is unitary, and there is a certain twist of ${\bf A}$ to make it unitary. These operators ${\bf T}$ and ${\bf A}$ are unique up to unitary conjugation.
\end{theorem}
This time, unlike the case in \cite{FKi}, the operator ${\bf T}_{jk}$ involves the modular double {\it compact} quantum dilogarithm functions $\Phi^{\pm {\rm i} \hbar}$. One can interpret this result as saying that natural intertwiners between representations of the quantum pseudo-K\"ahler plane algebra yield the unitary operators for a not-yet-existing Kashaev-type 3d quantum gravity theory for the positive cosmological constant case, whose existence is being suggested by this phenomenon. In other words,  we realized the Hilbert spaces of states for some to-be-found 3d quantum gravity theory for positive cosmological constant as the space of intertwiners of the newly-introduced quantum pseudo-K\"ahler plane Hopf algebra $\mathcal{C}_{\bf q}$.

\vs

The present work immediately calls for investigation of various further problems, including 1) similar questions to the ones regarding the quantum universal Teichm\"uller space as raised in \cite{FKi}, 2) a functional analytic characterization in the style of Schm\"udgen \cite{Sch92} or Ponsot-Teschner \cite{PT} of the class of representations of the quantum pseudo-K\"ahler plane, 3) the construction of a Kashaev-type quantization of 3d gravity whose existence is hinted by our main theorem, which in particular would require the study of a new coordinate system on a moduli space of 3d spacetimes that should mimic Penner's `lambda lengths', 4) the extension of the results to bigger quantum groups, either to a suitable `complexified' version of $\mathcal{U}_q(\frak{sl}_2)$ as an analog of the work of Nidaiev and Teschner \cite{NT}, or to Borel subalgebras of higher rank quantum groups as an analog of Ip's work \cite{I14}, 5) a possible relationship to some 2d conformal field theory like \cite{Te}, 6) the search for a finite-dimensional analog in the vein of \cite{Kim_JPAA}, and 7) all constructions and same questions for the zero cosmological constant case.

\vs

\emph{Acknowledgments.} This research was supported by Basic Science Research Program through the National Research Foundation of Korea(NRF) funded by the Ministry of Education(grant number 2017R1D1A1B03030230). This work was supported by the National Research Foundation of Korea(NRF) grant funded by the Korea government(MSIT) (No. 2020R1C1C1A01011151). H.K. has been supported by a KIAS Individual Grant (MG047203) at Korea Institute for Advanced Study. H.K. acknowledges the support he received from the Associate Member Program of Korea Institute for Advanced Study (Open KIAS program) during 2017--2020. H.K. thanks Dylan Allegretti and Kenji Hashimoto for helpful discussions. H.K. thanks the editors and the anonymous referees for their effort.

\vs

\emph{Conflict of interest.} The author declares that he has no conflict of interests.

\section{The Kashaev quantization of Teichm\"uller spaces}
\label{sec:Kashaev_quantization}

We briefly recall the results of Kashaev's quantum Teichm\"uller theory \cite{Ka1} \cite{Te} \cite{K16}. 

\vs

\begin{figure}[htbp!]
$\begin{array}{llll}
\begin{pspicture}[showgrid=false,linewidth=0.5pt,unit=6mm](-1,-1.2)(2.0,1.3)
%
\psarc[arcsep=0.5pt](1.091,-3.055){2.182}{61.7}{158.0}
\psarc[arcsep=0.5pt](-7.2,0.4){6.8}{-22}{16.3}
\psarc[arcsep=0.5pt](4.8,3.9){5.7}{-163.7}{-118.3}
%
%
\rput{127.0}(1.2,1.0){\fontsize{19}{19} $\cdots$}
\rput{86.0}(-1.4,0.2){\fontsize{19}{19} $\cdots$}
\rput{20.0}(0.58,-1.65){\fontsize{17}{17} $\cdots$}
%
\rput(-0.5,-1.3){\fontsize{11}{11} $\bullet$}
%
\rput(0.1,-0.3){\fontsize{11}{11} $j$}
%
\rput[l](3.2,0){\pcline[linewidth=0.7pt, arrowsize=2pt 4]{->}(0,0)(1;0)\Aput{\,$A_{j}$}}
\end{pspicture}
%
%
%
%
& \begin{pspicture}[showgrid=false,linewidth=0.5pt,unit=6mm](-2.0,-1.2)(2.4,1.3)
%
\psarc[arcsep=0.5pt](1.091,-3.055){2.182}{61.7}{158.0}
\psarc[arcsep=0.5pt](-7.2,0.4){6.8}{-22}{16.3}
\psarc[arcsep=0.5pt](4.8,3.9){5.7}{-163.7}{-118.3}
%
\rput(1.1,-0.7){\fontsize{11}{11} $\bullet$}
%
\rput(0.1,-0.3){\fontsize{11}{11} $j$}
%
\rput{127.0}(1.2,1.0){\fontsize{19}{19} $\cdots$}
\rput{86.0}(-1.4,0.2){\fontsize{19}{19} $\cdots$}
\rput{20.0}(0.58,-1.65){\fontsize{17}{17} $\cdots$}
\end{pspicture}
%
%
%
%
& \begin{pspicture}[showgrid=false,linewidth=0.5pt,unit=6mm](-2.0,-1.2)(2.0,1.3)
%
\psarc[arcsep=0.5pt](0.343,-4.457){3.771}{61.7}{126.4}
\psarc[arcsep=0.5pt](-3.36,0.48){2.4}{-53.6}{37.4}
\psarc[arcsep=0.5pt](-0.185,2.862){1.569}{-142.6}{-29}
\psarc[arcsep=0.5pt](3.220,0.937){2.341}{151}{-118.3}
%
\psarc[arcsep=0.5pt](7.2,8.4){10.8}{-142.6}{-118.3}
%
\rput(-0.95,1.00){\fontsize{10}{10} $\bullet$}
\rput(0.72,1.42){\fontsize{10}{10} $\bullet$}
%
\rput(-0.5,-0.1){\fontsize{10}{10} $j$}
\rput(0.3,0.8){\fontsize{10}{10} $k$}
%
\rput{104.0}(1.7,0.5){\fontsize{15}{15} $\cdots$}
\rput{5.0}(-0.10,1.94){\fontsize{15}{15} $\cdots$}
\rput{85.0}(-1.7,0.3){\fontsize{15}{15} $\cdots$}
\rput{4.0}(0.08,-1.50){\fontsize{17}{17} $\cdots$}
%
\rput[l](3.2,0){\pcline[linewidth=0.7pt, arrowsize=2pt 4]{->}(0,0)(1;0)\Aput{$T_{jk}$}}
\end{pspicture}
%
%
%
%
& \begin{pspicture}[showgrid=false,linewidth=0.5pt,unit=6mm](-2.0,-1.2)(2.4,1.3)
%
\psarc[arcsep=0.5pt](0.343,-4.457){3.771}{61.7}{126.4}
\psarc[arcsep=0.5pt](-3.36,0.48){2.4}{-53.6}{37.4}
\psarc[arcsep=0.5pt](-0.185,2.862){1.569}{-142.6}{-29}
\psarc[arcsep=0.5pt](3.220,0.937){2.341}{151}{-118.3}
%
%
\psarc[arcsep=0.5pt](-8.8,7.733){11.467}{-53.0}{-29.7}
%
\rput(-1.00,1.30){\fontsize{10}{10} $\bullet$}
\rput(0.67,1.10){\fontsize{10}{10} $\bullet$}
%
\rput(-0.50,0.8){\fontsize{10}{10} $j$}
\rput(0.3,-0.1){\fontsize{10}{10} $k$}
%
\rput{104.0}(1.7,0.5){\fontsize{15}{15} $\cdots$}
\rput{5.0}(-0.10,1.94){\fontsize{15}{15} $\cdots$}
\rput{85.0}(-1.7,0.3){\fontsize{15}{15} $\cdots$}
\rput{4.0}(0.08,-1.50){\fontsize{17}{17} $\cdots$}
\end{pspicture}
\end{array} $
\caption{The actions of $A_j$ and $T_{jk}$ on dotted ideal triangulations}
\label{fig:action_on_dotted_ideal_triangulations}
\end{figure}
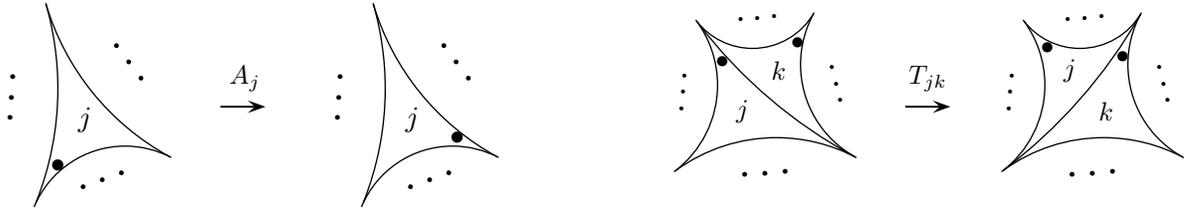

First, let $S$ be an oriented real surface with punctures. An \emph{ideal triangulation} of $S$ is a collection of isotopy classes of unoriented paths running between punctures, called \emph{edges}, so that the complementary region in $S$ is the disjoint union of \emph{ideal triangles}, i.e. regions bounded by three edges. A \emph{dotted ideal triangulation} of $S$, used in the Kashaev quantization is an ideal triangulation of $S$ together with the choice of a distinguished corner for each ideal triangle indicated by a dot $\bullet$ in pictures, and the choice of labeling of the ideal triangles by some fixed index set $I$. There are three types of \emph{elementary} transformations of dotted ideal triangles, namely 1) the dot-change $A_j$ for a triangle $j\in I$, moving the dot $\bullet$ in $j$ counterclockwise, 2) the flip $T_{jk}$ for two distinct adjacent triangles $j,k\in I$ whose dots are configured like in Fig.\ref{fig:action_on_dotted_ideal_triangulations} \footnote{We note that the basic source code for Fig.\ref{fig:action_on_dotted_ideal_triangulations} is taken from \cite{K16}.}, replacing the common edge of $j$ and $k$ by the new one obtained by rotating it $90^\circ$ clockwise, and 3) the label change $P_\sigma$ for ideal triangles associated a permutation $\sigma$ of $I$, relabeling each triangle $j\in I$ by $\sigma(j)$. 


\begin{proposition}[{\cite[\S3]{Ka2} \cite[Appendix.C]{Te} \cite[Thm.2.23]{K16} \cite[\S4.1]{K16b}}] 
Any two dotted ideal triangles are connected by a finite sequence of these elementary transformations, and all algebraic relations satisfied by the elementary transformations are generated by the following list:
$$
A_j^3 = \mathrm{id}, \quad T_{k\ell} T_{jk} = T_{jk} T_{j\ell} T_{k\ell}, \quad A_j T_{jk} A_k = A_k T_{kj} A_j, \quad T_{jk} A_j T_{kj} = A_j A_k P_{(jk)},
$$
for mutually distinct $j,k,\ell\in I$, together with the trivial relations regarding the index permutations and commutativity of the factors with mutually-non-intersecting sets of indices.
\end{proposition}

\vs

Given any dotted ideal triangulation, Kashaev constructs a coordinate system consisting of the \emph{ratio coordinate functions} associated to ideal triangles, parametrizing a space related to the Teichm\"uller space of $S$. One advantage of using Kashaev's ratio coordinates is that they provide an analog of a Darboux-like `diagonal' system for the Poisson/symplectic structure on the Teichm\"uller space of $S$. An elementary transformation of dotted ideal triangulations induces a change of ratio coordinate functions, which are related to the previous ones by certain simple rational formulas.

\vs

To obtain a quantization we must replace the ratio coordinate functions by self-adjoint operators on a Hilbert space that depend real analytically on a real quantization parameter $\hbar$, so that the operator commutators recover the classical Poisson structure in the first order of $\hbar$ as $\hbar\to 0$, i.e. up to $o(\hbar)$. There is a non-unique but somewhat canonical way of doing this, for each chosen dotted ideal triangulation. The main nontrivial question is to remove the dependency of such a quantization on the choice of a dotted ideal triangulation. Namely, for each transformation of dotted ideal triangulations we must find a unitary map between the corresponding Hilbert spaces such that 1) it intertwines the two sets of quantum ratio-coordinate operators that are supposed to be related by some quantum version of the classical coordinate change formulas, and that 2) such an assignment of unitary maps to transformations of dotted ideal triangulations must be consistent, i.e. preserves the compositions up to constants. Kashaev's answer \cite{Ka1} can be described in terms of the following unitary intertwining operators for the corresponding elementary transformations $A_j$, $T_{jk}$, $P_\sigma$ of dotted ideal triangulations
\begin{align}
\label{eq:Kashaev_operators}
\til{\mathbf{A}}_j = e^{-\pi {\rm i}/3} \, e^{3\pi {\rm i}x_j^2} \, e^{\pi{\rm i}\,(p_j+x_j)^2}, \quad \til{\mathbf{T}}_{jk} = e^{2\pi{\rm i}\,p_j x_k} \, e_b(x_j+p_k-x_k)^{-1}, \quad \til{\mathbf{P}}_{\sigma},
\end{align}
on a single Hilbert space
\begin{align}
\label{eq:Kashaev_space}
\mathcal{H} = L^2(\mathbb{R}^I, \, \textstyle \bigwedge_{j\in I} dx_j),
\end{align}
where $\til{\mathbf{P}}_\sigma$ is just the permutation of variables
$$
(\til{\mathbf{P}}_\sigma f)(\{x_j\}_{j\in I}) = f(\{x_{\sigma(j)}\}_{j\in I}), \qquad \forall f\in \mathcal{H},
$$
and $x_j,p_j$ are self-adjoint operators which can be described as essentially self-adjoint operators given on a dense subspace by
$$
\mbox{$x_j=$multiplication by $x_j$,\quad and $p_j = \frac{1}{2\pi {\rm i}} \, \frac{\partial}{\partial x_j}$},
$$
while the function $e_b$ is one version of the non-compact quantum dilogarithm function introduced by Faddeev and Kashaev \cite{FK} related to $\Phi^\hbar$ in eq.\eqref{eq:Phi_hbar_first_time} by
$$
e_b(z) = \Phi^\hbar(2\pi b z)^{-1} \quad\mbox{with}\quad \hbar = b^2, ~ b>0;
$$
the expressions in eq.\eqref{eq:Kashaev_operators} are made sense by using the functional calculus of self-adjoint operators. 

\vs

In particular, these operators satisfy
\begin{align}
\label{eq:Kashaev_operators_relations}
\til{\mathbf{A}}_j^3 = \mathrm{id}, \quad \til{\mathbf{T}}_{k\ell} \, \til{\mathbf{T}}_{jk} = \til{\mathbf{T}}_{jk} \, \til{\mathbf{T}}_{j\ell} \, \til{\mathbf{T}}_{k\ell}, \quad \til{\mathbf{A}}_j \til{\mathbf{T}}_{jk} \til{\mathbf{A}}_k = \til{\mathbf{A}}_k \til{\mathbf{T}}_{kj} \til{\mathbf{A}}_j, \quad \til{\mathbf{T}}_{jk} \til{\mathbf{A}}_j \til{\mathbf{T}}_{kj} = \zeta \til{\mathbf{A}}_j \til{\mathbf{A}}_k \til{\mathbf{P}}_{(jk)},
\end{align}
and also the trivial relations as well, where
$$
\zeta = e^{-2\pi {\rm i} (b+b^{-1})^2/24}.
$$
But of course, when seen for the first time, the formulas in eq.\eqref{eq:Kashaev_operators} may not come to one's mind swiftly as being some natural solutions to the problem. In other words, this quantization problem might seem to be a very difficult one to solve, judging by the complexity of the solution. So one might seek for some other natural and easier explanation that leads to this solution.

\section{The Frenkel-Kim construction for intertwiners of the quantum plane}
\label{sec:algebraic_story}

The pentagon identity $\til{\mathbf{T}}_{k\ell} \, \til{\mathbf{T}}_{jk} = \til{\mathbf{T}}_{jk} \, \til{\mathbf{T}}_{j\ell} \, \til{\mathbf{T}}_{k\ell}$ arises naturally in other areas of mathematics as well, one example being the consistency axiom for the associativity morphisms $(A\otimes B) \otimes C \to A\otimes(B\otimes C)$ in a tensor category. So one may ask if there exists a natural tensor category of certain Hilbert spaces whose associativity morphisms coincide with Kashaev's operators $\til{\mathbf{T}}_{jk}$. Frenkel and the author \cite{FKi} construct such a category, as the category of some nice class of representations of a rather basic Hopf algebra, namely the quantum plane $\mathcal{B}_q$ as introduced in \S\ref{sec:introduction}. Moreover, it is shown in \cite{FKi} that this is in fact a rigid tensor category in a sense, and from certain maps between dual objects they also recover Kashaev's operator $\til{\mathbf{A}}_j$, while all relations in eq.\eqref{eq:Kashaev_operators_relations} are proven essentially just by the representation theory of the quantum plane. In the present section we briefly review this work \cite{FKi}, of which the entire framework is mimicked in the present paper.

\vs

An \emph{integrable representation} of the quantum plane $\mathcal{B}_q$ is a pair $(\mathscr{H},\pi)$ consisting of a separable complex Hilbert space $\mathscr{H}$ and positive-definite (densely-defined) self-adjoint operators $\pi(X)$ and $\pi(Y)$ on $\mathscr{H}$ satisfying
$$
\pi(X)^{{\rm i} \alpha} \, \pi(Y)^{{\rm i} \beta} = e^{-2\alpha\beta \pi{\rm i} b^2} \, \pi(Y)^{{\rm i}\beta} \, \pi(X)^{{\rm i}\alpha}, \quad \forall \alpha,\beta\in\mathbb{R},
$$
where $q = e^{\pi {\rm i} b^2}$, and the unitary operators $\pi(X)^{{\rm i} \alpha}$ and $\pi(Y)^{{\rm i} \beta}$ are constructed by functional calculus. Then, on a dense subspace we have the sought-for algebraic relation $\pi(X) \pi(Y) = e^{2\pi {\rm i} b^2} \pi(Y) \pi(X)$. We say $(\mathscr{H},\pi)$ is \emph{irreducible} if $0$ and $\mathscr{H}$ are the only closed $\mathbb{C}$-vector subspaces that are invariant under $\pi(X)^{{\rm i}\alpha}$ and $\pi(Y)^{{\rm i}\beta}$ for all $\alpha,\beta\in\mathbb{R}$; this is a different but better definition of irreducibility than the one used in \cite{FKi}. Schm\"udgen \cite{Sch92} showed that there is a unique irreducible integrable representation of $\mathcal{B}_q$ up to unitary equivalence, given by
$$
\mathscr{H}=L^2(\mathbb{R},dx), \qquad \pi(X) = e^{-2\pi bp}, \qquad \pi(Y) = e^{2\pi bx},
$$
where $p = \frac{1}{2\pi{\rm i}} \frac{d}{dx}$. In particular, $\pi(X)$ acts as a complex shift operator $(\pi(X) f)(x) = f(x+{\rm i}b)$ on a dense subspace of $\mathscr{H}$. We omit $\pi$ from $(\mathscr{H},\pi)$ and just write $\mathscr{H}$, whenever it is clear.

\vs

Then the tensor squared representation $\mathscr{H} \otimes \mathscr{H}$ is considered, where each element $u$ of $\mathcal{B}_q$ acts via the coproduct $(\pi\otimes\pi)(\Delta u)$. Throughout the paper, whenever we take the tensor product of two separable Hilbert spaces $V_1$ and $V_2$ (with inner products $\langle \cdot,\cdot \rangle_1$ and $\langle \cdot,\cdot \rangle_2$), we mean the Hilbert space tensor product (see {\cite[\S II.4]{RS}}), which is defined as the completion of the algebraic tensor product $V_1 \otimes V_2$ equipped with the inner product given by $\langle \varphi_1 \otimes \varphi_2, \psi_1 \otimes \psi_2 \rangle = \langle \varphi_1, \psi_1 \rangle_{1} \cdot \langle \varphi_2,\psi_2\rangle_{2}$; for convenience, we still denote the resulting Hilbert space as $V_1 \otimes V_2$. This in particular applies to $\mathscr{H} \otimes \mathscr{H}$ here. One can check that this representation ($\mathscr{H} \otimes \mathscr{H}, \pi\otimes \pi)$ is integrable, and hence try to decompose into the direct sum of irreducibles, i.e. of $\mathscr{H}$. It turns out that it can only be decomposed as a direct integral $\mathscr{H} \otimes \mathscr{H} \cong \int^\oplus_\mathbb{R} \mathscr{H}$, in a certain sense; to be precise, a unitary intertwiner
\begin{align}
\label{eq:FKi_bf_F}
\mathbf{F} : \mathscr{H}\otimes \mathscr{H} \to M \otimes \mathscr{H}
\end{align}
is constructed, where the {\it multiplicity module} $M$ is realized as $L^2(\mathbb{R})$ as a Hilbert space and as the trivial module of the Hopf algebra $\mathcal{B}_q$. The intertwining equations that $\mathbf{F}$ must solve boil down to certain `difference equations' of a function $F$ expressing $\mathbf{F}$ by functional calculus, in the sense that some shift in the argument of $F$ results in gaining some factor, like in the case of the classical Gamma function $\Gamma(z+1)= z \Gamma(z)$; one could have expected such a situation, for $\pi(X)$ is a shift operator. The difference equation for $F$ was recognized as the one satisfied by the non-compact quantum dilogarithm function $e_b$ of Faddeev-Kashaev, leading to a solution for $F$. 

\vs

Next, the tensor cubed reprepresentation $\mathscr{H}\otimes \mathscr{H}\otimes \mathscr{H}$ is considered. Again this decomposes into a direct integral of $\mathscr{H}$'s, but in two\footnote{the third Catalan number} ways: one by $(\mathscr{H}\otimes \mathscr{H}) \otimes \mathscr{H} \to M\otimes \mathscr{H}\otimes\mathscr{H} \to M\otimes M \otimes \mathscr{H}$ and the other by $\mathscr{H}\otimes (\mathscr{H}\otimes \mathscr{H}) \to \mathscr{H}\otimes M \otimes \mathscr{H} \to M\otimes M \otimes \mathscr{H}$, where each arrow uses $\mathbf{F}$. Composing these four arrows one obtains a unitary intertwining map $M\otimes M \otimes \mathscr{H} \to M\otimes M \otimes \mathscr{H}$ which encodes the `$\mathscr{H}$-level' identity intertwining map $(\mathscr{H}\otimes\mathscr{H})\otimes\mathscr{H} \to \mathscr{H}\otimes(\mathscr{H}\otimes\mathscr{H})$, i.e. encodes the associativity of tensor products induced by the coassociativity of the coproduct, or can be viewed simply as encoding the change of parentheses. Either by using the irreducibility of $\mathscr{H}$ or by an explicit computation, one finds out that this map $M\otimes M\otimes \mathscr{H} \to M\otimes M\otimes \mathscr{H}$ in fact acts nontrivially only on the first two tensor factors, by a unitary operator
$$
\mathbf{T} : M\otimes M \to M\otimes M,
$$
which is expressed by functional calculus of some self-adjoint operator applied to the non-compact quantum dilogarithm function $e_b$, just like $\mathbf{F}$ is. Then, consider the tensor quadrupled representation $\mathscr{H} \otimes \mathscr{H} \otimes \mathscr{H} \otimes \mathscr{H}$. There are five\footnote{the fourth Catalan number} different ways of parenthesizing, each of which tells us how to decompose this representation into $M\otimes M \otimes M \otimes \mathscr{H}$. The identity maps in the $\mathscr{H}$-level are then carried over to five $\mathbf{T}$'s among $M\otimes M \otimes M$, in particular leading to the pentagon identity
$$
\mathbf{T}_{23} \mathbf{T}_{12} = \mathbf{T}_{12} \mathbf{T}_{13} \mathbf{T}_{23},
$$
where $\mathbf{T}_{ab}$ means the operator $\mathbf{T}:M\otimes M\to M\otimes M$ acting on the $a$-th and the $b$-th tensor factors in $M\otimes M\otimes M$. One can view this pentagon identity as being the consistency equation for the `associativity map' $\mathbf{T}$. Consistency equations coming from higher tensor powers are all consequences of this pentagon identity.

\vs

Meanwhile, the situation in eq.\eqref{eq:FKi_bf_F} is naturally interpreted as
\begin{align}
\label{eq:FKi_M}
M \cong \mathrm{Hom}_{\mathcal{B}_q}(\mathscr{H},\mathscr{H}\otimes \mathscr{H}),
\end{align}the space of intertwiners between $\mathscr{H}$ and $\mathscr{H}\otimes \mathscr{H}$. Then the left and the right dual representations $\mathscr{H}'$ and ${}'\mathscr{H}$ are studied, which are defined on the Hilbert space $L^2(\mathbb{R})$ as the transpose of the action on $\mathscr{H} = L^2(\mathbb{R})$ by the antipode or by the inverse of the antipode, as usually studied in the ordinary finite-dimensional representation theory of Hopf algebras. The unique non-unitary intertwining maps $\mathscr{H} \to \mathscr{H}'$ and $\mathscr{H} \to {}'\mathscr{H}$ are constructed. Then, by composing natural intertwiners and vector space isomorphisms among $\mathrm{Hom}_\mathbb{C}(\mathscr{H},\mathscr{H}\otimes\mathscr{H})$, ${}'\mathscr{H} \otimes (\mathscr{H} \otimes \mathscr{H})$, and $(\mathscr{H}\otimes\mathscr{H})\otimes \mathscr{H}'$ (for a careful treatment of the symbol ${\rm Hom}$, see \S\ref{sec:dual_representation}), with the help of the above maps $\mathscr{H}\to \mathscr{H}'$ and $\mathscr{H} \to {}' \mathscr{H}$, one can construct a map
$$
\mathbf{A} : M \to M
$$
which cyclically permutes the three $\mathscr{H}$'s in the right hand side of eq.\eqref{eq:FKi_M}. One can show
$$
\mathbf{A}^3 = \mathrm{id},
$$
and moreover also
$$
\mathbf{A}_1 \mathbf{T}_{12} \mathbf{A}_2 = \mathbf{A}_2 \mathbf{T}_{21} \mathbf{A}_1, \qquad
\mathbf{T}_{12} \mathbf{A}_1 \mathbf{T}_{21} = \mathbf{A}_1 \mathbf{A}_2 \, \mathbf{P}_{(12)}, 
$$
with $\mathbf{P}_{(12)}:M\otimes M \to M\otimes M$ being the permutation of the two factors, where each of these two equations comes from the consistency of certain natural maps between spaces of intertwiners, i.e. spaces looking like $\mathrm{Hom}_{\mathcal{B}_q}(V,W)$ with $V$ and $W$ being tensor powers of $\mathscr{H}$. 

\vs

Such an intertwiner space $\mathrm{Hom}_{\mathcal{B}_q}(V,W)$ is described pictorially as a polygon, with a dotted triangulation determined by the choices of parentheses for $V$ and $W$; each dotted triangle corresponds to one $M$, so that one can identify $\mathrm{Hom}_{\mathcal{B}_q}(V,W)$ as the tensor product of as many $M$'s as there are triangles. Changes of parentheses and permutations of $\mathscr{H}$-factors then correspond to changes of dotted triangulations of the polygon, and also to the maps on the tensor product of $M$'s given by compositions of $\mathbf{T}$'s, $\mathbf{A}$'s, and tensor-factor-permutation $\mathbf{P}$'s. Hence, like in Kashaev's quantum Teichm\"uller theory, to elementary changes of dotted triangulations are assigned maps between Hilbert spaces, which satisfy the algebraic consistency relations among them. Moreover, a unitary map $\mathbf{U} : M \to L^2(\mathbb{R})$ such that 
$$
\mathbf{U}^{\otimes n} \, \mathbf{A}_j^{(0)} \, (\mathbf{U}^{\otimes n})^{-1} = \til{\mathbf{A}}_j, \qquad
\mathbf{U}^{\otimes n} \, \mathbf{T}_{jk} \, (\mathbf{U}^{\otimes n})^{-1} = \til{\mathbf{T}}_{jk}, \qquad
\mathbf{U}^{\otimes n} \, \mathbf{P}_\sigma \, (\mathbf{U}^{\otimes n})^{-1} = \til{\mathbf{P}}_\sigma
$$
is found, where $n$ is the number of triangles or equivalently that of $M$-factors, and the unitary operator
$$
\mathbf{A}^{(0)} := e^{\pi{\rm i}(b+b^{-1})^2/12} \, \mathbf{A} \, e^{-\pi (b+b^{-1})p}
$$
is a slight twist of the representation-theoretic $\mathbf{A}$. So one may summarize by saying that the authors of \cite{FKi} realized the Hilbert spaces $\mathscr{H}$ \eqref{eq:Kashaev_space} of states of Kashaev's quantum Teichm\"uller theory for `polygons' in terms of the spaces of intertwiners of the quantum plane Hopf algebra $\mathcal{B}_q$.

\vs

The goal of the present paper is to do what was done in \cite{FKi} as just described, for a different Hopf algebra and for a different quantum geometry, both of which involve a new special function different from the usual non-compact quantum dilogarithm.

\vs

To finish the section, we mention a few important aspects which are studied in \cite{FKi} but omitted above for the sake of the simplicity of discussions. One is about the non-uniqueness of $\mathbf{F}$; any other unitary intertwiner $\mathbf{F}' : \mathscr{H}\otimes \mathscr{H} \to M \otimes \mathscr{H}$ is related by $\mathbf{F}' = (U\otimes 1) \circ \mathbf{F}$ for a unitary map $U:M\to M$ and vice versa, resulting in the corresponding new operators $\mathbf{A}', \mathbf{T}'$ related to the previous ones by $\mathbf{A}' = U \, \mathbf{A} \, U^{-1}$ and $\mathbf{T}' = U^{\otimes 2} \, \mathbf{T} \, (U^{\otimes 2})^{-1}$; so we have a firm control over this non-uniqueness of $\mathbf{F}$. The other is about the `modular double' of the quantum plane algebra. Due to the mysterious built-in invariance property $b\leftrightarrow b^{-1}$ of the non-compact quantum dilogarithm function $e_b$ used in the expression of $\mathbf{F}$, the map $\mathbf{F}$ not only intertwines the representations of the quantum plane $\mathcal{B}_q$ for $q = e^{\pi{\rm i} b^2}$ but automatically also those of another copy of the quantum plane $\mathcal{B}_{\til{q}}$ for $\til{q} = e^{\pi {\rm i} /b^2}$, regarded as the `$1/b^2$-th' power of $\mathcal{B}_q$ in terms of representations, in a certain sense. These phenomena also appear in the following sections, so we do not explain them further at the moment.

\section{The quantum pseudo-K\"ahler plane algebra $\mathcal{C}_\mathbf{q}$}
\label{sec:quantum_pK_plane_algebra}

The main original contents of the present paper begin here with the following new Hopf $*$-algebra.

\begin{definition}[the quantum pseudo-K\"ahler plane]
Let
$$
\hbar \in \mathbb{R}, \qquad \hbar>0,
$$
and
$$
\mathbf{q} = e^{-\pi\hbar}.
$$
Define the \ul{\em quantum pseudo-K\"ahler plane algebra} $\mathcal{C}_{\mathbf{q}}$ as the $*$-algebra over $\mathbb{Z}[\mathbf{q},\mathbf{q}^{-1}]$ generated by
$$
\mbox{$Z_1^{\pm 1}$, $Z_2$, $(Z_1^*)^{\pm 1}$, $Z_2^*$,}
$$
with the $*$-structure ${}^* : \mathcal{C}_\mathbf{q} \to \mathcal{C}_\mathbf{q}$ being the $\mathbb{Z}$-algebra anti-isomorphism given by
\begin{align}
{}^* ~ : ~ Z_1 \mapsto Z_1^*, \quad Z_2 \mapsto Z_2^*, \quad Z_1^* \mapsto Z_1, \quad Z_2^* \mapsto Z_2, \quad \mathbf{q} \mapsto \mathbf{q},
\end{align}
mod out by the relations
\begin{align}
\label{eq:C_q_relation}
Z_1 \, Z_2 = \mathbf{q}^2 \, Z_2 \, Z_1, \quad
\mbox{and each of $Z_1,Z_2$ commutes with each of $Z_1^*, Z_2^*$}.
\end{align}
\end{definition}
In particular, we have $Z_1^* \, Z_2^* = \mathbf{q}^{-2} \, Z_2^* \, Z_1^*$ as a consequence.

\vs

Maybe, it might be better to define the algebra with slightly different relations: namely, to replace the set of relations as in eq.\eqref{eq:C_q_relation} by the following modified version:
$$
Z_1 \, Z_2^* = \mathbf{q} \, Z_2^* \, Z_1, \quad \mbox{and each of $Z_1$, $Z_2^*$ commutes with each of $Z_1^*$, $Z_2$}.
$$
Then, the algebra could be thought of as a deformation quantization of the \emph{pseudo-K\"ahler plane}, i.e. the complex surface $\mathbb{C}^2 = \{(z_1,z_2) : z_1,z_2\in\mathbb{C}\}$ equipped with the pseudo-K\"ahler metric $h = dz_1\otimes d\ol{z}_2 + dz_2 \otimes d\ol{z}_1$. The associated Hermitian form is $\frac{{\rm i}}{2} (dz_1 \wedge d\ol{z}_2 + dz_2 \wedge d\ol{z}_1)$, so that normal operators $\hat{z}_1$ and $\hat{z}_2$ satisfying\footnote{The normality gives $[\hat{z}_1,\hat{z}_1^*]=0=[\hat{z}_2,\hat{z}_2^*]$.} the Heisenberg relations $[\hat{z}_1,\hat{z}_2]=0=[\hat{z}_1^*,\hat{z}_2^*]$ and $[\hat{z}_1,\hat{z}_2^*] = -\pi\hbar = 2\pi{\rm i}\,\hbar \cdot \frac{{\rm i}}{2} =  [\hat{z}_2, \hat{z}_1^*]$ can indeed be viewed as a deformation quantization of the pseudo-K\"ahler plane. Then, one can figure out that their exponentials $\hat{Z}_1 = e^{\hat{z}_1}$, $\hat{Z}_1^* = e^{\hat{z}_1^*}$, $\hat{Z}_2 = e^{\hat{z}_2}$, $\hat{Z}_2^* = e^{\hat{z}_2^*}$ would satisfy the relations $\hat{Z}_1 \, \hat{Z}_2^* = \mathbf{q} \, \hat{Z}_2^* \, \hat{Z}_1$, etc, as just mentioned, e.g. by formally using the Baker-Campbell-Hausdorff (BCH) formula
\begin{align}
\label{eq:BCH}
e^A e^B = e^{A+B+\frac{1}{2}[A,B]+\frac{1}{12}([A,[A,B]]+[B,[B,A]])+\cdots}.
\end{align}
Indeed, for example, $e^{\hat{z}_1} \, e^{\hat{z}_2^*} = e^{-\pi\hbar/2} \, e^{\hat{z}_1 + \hat{z}_2^*}$ and $e^{\hat{z}_2^*} \, e^{\hat{z}_1} = e^{\pi\hbar/2} \, e^{\hat{z}_2^* + \hat{z}_1}$, so $e^{\hat{z}_1} \, e^{\hat{z}_2^*} = e^{-\pi\hbar} \, e^{\hat{z}_2^*} \, e^{\hat{z}_1}$, etc. This justisfies the name `quantum pseudo-K\"ahler plane' for $\mathcal{C}_\mathbf{q}$. However, for the notational simplicity, let us stick to eq.\eqref{eq:C_q_relation} as the defining set of relations for $\mathcal{C}_\mathbf{q}$. We now move on to the Hopf algebra structure on the algebra $\mathcal{C}_\mathbf{q}$.

\begin{lemma}
The coproduct $\Delta : \mathcal{C}_\mathbf{q} \to \mathcal{C}_\mathbf{q} \otimes \mathcal{C}_\mathbf{q}$, counit $\epsilon : \mathcal{C}_{\mathbf{q}}\to \mathbb{Z}[\mathbf{q},\mathbf{q}^{-1}]$, and the antipode $S : \mathcal{C}_{\mathbf{q}} \to \mathcal{C}_{\mathbf{q}}$ given on the generators by
\begin{align*}
\begin{array}{ll}
\Delta(Z_1) = Z_1 \otimes Z_1, & \quad \Delta(Z_2) = Z_2 \otimes Z_1 + 1 \otimes Z_2, \\
\Delta(Z_1^*) = Z_1^* \otimes Z_1^*, & \quad \Delta(Z_2^*) = Z_2^* \otimes Z_1^* + 1 \otimes Z_2^*, \\
\epsilon(Z_1) = 1 = \epsilon(Z_1^*), & \quad \epsilon(Z_2) = 0 = \epsilon(Z_2^*),  \\
S(Z_1) = Z_1^{-1}, &  \quad S(Z_2) = -  Z_2 \, Z_1^{-1}, \\
S(Z_1^*) = (Z_1^*)^{-1}, &  \quad S(Z_2^*) = -  Z_2^* \, (Z_1^*)^{-1},
\end{array}
\end{align*}
define a Hopf $*$-algebra structure on $\mathcal{C}_\mathbf{q}$. \qed
\end{lemma}
Here, a Hopf $*$-algebra means a Hopf algebra whose coproduct and counit are $*$-maps; see e.g. \cite{CP} for basic definitions. We omit a proof as it follows almost immediately from a well-known corresponding statement for the quantum plane Hopf algebra $\mathcal{B}_q$ which appeared in \S\ref{sec:introduction}, whose counit and antipode are given by $\epsilon(X)=1$, $\epsilon(Y)=0$ and $S(X) = X^{-1}$, $S(Y) = -YX^{-1}$ respectively, which in particular is a Hopf $*$-subalgebra of the split real quantum group $\mathcal{U}_q(\frak{sl}(2,\mathbb{R}))$; see \cite{FKi}.

\begin{remark}[about the name]
The Hopf algebra $\mathcal{C}_\mathbf{q}$ can be viewed as a certain `complexified' version of the quantum plane Hopf algebra $\mathcal{B}_q$. The letter $\mathcal{C}$ stands for the adjective `c'omplexified. We avoided using the term `complex quantum plane' or `quantum complex plane' for $\mathcal{C}_\mathbf{q}$, as some others used that term for different algebras; see e.g. \cite{So} \cite{CW}.
\end{remark}

We will also consider the `(Langlands) modular double' $\mathcal{C}_{\mathbf{q},\mathbf{q}^\vee}$ of this quantum pseudo-K\"ahler plane Hopf $*$-algebra $\mathcal{C}_\mathbf{q}$, which we define to be
\begin{align}
\label{eq:C_q_q_vee}
\mathcal{C}_{\mathbf{q},\mathbf{q}^\vee} := \mathcal{C}_{\mathbf{q}} \otimes \mathcal{C}_{1/\mathbf{q}^\vee}
\end{align}
as an algebra, where the `(Langlands) modular dual' part $\mathcal{C}_{1/\mathbf{q}^\vee}$ is another copy of the quantum pseudo-K\"ahler plane algebra with the parameter $1/{\bf q}^\vee$, where
$$
\mathbf{q}^\vee := e^{-\pi/\hbar}.
$$
Hence, $\mathcal{C}_{\mathbf{q},\mathbf{q}^\vee}$ can be viewed as the algebra over $\mathbb{Z}[\mathbf{q}^{\pm1}, (\mathbf{q}^{\vee})^{\pm 1}]$ generated by
$$
Z_1^{\pm 1}, ~ Z_2, ~ (Z_1^\vee)^{\pm 1}, ~ Z_2^\vee, \quad  (Z_1^*)^{\pm 1}, ~ Z_2^*, ~ ((Z_1^\vee)^*)^{\pm 1}, ~ (Z_2^\vee)^*,
$$
mod out by the relations
\begin{align}
\label{eq:C_q_q_vee_relation}
\left\{
\begin{array}{l}
 Z_1 \, Z_2 = \mathbf{q}^2\, Z_2\, Z_1, \quad Z_1^\vee Z_2^\vee = (\mathbf{q}^\vee)^{-2} \, Z_2^\vee \, Z_1^\vee, \\ 
\mbox{each of $Z_1$, $Z_2$ commutes with each of $Z_1^*$, $Z_2^*$,} \\
\mbox{each of $Z_1^\vee$, $Z_2^\vee$ commutes with each of $(Z_1^\vee)^*$, $(Z_2^\vee)^*$,} \\
 \mbox{each of $Z_1,Z_2,Z_1^*,Z_2^*$ commutes with each of $Z_1^\vee,Z_2^\vee,(Z_1^\vee)^*,(Z_2^\vee)^*$;}
  \end{array} \right.
\end{align}
the $*$-structure and the Hopf algebra structure come from those of $\mathcal{C}_{\mathbf{q}}$ and $\mathcal{C}_{1/\mathbf{q}^\vee}$. In addition, we shall consider the following `transcendental' relations
\begin{align}
\label{eq:transcendental_relation}
Z_1^{1/({\rm i}\hbar)} = Z_1^\vee, \qquad Z_2^{1/({\rm i}\hbar)} = Z_2^\vee,
\end{align}
where the $\frac{1}{{\rm i}\hbar}$-th power on the left hand side should be made a suitable sense when dealing with representations, i.e. when the generators are represented as operators on a Hilbert space.

\begin{remark}
Maybe it is better to define $\mathbf{q}^\vee$ as $e^{\pi/\hbar}$ instead, in order to indicate the `modularity' phenomenon $\hbar \leftrightarrow -1/\hbar$ with respect to ${\bf q} = e^{-\pi \hbar}$  more clearly. The advantage of defining ${\bf q}^\vee$ as $e^{-\pi/\hbar}$ for $\hbar >0$ is that then we have $|{\bf q}^\vee|<1$, being aligned with the corresponding condition $|{\bf q}|<1$. 
\end{remark}

\section{The integrable representation $(\mathscr{H},\pi)$ of $\mathcal{C}_\mathbf{q}$}
\label{sec:integrable_representation_of_C_q}

Consider the separable complex Hilbert space
\begin{align}
\label{eq:H}
\mathscr{H} := L^2(\mathbb{R}^2, dt \, ds),
\end{align}
and $p_t, p_s, q_t, q_s$ be the operators on $\mathscr{H}$ given as
\begin{align}
\label{eq:p_j_and_q_j}
p_t = \pi{\rm i}\hbar \, \frac{\partial}{\partial t}, \quad
p_s = \pi{\rm i}\hbar \, \frac{\partial}{\partial s}, \quad
q_t = \mbox{multiplication by $t$}, \quad
q_s = \mbox{multiplication by $s$},
\end{align}
on a dense subspace of $\mathscr{H}$, e.g. 
\begin{align}
\label{eq:D}
\mathscr{D} := \mathrm{span}_\mathbb{C}\left\{ e^{-\alpha_1 \, t^2 - \alpha_2 \, s^2 + \beta_1 \, t + \beta_2 \, s } P(t,s) \, \left|
\begin{array}{l}
 \alpha_j, \beta_j \in \mathbb{C}, ~ \mathrm{Re}(\alpha_j)>0, ~ \mbox{and $P$ is a }\\
\mbox{polynomial in two variables over $\mathbb{C}$}
\end{array}
\, \right. \right\}.
\end{align}
This space $\mathscr{D}$ is a slight variant of the version as used by Fock and Goncharov e.g. in \cite{FG09}; see also \cite{KS}. It is well known that each of these $p_t,p_s,q_t,q_s$ defines an essentially self-adjoint operator on $\mathscr{D}$, hence has a unique self-adjoint extension in $\mathscr{H}$; let us just denote by the same symbols $p_t,p_s,q_t,q_s$ the corresponding self-adjoint operators, by abuse of notation. They satisfy the Heisenberg relations
\begin{align}
\label{eq:p_j_and_q_k_Heisenberg_relations}
[p_j,p_k]=0=[q_j,q_k],\qquad
[p_j,q_k] = \pi{\rm i}\hbar \cdot \delta_{j,k} \cdot \mathrm{id}, \qquad \mbox{for}\quad j,k\in \{t,s\},
\end{align}
say on $\mathscr{D}$, where $\delta_{j,k}$ is the Kronecker delta; more precisely, they satisfy the \emph{Weyl-relations} version of the above Heisenberg relations. Here, for any two self-adjoint operators $A,B$, by the Weyl-relations version of the Heisenberg commutation relation
$$
[A,B] = {\rm i} \, c \cdot {\rm id}
$$
for some constant $c \in \mathbb{R}$ we mean the family of equalities
$$
e^{{\rm i} \alpha A} e^{{\rm i} \beta B} = e^{-\alpha \beta {\rm i} c} e^{{\rm i} \beta B} e^{{\rm i} \alpha A}, \qquad \forall \alpha,\beta \in \mathbb{R},
$$
of unitary operators, where $e^{{\rm i}\alpha A}$ and $e^{{\rm i}\beta B}$ are defined by the functional calculus of $A,B$. It is well known that any $\mathbb{R}$-linear combination of the above $p_j,q_j$'s is essentially self-adjoint on $\mathscr{D}$, and that they satisfy the Weyl-relations version of the corresponding Heisenberg commutation relations for them; see e.g. \cite{H13}. 

\vs

On the other hand, a densely-defined operator $A$ on a Hilbert space $\mathcal{H}$ is called {\em normal} if
\begin{align}
\label{eq:normal_definition}
\mathrm{Dom}(A) = \mathrm{Dom}(A^*) \quad \mbox{and} \quad ||A\varphi|| = ||A^* \varphi|| \mbox{ for all $\varphi\in \mathrm{Dom}(A)$},
\end{align}
where $\mathrm{Dom}$ stands for the domain of a densely-defined operator. This definition eq.\eqref{eq:normal_definition} of a normal operator is taken from \cite[Def.3.8]{Sch}; it is equivalent to the condition that $A$ is closed and $A^*A = A A^*$, as shown in \cite[Prop.3.25]{Sch}. Here, for a densely-defined operator $A$, its {\em adjoint} $A^*$ is defined on its domain
$$
\mathrm{Dom}(A^*) := \{ \varphi \in \mathcal{H} \, | \, \psi \mapsto \langle A\psi, \varphi \rangle \mbox{ is a continuous linear functional on $\mathrm{Dom}(A)$} \},
$$
by the formula
$$
\langle A\psi,\varphi\rangle = \langle \psi, A^* \varphi\rangle, \qquad \forall \psi \in \mathrm{Dom}(A), 
$$
which determines $A^* \varphi$ uniquely for each $\varphi \in \mathrm{Dom}(A^*)$, thanks to the Riesz Representation Theorem and the dense-ness of $\mathrm{Dom}(A)$. 

\vs

Let $B$ and $C$ be self-adjoint operators on a Hilbert space $\mathcal{H}$ that strongly commute with each other, i.e. the spectral projections of $B$ commute with those of $C$; one equivalent condition is that $B$ and $C$ satisfy the Weyl-relations version of the Heisenberg relation $[B,C]=0$ (see e.g. \cite[Thm.VIII.13]{RS}). Then one can show that the intersection of the domain $\mathrm{Dom}(B)$ of $B$ and the domain $\mathrm{Dom}(C)$ of $C$ is dense. Now define the operator
$$
A = B + {\rm i} \, C
$$
on its dense domain $\mathrm{Dom}(A) := \mathrm{Dom}(B) \cap \mathrm{Dom}(C)$ by the obvious formula. Then its closure is normal, and one notes that $A^*$ equals the closure of $B - {\rm i} \, C$ \footnote{These facts follow e.g. from {\cite[Thm,4.16.(ii)]{Sch}} applied to the two-variable functional calculus for $B,C$.}. The `essentially normal-ness' of operators of the form $B + {\rm i} \, C$ was studied e.g. in \cite{SS1}; in the present paper, by the expression $B+{\rm i}\,C$ we mean its closure. Moreover, there is a well-known functional calculus of normal operators, which allows us to plug a normal operator into a function $\mathbb{C} \to \mathbb{C}$, resulting in a normal operator. One convenient way of thinking about it is as a special case of two-variable functional calculus of two strongly commuting self-adjoint operators; we will implicitly have this in mind throughout the paper. More details on the normal operators and these functional calculus can be found e.g. in \cite[\S4.3, \S5.5]{Sch}. 

\vs

Back to our case. Define the following normal operators on $\mathscr{H}$
\begin{align}
\label{eq:P_and_Q_as_p_and_q}
P := q_t + {\rm i} \, p_s, \qquad
Q := q_s - {\rm i} \, p_t,
\end{align}
constructed by the above discussion; it is possible because $q_t$ and $p_s$ strongly commute with each other, while $q_s$ and $p_t$ strongly commute with each other. Indeed, $q_t,q_s,p_t,p_s$ satisfy the Weyl-relations version of the Heisenberg relations $[q_t,p_s]=0$ and $[q_s,p_t]=0$ (and more) which can be deduced from eq.\eqref{eq:p_j_and_q_k_Heisenberg_relations}.

\vs

In particular, we have
$$
P^* = q_t - {\rm i}\, p_s, \qquad
Q^* = q_s + {\rm i} \, p_t.
$$
We then consider the representation $\pi$ of $\mathcal{C}_{\mathbf{q}}$ on the Hilbert space $\mathscr{H}$:
\begin{align}
\label{eq:representation_q}
\left\{  \begin{array}{ll}
Z_1^{\pm 1} \mapsto \pi(Z_1^{\pm 1}):=\mathbf{Z}_1^{\pm 1} := e^{\pm P}, &
Z_2 \mapsto \pi(Z_2):=\mathbf{Z}_2 := e^Q, \\
(Z_1^*)^{\pm 1} \mapsto \pi((Z_1^*)^{\pm 1}):=(\mathbf{Z}_1^{\pm 1})^* =e^{\pm P^*}, &
Z_1^* \mapsto \pi(Z_2^*):=\mathbf{Z}_2^* =e^{Q^*},    
  \end{array} \right.
\end{align}
where the normal operators $e^P$ and $e^Q$ are constructed by the functional calculus for the normal operators $P$ and $Q$ applied to the function $\lambda \mapsto e^\lambda$. Again, these can also be thought of as
\begin{align}
\label{eq:bf_Z_1_bf_Z_2_explicit}
\mathbf{Z}_1 = e^{q_t + {\rm i}\,p_s}, \qquad
\mathbf{Z}_2 = e^{q_s - {\rm i}\,p_t},
\end{align}
each being obtained by the two-variable functional calculus for two strongly commuting self-adjoint operators. When restricted to the subspace $\mathscr{D}$, one observes that $P,Q,\mathbf{Z}_1^{\pm 1},\mathbf{Z}_2$ and their adjoints preserve $\mathscr{D}$, and satisfy the algebraic relation
\begin{align}
\label{eq:bf_Z_1_bf_Z_2_relation}
\mathbf{Z}_1 \, \mathbf{Z}_2 = \mathbf{q}^2 \, \mathbf{Z}_2 \, \mathbf{Z}_1
\end{align}
when applied to elements of $\mathscr{D}$. 

\vs

A priori, eq.\eqref{eq:bf_Z_1_bf_Z_2_relation} should not be regarded as a full-fledged equation of normal operators; usually, composition and equality of densely-defined operators are delicate matters, because of the domain issues. However, there might be a `stronger' version of eq.\eqref{eq:bf_Z_1_bf_Z_2_relation} that better characterizes the operators $\mathbf{Z}_1$ and $\mathbf{Z}_2$ defined in eq.\eqref{eq:representation_q}, as in the case of the quantum plane; namely, the positive-definite self-adjoint operators $\mathbf{X} = e^{2\pi{\rm i}\,\hbar \frac{d}{d t}}$ and $\mathbf{Y} = e^x$ on $L^2(\mathbb{R},dt)$ are completely characterized up to unitary equivalence by the family of equations of unitary operators
\begin{align}
\label{eq:strong_equations_for_quantum_plane}
\mathbf{X}^{{\rm i}\,\alpha} \, \mathbf{Y}^{{\rm i}\,\beta} = e^{-\alpha\beta\, 2\pi{\rm i}\,\hbar} \, \mathbf{Y}^{{\rm i}\,\beta} \, \mathbf{X}^{{\rm i}\,\alpha}, \qquad \forall \alpha,\beta\in\mathbb{R},
\end{align}
which can be thought of as being a stronger version of the less-well-defined relation $\mathbf{X} \, \mathbf{Y} = e^{2\pi{\rm i}\,\hbar} \, \mathbf{Y} \, \mathbf{X}$. See \cite{Sch92} for this result on the quantum plane algebra representation, which is also used in \cite{FKi}; the positive-definite self-adjoint operators $\mathbf{X}$ and $\mathbf{Y}$ satisfying eq.\eqref{eq:strong_equations_for_quantum_plane} is called an `integrable' representation of the quantum plane $\mathcal{B}_{\exp(\pi{\rm i}\hbar)}$. One way of obtaining such a characterization is to deal with the $\log$ of the generators, i.e. to just characterize $P$ and $Q$, or more easily, to characterize $q_t,p_t,q_s,p_s$:
\begin{proposition}[(higher) Stone-von Neumann Theorem; see e.g. {\cite[Thm.14.8]{H13}}]
Suppose $\wh{q}_t,\wh{p}_t,\wh{q}_s,\wh{p}_s$ are self-adjoint operators on a separable Hilbert space $\mathcal{H}$ that satisfy the Weyl-relations version of the relations $[\wh{q}_t,\wh{q}_s]=0$, $[\wh{p}_s,\wh{p}_t]=0$, and $[\wh{p}_j,\wh{q}_k] = \pi{\rm i}\,\hbar \cdot\delta_{j,k} \cdot \mathrm{id}$ in eq.\eqref{eq:p_j_and_q_k_Heisenberg_relations} for $j,k\in \{t,s\}$, where $\hbar\in \mathbb{R}_{>0}$. Suppose that $\mathcal{H}$ is \emph{irreducible} in the sense that $0$ and $\mathcal{H}$ are the only closed $\mathbb{C}$-vector subspaces of $\mathcal{H}$ that are invariant under each $e^{{\rm i}\, \alpha \, \wh{q}_j}$ and $e^{{\rm i}\,\alpha \, \wh{p}_j}$, $\forall \alpha\in \mathbb{R}$, $\forall j \in \{t,s\}$. Then there exists a unitary map $U : \mathcal{H} \to \mathscr{H} = L^2(\mathbb{R}^2,dt\,ds)$ such that $U \, \wh{q}_j \, U^{-1} = q_j$ and $U\, \wh{p}_j \, U^{-1} = p_j$ for each $j\in\{t,s\}$, with $q_j$ and $p_j$ as in eq.\eqref{eq:p_j_and_q_j}, and this unitary map $U$ is unique up to a multiplicative constant of modulus $1$. \qed
\end{proposition}
For example, if we define the word \emph{integrable} representation of the quantum pseudo-K\"ahler plane $\mathcal{C}_\mathbf{q}$ as the formulas in eq.\eqref{eq:representation_q} and eq.\eqref{eq:P_and_Q_as_p_and_q}, with $q_t,p_t,q_s,p_s$ being self-adjoint operators on a separable Hilbert space satisfying the Weyl-relations version of eq.\eqref{eq:p_j_and_q_k_Heisenberg_relations}, then our representation $(\mathscr{H},\pi)$ would be the unique irreducible integrable representation of $\mathcal{C}_{\bf q}$.  

\vs

Now we shall try to deal with the exponential normal operators $\mathbf{Z}_1^{\pm 1},\mathbf{Z}_2,(\mathbf{Z}_1^{\pm 1})^*,\mathbf{Z}_2^*$, not their $\log$ versions. If one is happy about formal manipulations, one may formally use the BCH formula in eq.\eqref{eq:BCH} to figure out relations like in eq.\eqref{eq:bf_Z_1_bf_Z_2_relation}. Indeed, 
\begin{align}
\label{eq:PQ_commutator}
[P,Q] = [q_t+{\rm i}\,p_s, \, q_s-{\rm i}\,p_t] = [q_t,q_s] + [p_s,p_t] + {\rm i}\,(-[q_t,p_t]+[p_s,q_s]) = -2\pi\hbar \cdot \mathrm{id},
\end{align}
which makes sense either formally, or as operators $\mathscr{D}\to \mathscr{D}$; thus $e^P e^Q = e^{-\pi\hbar} e^{P+Q}$ and $e^Q e^P = e^{\pi\hbar} e^{Q+P}$, hence $e^P e^Q = e^{-2\pi\hbar} e^Q e^P$. The other algebraic relations are obtained similarly by the BCH formula and the following computational results for commutators
\begin{align}
\label{eq:other_commutators}
[P,Q^*]=[P^*,Q]=0, \qquad [P^*,Q^*] = 2\pi\hbar\cdot\mathrm{id}.
\end{align}

\vs

By a property of functional calculus saying that taking the complex conjugate of the function has the effect of taking the adjoint of the resulting normal operator (see e.g. \cite[Thm4.16.(i)]{Sch}), we have
$$
(\mathbf{Z}_1^{\pm 1})^* = e^{\pm P^*} = e^{\pm (q_t - {\rm i}\,p_s)} \qquad\mbox{and}\qquad
\mathbf{Z}_2^* = e^{Q^*} = e^{q_s + {\rm i}\,p_t}.
$$
Anyways, the $*$-structure is preserved by $\pi$, i.e. $\pi(u^*) = \pi(u)^*$ for all generators $u$. Finally, e.g. using the BCH formula, one finds that $\mathbf{Z}_1,\mathbf{Z}_2,\mathbf{Z}_1^*,\mathbf{Z}_2^*$ satisfy all the algebraic relations as in eq.\eqref{eq:C_q_relation}.

\vs

In terms of the genuine algebra elements of $\mathcal{C}_\mathbf{q}$, instead of the $\log$ of generators, we may at the moment be content just with the following weak notion of a representation:
\begin{definition}[densely-defined representation]
\label{def:densely-defined_representation}
A {\em densely-defined representation} of a Hopf $*$-algebra $\mathcal{C}$ over $\mathbb{C}$ is a separable complex Hilbert space together with a dense subspace on which $\mathcal{C}$ acts, preserving the $*$-structure.
\end{definition}
So our $(\mathscr{H},\pi)$ can be viewed as a densely-defined representation of $\mathcal{C}_\mathbf{q}$. From now on, by a representation we would mean a densely-defined representation, unless otherwise noted.

\section{An integrable representation of the modular double $\mathcal{C}_{\mathbf{q},\mathbf{q}^\vee}$}
\label{sec:integral_representation_of_the_modular_double}

We enhance the representation $(\mathscr{H},\pi)$ of the quantum pseudo-K\"ahler plane $\mathcal{C}_\mathbf{q}$ studied in \S\ref{sec:integrable_representation_of_C_q} to a representation of its modular double $\mathcal{C}_{\mathbf{q},\mathbf{q}^\vee}$ in eq.\eqref{eq:C_q_q_vee} considered in \S\ref{sec:quantum_pK_plane_algebra}; we still denote it by $(\mathscr{H},\pi)$. On the space $\mathscr{H}$ defined in eq.\eqref{eq:H}, we assign (densely-defined) normal operators on $\mathscr{H}$ to generators $Z_1^{\pm 1},(Z_1^*)^{\pm 1},Z_2,Z_2^*$ of the subalgebra $\mathcal{C}_\mathbf{q} \cong \mathcal{C}_{\mathbf{q}} \otimes 1 \subset \mathcal{C}_{\mathbf{q},\mathbf{q}^\vee}$ as done in eq.\eqref{eq:representation_q}. To generators $Z_1^\vee, (Z_1^\vee)^*, Z_2^\vee, (Z_2^\vee)^*$ of the sualgebra $\mathcal{C}_{1/\mathbf{q}^\vee} \cong 1\otimes \mathcal{C}_{1/\mathbf{q}^\vee} \subset \mathcal{C}_{\mathbf{q},\mathbf{q}^\vee}$ we assign the normal operators on $\mathscr{H}$ as follows:
\begin{align}
\label{eq:representation_q_vee}
\left\{  \begin{array}{l}
Z_1^\vee \mapsto \pi(Z_1^\vee):=\mathbf{Z}_1^\vee := e^{P/({\rm i}\, \hbar)}, \qquad\qquad
Z_2^\vee \mapsto \pi(Z_2^\vee):=\mathbf{Z}_2^\vee := e^{Q/({\rm i}\, \hbar)}, \\
(Z_1^\vee)^* \mapsto \pi((Z_1^\vee)^*):=(\mathbf{Z}_1^\vee)^* = e^{-P^*/({\rm i}\,\hbar)}, \\
(Z_2^\vee)^* \mapsto \pi((Z_2^\vee)^*):=(\mathbf{Z}_2^\vee)^* = e^{-Q^*/({\rm i}\,\hbar)}.
  \end{array} \right.
\end{align}
Readers can easily deduce the operators for $(Z_1^\vee)^{-1}$ and $((Z_1^\vee)^*)^{-1}$.

\vs

E.g. by using BCH formula in eq.\eqref{eq:BCH} and the commutation relations $[\frac{P}{{\rm i}\,\hbar}, \frac{Q}{{\rm i}\,\hbar}] = 2\pi\frac{1}{\hbar}$, etc, one checks that eq.\eqref{eq:representation_q_vee} indeed provides a representation of the algebra $\mathcal{C}_{1/\mathbf{q}^\vee}$, say on the dense subspace $\mathscr{D}$ \eqref{eq:D} of $\mathscr{H}$. Moreover, one can check that the operators in eq.\eqref{eq:representation_q} of $\mathcal{C}_\mathbf{q}$ and those in eq.\eqref{eq:representation_q_vee} of $\mathcal{C}_{1/\mathbf{q}^\vee}$ `weakly' commute, say as operators on $\mathscr{D}$; for example, since $[P, \frac{Q}{{\rm i}\,\hbar}] = 2\pi{\rm i}$, from the BCH formula in eq.\eqref{eq:BCH} we get
$$
\mathbf{Z}_1 \, \mathbf{Z}_2^\vee = e^P \, e^{Q/({\rm i}\,\hbar)} = e^{2\pi{\rm i}} \, e^{Q/({\rm i}\,\hbar)} \, e^P = \mathbf{Z}_2^\vee \, \mathbf{Z}_1,
$$
either as formally, or as equality of operators on $\mathscr{D}$. This in particular hints that in a certain sense, the representation $\pi$ of $\mathcal{C}_\mathbf{q}$ on $\mathscr{H}$ would not be strongly irreducible, as there exist many operators on $\mathscr{H}$ that weakly commute with the operators representing elements of $\mathcal{C}_\mathbf{q}$ via $\pi$, namely the operators in eq.\eqref{eq:representation_q_vee} representing $\mathcal{C}_{1/\mathbf{q}^\vee}$ via $\pi$. The following conjecture on a certain irreducibility of the representation $(\mathscr{H},\pi)$ is an analog of the corresponding statement for the modular double quantum plane $\mathcal{B}_q \otimes \mathcal{B}_{\til{q}}$; see \cite{FKi}. In the present paper, we will not really depend on this conjecture.
\begin{conjecture}
\label{conj:strong_irreducibility}
We conjecture that the operator algebra $\pi(\mathcal{C}_{1/\mathbf{q}^\vee})$ generated by eq.\eqref{eq:representation_q_vee} is the `full commutant' of the operator algebra $\pi(\mathcal{C}_\mathbf{q})$ generated by eq.\eqref{eq:representation_q} on $\mathscr{H}$, that is to say, if a (densely-defined) linear operator on $\mathscr{H}$ (weakly) commutes with all operators in eq.\eqref{eq:representation_q} and with those in eq.\eqref{eq:representation_q_vee}, then it must be a scalar times the identity operator. 
\end{conjecture}

\vs

On the other hand, like in eq.\eqref{eq:bf_Z_1_bf_Z_2_explicit} one can write these operators more explicitly as
\begin{align}
\left\{ \begin{array}{l}
\mathbf{Z}_1^\vee = e^{\frac{1}{\hbar}\, p_s - {\rm i} \, \frac{1}{\hbar} \, q_t }, \quad\qquad
\mathbf{Z}_2^\vee = e^{ - \frac{1}{\hbar} \, p_t - {\rm i} \, \frac{1}{\hbar} \, q_s}, \\
(\mathbf{Z}_1^\vee)^* = e^{\frac{1}{\hbar}\, p_s + {\rm i} \, \frac{1}{\hbar} \, q_t }, \qquad
(\mathbf{Z}_2^\vee)^* = e^{ - \frac{1}{\hbar} \, p_t + {\rm i} \, \frac{1}{\hbar} \, q_s}.
\end{array} \right.
\end{align}

\section{The modular double compact quantum dilogarithm}
\label{sec:Psi_hbar}

In this section we investigate some properties of the modular double compact quantum dilogarithm function $\Phi^{\pm{\rm i}\hbar}$ defined in eq.\eqref{eq:Psi_hbar_first_time} of  \S\ref{sec:introduction} for $\hbar>0$ which was first systematically studied in the author's previous joint work \cite{KS} with Scarinci on 3d quantum gravity. We begin with the following well-known fact:
\begin{lemma}[compact quantum dilogarithm $\psi^q$]
\label{lem:psi_q}
For each complex number $q$ with $|q|<1$, the infinite product formula in eq.\eqref{eq:psi_q_first_time} for $z\in \mathbb{C}$ converges uniformly on compact subsets of $\mathbb{C}$, defining a meromorphic function $\psi^q(z)$ in the complex $z$-plane, with no zeros and whose poles are at $- q^{-2n-1}$, $n=0,1,2,\ldots$, which are all simple. \qed
\end{lemma}
It is straightforward to observe the functional relation:
$$
\psi^q(q^2 z) = (1+qz) \psi^q(z),
$$
which is a characteristic property of $\psi^q$. As in eq.\eqref{eq:Phi_ihbar_as_ratio}, $\Phi^{\pm {\rm i}\hbar}$ can be expressed as the ratio of two compact quantum dilogarithm functions. Also, as in eq.\eqref{eq:Psi_hbar_first_time}, $\Phi^{\pm {\rm i}\hbar}$ can be expressed as the contour integral, in the style of the usual non-compact quantum dilogarithm $\Phi^\hbar$; in a sense, $\Phi^{\pm {\rm i} \hbar}$ can be thought of as being obtained by analytically continuing the $\hbar$-parameter of Faddeev-Kashaev's non-compact quantum dilogarithm $\Phi^\hbar$. Here are some useful properties of $\Phi^{\pm{\rm i}\hbar}$.

\begin{proposition}[properties of the modular double compact quantum dilogarithm; \cite{KS}]
\label{prop:properties_of_mdqc}
For each $\epsilon \in \{+1,-1\}$ and $\hbar>0$, the function $\Phi^{\epsilon {\rm i}\hbar}(z)$ is meromorphic on the complex plane, and satisfies:
\begin{itemize}
\item[\rm (1)] The zeros and poles are at
\begin{align*}
\mbox{the set of zeros} & = \{(2n+1) \pi{\rm i} - (2m+1) \pi \epsilon \hbar \, | \, n,m\in \mathbb{Z}_{\ge 0} \},\\
\mbox{the set of poles} & = \{-(2n+1) \pi{\rm i} + (2m+1) \pi \epsilon \hbar \, | \, n,m\in \mathbb{Z}_{\ge 0} \}.
\end{align*}
All the zeros and poles are simple.

\item[\rm (2)] (difference equations) Each of the functional relations
\begin{align*}
\left\{
  \begin{array}{lcr}
\Phi^{\epsilon {\rm i} \hbar}(z-2\pi \epsilon \hbar) & = & (1+e^{-\pi \epsilon \hbar} \, e^z) \, \Phi^{\epsilon {\rm i} \hbar}(z), \\
\Phi^{\epsilon {\rm i} \hbar}(z+2\pi{\rm i}) & = & (1+e^{\pi/( \epsilon \hbar)} \, e^{z/(\epsilon{\rm i}\hbar)}) \, \Phi^{\epsilon {\rm i} \hbar}(z)
  \end{array}
\right.
\end{align*}
holds, whenever the two arguments of $\Phi^{\epsilon {\rm i} \hbar}$ are not poles.

\item[\rm (3)] (involutivity) One has
$$
\Phi^{\epsilon {\rm i} \hbar}(z) \, \Phi^{\epsilon {\rm i} \hbar}(-z) = c_{\epsilon {\rm i}\hbar} \, \exp\left( {z^2}/({- 4\pi \epsilon\hbar}) \right),
$$
whenever $z$ and $-z$ are not poles of $\Phi^{\epsilon {\rm i} \hbar}$, where
$$
c_{\epsilon{\rm i}\hbar} := e^{\frac{\pi }{12} \epsilon (\hbar-\hbar^{-1})} \in \mathbb{R}_{>0}.
$$

\end{itemize}
\end{proposition}

\vs

What is actually used in \cite{KS} for the quantization of the moduli space of 3d spacetimes for a positive cosmological constant is the following function $\Psi^\hbar$ which is a certain combination of $\Phi^{\pm {\rm i}\hbar}$, which is sort of a `generalized-complex double' of $\Phi^{{\rm i}\hbar}$ (see \cite{KS} for a detailed discussion of this `double'):
\begin{align*}
& \Psi^\hbar ~:~ \mathbb{C} \to \mathbb{C} \\
& \Psi^\hbar(x+{\rm i}y) := \Phi^{{\rm i}\hbar}(x+{\rm i}y) \, \Phi^{-{\rm i}\hbar}(x-{\rm i}y)
= \Phi^{{\rm i}\hbar}(x+{\rm i}y) / \ol{ \Phi^{{\rm i}\hbar}(x+{\rm i}y) }, \qquad \forall x,y\in \mathbb{R}.
\end{align*}
In particular, one can show that $\Psi^\hbar$ is indeed well-defined on the whole $\mathbb{C}$, and that it is unitary in the sense that
\begin{align}
\label{eq:Psi_unitarity}
|\Psi^\hbar(x+{\rm i}y)| = 1,\quad \forall x,y\in \mathbb{R}.
\end{align}

An immediate corollary of the involutivity of $\Phi^{\pm {\rm i}\hbar}$ as in Prop.\ref{prop:properties_of_mdqc}(3) is the following (see \cite{KS}):
\begin{align}
\label{eq:involutivity_Psi}
\Psi^\hbar(x+{\rm i}y) \Psi^\hbar(-x-{\rm i}y) = \exp( xy/(\pi {\rm i} \hbar) ), \qquad \forall x,y\in \mathbb{R}.
\end{align}
What will be crucially used in the present paper is the operator version of the difference equations as in Prop.\ref{prop:properties_of_mdqc}(2):
\begin{proposition}[the operator version of the difference equations for $\Psi^\hbar$; \cite{KS}]
\label{prop:operator_version_of_difference_equations_for_Psi}
Let ${\bf x}, {\bf y}, {\bf x}', {\bf y}'$ be self-adjoint operators on a separable complex Hilbert space satisfying the Weyl-relations version of the Heisenberg commutation relations 
$$
[{\bf x},{\bf x}']=0=[{\bf y},{\bf y}'] = [{\bf x},{\bf y}] = [{\bf x}',{\bf y}'], \qquad
[{\bf x},{\bf y}'] = \pi {\rm i} \cdot {\rm id} = [{\bf y},{\bf x}'].
$$
\begin{itemize}
\item[\rm (1)] For each $\epsilon \in \{+1,-1\}$, the following operator identity holds:
\begin{align}
\nonumber
\Psi^\hbar({\bf x} + {\rm i} \hbar {\bf y}) e^{{\bf x}' + {\rm i} \epsilon \hbar {\bf y}'} = e^{{\bf x}' + {\rm i} \epsilon \hbar {\bf y}'} ( 1+ e^{-\pi \epsilon \hbar} e^{{\bf x}+{\rm i} \epsilon \hbar {\bf y}}) \Psi^\hbar({\bf x} + {\rm i} \hbar {\bf y})
\end{align}

\item[\rm (2)] For each $\epsilon \in \{+1,-1\}$, the following operator identity holds:
\begin{align}
\nonumber
\Psi^\hbar({\bf x} + {\rm i} \hbar {\bf y}) e^{({\bf x}' + {\rm i} \epsilon \hbar {\bf y}')/(\epsilon {\rm i} \hbar)} = e^{({\bf x}' + {\rm i} \epsilon \hbar {\bf y}')/(\epsilon {\rm i} \hbar)} ( 1+ e^{\pi /(\epsilon \hbar)} e^{({\bf x}+{\rm i} \epsilon \hbar {\bf y})/(\epsilon {\rm i} \hbar)}) \Psi^\hbar({\bf x} + {\rm i} \hbar {\bf y})
\end{align}
\end{itemize}
\end{proposition}
Prop.\ref{prop:operator_version_of_difference_equations_for_Psi}(1) was proved in \cite{KS} by applying the analytic continuation for $\hbar$ to the corresponding operator identity for the usual non-compact quantum dilogarithm $\Phi^\hbar$, which in turn was proved e.g. in \cite{FG09} \cite{G}. We note that Prop.\ref{prop:operator_version_of_difference_equations_for_Psi}(2) can be proved in a similar manner. Another important operator identity for $\Psi^\hbar$ is the pentagon identity, which again was proved in \cite{KS} by applying the analytic continuation for $\hbar$ to the corresponding result for the function $\Phi^\hbar$, which in turn was proved in \cite{G} \cite{FKV} \cite{W}.
\begin{proposition}[the pentagon identity for $\Psi^\hbar$; \cite{KS}]
\label{prop:pentagon_identity_for_Psi}
Let ${\bf x},{\bf y},{\bf x}',{\bf y}'$ be as in Prop.\ref{prop:operator_version_of_difference_equations_for_Psi}. Then
$$
\Psi^\hbar({\bf x}+{\rm i} \hbar {\bf y}) \Psi^\hbar({\bf x}'+{\rm i} \hbar {\bf y}')
= \Psi^\hbar({\bf x}'+{\rm i} \hbar {\bf y}')
\Psi^\hbar(({\bf x}+{\bf x}')+{\rm i} \hbar ({\bf y}+{\bf y}'))
\Psi^\hbar({\bf x}+{\rm i} \hbar {\bf y}).
$$
\end{proposition}
We refer the readers to \cite{KS} for more details about $\Psi^\hbar$.

\section{Decomposing the tensor square representation: operator $\mathbf{F}$}
\label{sec:F}

For an index $j$, which will be ranging in $\{1,2,3\}$ in the present section but can be more general, we denote by $(\mathscr{H}_j, \pi_j)$ the representation on the Hilbert space $\mathscr{H}$ studied in the previous sections. Write
$$
\mathscr{H}_j = L^2(\mathbb{R}^2, dt_j \, ds_j)
$$
as Hilbert spaces, and write the algebra action $\pi_j$ of $\mathcal{C}_{\mathbf{q},\mathbf{q}^\vee}$ as
\begin{align*}
\pi_j(Z_1) =\mathbf{Z}_1 = e^{P_j}, \qquad \pi_j(Z_2) = \mathbf{Z}_2 = e^{Q_j}, \qquad \pi_j(Z_1^*) = \mathbf{Z}_1^* = e^{P^*_j}, \qquad \pi_j(Z_2^*) = \mathbf{Z}_2^* = e^{Q^*_j},
\end{align*}
for the generators of $\mathcal{C}_\mathbf{q}$, where
\begin{align}
\label{eq:P_and_Q_j}
& \hspace{-3mm} P_j = q_{t_j} + {\rm i} \, p_{s_j}, \quad Q_j = q_{s_j} - {\rm i} \, p_{t_j}, \quad
P^*_j = q_{t_j} - {\rm i} \, p_{s_j}, \quad Q^*_j = q_{s_j} + {\rm i} \, p_{t_j}, \\
\label{eq:p_and_q_j}
& \hspace{-3mm} q_{t_j} = t_j, \qquad q_{s_j} = s_j, \qquad p_{t_j} = \pi{\rm i}\hbar\,\frac{\partial}{\partial t_j}, \qquad p_{s_j} = \pi{\rm i}\hbar\,\frac{\partial}{\partial s_j},
\end{align}
and likewise for the generators for $\mathcal{C}_{1/\mathbf{q}^\vee}$ as
\begin{align*}
& \pi_j(Z_1^\vee) = \mathbf{Z}_1^\vee = e^{P_j/({\rm i}\hbar)}, \qquad\qquad
\pi_j(Z_2^\vee) = \mathbf{Z}_2^\vee = e^{Q_j/({\rm i}\hbar)}, \\
& \pi_j((Z_1^\vee)^*) = (\mathbf{Z}_1^\vee)^* = e^{-P_j^*/({\rm i}\hbar)}, \qquad
\pi_j((Z_2^\vee)^*) = (\mathbf{Z}_2^\vee)^* = e^{-Q_j^*/({\rm i}\hbar)}.
\end{align*}

\vs

We will now consider the tensor product representation on the tensor square of $\mathscr{H}$, on which the algebra action is defined using the coproduct, as is on the tensor product of any two representations of a Hopf algebra. We shall see that the representation $\mathscr{H} \otimes \mathscr{H}$ breaks into the direct integral of the same copies of $\mathscr{H}$:
$$
\mathscr{H} \otimes \mathscr{H} \cong \int^{\oplus}_{\mathbb{R}} \mathscr{H}.
$$
The above isomorphism makes only a heuristic sense, and we shall realize this rigorously as the following unitary isomorphism of representations
\begin{align}
\label{eq:sought_for_decomposition_isomorphism}
\mathscr{H} \otimes \mathscr{H} \cong M \otimes \mathscr{H},
\end{align}
where $M$ is the trivial representation of the Hopf algebra $\mathcal{C}_{\mathbf{q},\mathbf{q}^\vee}$ which can be viewed as the {\it multiplicity space}. 

\vs

Now, let us be more precise. Consider the Hilbert space
$$
\mathscr{H} \otimes \mathscr{H} \cong 
\mathscr{H}_1 \otimes \mathscr{H}_2 =L^2(\mathbb{R}^2,dt_1\, ds_1) \otimes L^2(\mathbb{R}^2,dt_2\,ds_2) \cong L^2(\mathbb{R}^4, dt_1 \, ds_1 \, dt_2 \, ds_2),
$$
where $\otimes$ means the tensor product in the category of Hilbert spaces as mentioned in \S\ref{sec:algebraic_story}; we used the canonical isomorphism $L^2(\Omega_1, d\mu_1) \otimes L^2(\Omega_2,d\mu_2) \cong L^2(\Omega_1 \times \Omega_2, d\mu_1 \times d\mu_2)$ of Hilbert spaces, for measure spaces $(\Omega_1,\mu_1)$ and $(\Omega_2,\mu_2)$. The algebra action, denoted by $\pi^{(2)}$, is defined using the coproduct:
\begin{align}
\label{eq:pi_squared_definition}
\pi^{(2)}(u) := (\pi \otimes \pi)(\Delta(u)) = (\pi_1\otimes \pi_2)(\Delta(u)),
\end{align}
for all elements $u$ of the algebra $\mathcal{C}_{\mathbf{q}, \mathbf{q}^\vee}$. So
\begin{align}
\label{eq:pi_12}
\hspace{-7mm}
\left\{ \begin{array}{ll}
\pi^{(2)}(Z_1) = \mathbf{Z}_1 \otimes \mathbf{Z}_1 = e^{P_1} \, e^{P_2} = e^{P_1+P_2}, &
\pi^{(2)}(Z_1^*) = \mathbf{Z}_1^* \otimes \mathbf{Z}_1^* = e^{P_1^* + P_2^*}, \\
\pi^{(2)}(Z_2) = \mathbf{Z}_2 \otimes \mathbf{Z}_1 + 1 \otimes \mathbf{Z}_2 = e^{Q_1+P_2} + e^{Q_2}, \\
\pi^{(2)}(Z_2^*) = \mathbf{Z}_2^* \otimes \mathbf{Z}_1^* + 1 \otimes \mathbf{Z}_2^* = e^{Q_1^*+P_2^*} + e^{Q_2^*},
\end{array} \right.
\end{align}
and likewise
\begin{align}
\label{eq:pi_12_vee}
\hspace{-7mm}
\left\{ \begin{array}{ll}
\pi^{(2)}(Z_1^\vee) = \mathbf{Z}_1^\vee \otimes \mathbf{Z}_1^\vee = e^{(P_1+P_2)/({\rm i}\hbar)}, \\
\pi^{(2)}((Z_1^\vee)^*) = (\mathbf{Z}_1^\vee)^* \otimes (\mathbf{Z}_1^\vee)^* = e^{-(P_1^* + P_2^*)/({\rm i}\hbar)}, \\
\pi^{(2)}(Z_2^\vee) = \mathbf{Z}_2^\vee \otimes \mathbf{Z}_1^\vee + 1 \otimes \mathbf{Z}_2^\vee = e^{(Q_1+P_2)/({\rm i}\hbar)} + e^{Q_2/({\rm i}\hbar)}, \\
\pi^{(2)}((Z_2^\vee)^*) = (\mathbf{Z}_2^\vee)^* \otimes (\mathbf{Z}_1^\vee)^* + 1 \otimes (\mathbf{Z}_2^\vee)^* = e^{-(Q_1^*+P_2^*)/({\rm i}\hbar)} + e^{-Q_2^*/({\rm i}\hbar)}.
\end{array} \right.
\end{align}

In fact one has to be careful about domains of these operators $\pi^{(2)}(u)$.

\vs

We now turn into the multiplicity module $M$.
\begin{definition}[the trivial representation]
\label{def:trivial_representation}
Let $M$ be the Hilbert space
$$
L^2(\mathbb{R}^2, dt\,ds),
$$
and equip it with the trivial action $\pi_{\mathrm{triv}}$ of the Hopf algebra $\mathcal{C}_{\mathbf{q},\mathbf{q}^\vee}$, i.e.
\begin{align}
\label{eq:trivial_action}
\pi_{\mathrm{triv}}(u) = \epsilon(u) \cdot \mathrm{id},
\end{align}
for all elements $u$ of $\mathcal{C}_{\mathbf{q},\mathbf{q}^\vee}$, where $\epsilon$ is the counit of $\mathcal{C}_{\mathbf{q},\mathbf{q}^\vee}$.
\end{definition}

\begin{lemma}
\label{lem:trivial_representation}
Let $(M,\pi_{\mathrm{triv}})$ be as in Def.\ref{def:trivial_representation}, and $(V,\pi_V)$ be a representation of $\mathcal{C}_{\mathbf{q},\mathbf{q}^\vee}$. Then the algebra actions on $M{\otimes} V$ and $V {\otimes} M$ of elements $u$ of $\mathcal{C}_{\mathbf{q},\mathbf{q}^\vee}$ defined by the coproduct coincide respectively with $1\otimes \pi_V(u)$ and $\pi_V(u) \otimes 1$. In particular, the vector space isomorphism
\begin{align}
\label{eq:bf_P}
\mathbf{P} ~ : ~  M{\otimes} V \to V{\otimes} M, \quad \varphi\otimes \psi \mapsto \psi\otimes \varphi,
\end{align}
is an intertwining map.
\end{lemma}
{\it Proof.} This also works when $M$ is replaced by any trivial representation of the algebra, i.e. when the action on it is defined as in eq.\eqref{eq:trivial_action}. The action on $M{\otimes}V$ defined via the coproduct is $u \mapsto (\pi_{\mathrm{triv}}\otimes\pi_V)(\Delta(u)) = (1\otimes \pi_V)((\epsilon\otimes 1)(\Delta(u))) = (1\otimes \pi_V)(1\otimes u) = 1\otimes \pi_V(u)$; we used $(\epsilon\otimes 1)(\Delta(u)) = 1\otimes u$ which is one of the axioms of Hopf algebras (see e.g. \cite{CP}). One can likewise show that the action on $V{\otimes}M$ is $\pi_V(u)\otimes 1$ using $(1\otimes \epsilon)(\Delta(u)) = u\otimes 1$. \qed

\vs

What we aim to do is to realize the isomorphism in eq.\eqref{eq:sought_for_decomposition_isomorphism} by some unitary map
\begin{align}
\label{eq:bf_F_as_arrow}
\mathbf{F} ~ : ~ \mathscr{H}{\otimes}\mathscr{H} ~ \to ~ M{\otimes}\mathscr{H}.
\end{align}
We present one answer:
\begin{align}
\label{eq:bf_F}
\mathbf{F} := e^{\frac{-1}{2\pi\hbar}( P_1 Q_2 - P_1^* Q_2^*)} \, (\Phi^{{\rm i}\hbar}(Q_1+P_2-Q_2))^{-1} \, (\Phi^{-{\rm i}\hbar} (Q_1^* + P_2^* - Q_2^*) )^{-1}.
\end{align}
A better way of writing the second and the third factors is using the function $\Psi^\hbar$, using the functional calculus of the normal operator $Q_1+P_2-Q_2$:
\begin{align}
\label{eq:bf_F2}
\mathbf{F} = e^{\frac{-1}{2\pi\hbar}( P_1 Q_2 - P_1^* Q_2^*)} \, (\Psi^{\hbar}(Q_1+P_2-Q_2))^{-1}
\end{align}
As mentioned before, one can view the operator $\Psi^{\hbar}(Q_1+P_2-Q_2)$ as being defined as
$$
\Psi^{\hbar}( (q_{s_1}+q_{t_2}-q_{s_2}) + {\rm i} (- p_{t_1} + p_{s_2} + p_{t_2}) )
$$
using the two-variable functional calculus for the strongly commuting self-adjoint operators $q_{s_1}+q_{t_2}-q_{s_2}$ and $- p_{t_1} + p_{s_2} + p_{t_2}$. In particular, $\Psi^{\hbar}(Q_1+P_2-Q_2)$ is unitary, due to eq.\eqref{eq:Psi_unitarity}.

\vs

In order to justify our answer ${\bf F}$, we first investigate the first factor in detail. Let
\begin{align}
\label{eq:S_12}
\mathbf{S} = \mathbf{S}_{12} := e^{\frac{-1}{2\pi\hbar}( P_1 Q_2 - P_1^* Q_2^*)}.
\end{align}
The expression $e^{\frac{-1}{2\pi\hbar}(P_1Q_2 - P_1^* Q_2^*)}$ can be viewed as only being heuristic, and one rigorous way of defining it is to regard it as
\begin{align}
\label{eq:the first_factor_of_F_differently}
e^{\frac{-{\rm i}}{\pi\hbar} (p_{s_1} q_{s_2} - q_{t_1} p_{t_2})},
\end{align}
which in turn can be thought of as $e^{\frac{-{\rm i}}{\pi\hbar}\,p_{s_1}q_{s_2}} \, e^{\frac{{\rm i}}{\pi\hbar} \,q_{t_1}p_{t_2}}$, where each of $e^{\frac{-{\rm i}}{\pi\hbar}\,p_{s_1}q_{s_2}}$ and $e^{\frac{{\rm i}}{\pi\hbar} \,q_{t_1}p_{t_2}}$ is obtained by the two-variable functional calculus for two strongly commuting self-adjoint operators (\cite{Sch}), the first one for $p_{s_1}$ and $q_{s_2}$ while the second for $q_{t_1}$ and $p_{t_2}$. In particular, $e^{\frac{-{\rm i}}{\pi\hbar}\,p_{s_1}q_{s_2}}$ and $e^{\frac{{\rm i}}{\pi\hbar} \,q_{t_1}p_{t_2}}$ are unitary, hence so is their composition. It is sensible to write this operator as $e^{\frac{-1}{2\pi\hbar} (P_1 Q_2 - P_1^* Q_2^*)}$ because $P_1Q_2 - P_1^* Q_2^* = 2{\rm i}(p_{s_1} q_{s_2} - q_{t_1} p_{t_2})$ holds; this can either be thought of as an equality of operators $\mathscr{D} \to \mathscr{D}$, or as the equality ${\rm i}(P_1Q_2 - P_1^* Q_2^*) = -2(p_{s_1}q_{s_2} - q_{t_1}p_{t_2})$ of self-adjoint operators, each side being obtained by the four-variable functional calculus for four mutually strongly commuting self-adjoint operators $p_{s_1},q_{s_2},q_{t_1},p_{t_2}$. In fact, one can be even more explicit, as follows. Note first that 
\begin{align}
\label{eq:basic_shift}
(e^{\frac{-{\rm i}}{\pi\hbar} p_{s_1} q_{s_2}}f)(t_1,s_1,t_2,s_2) = f(t_1,s_1+s_2,t_2,s_2), \qquad \forall f\in L^2(\mathbb{R}^4, \, dt_1\,ds_1\,dt_2\,ds_2),
\end{align}
i.e. $e^{\frac{-{\rm i}}{\pi\hbar} p_{s_1} q_{s_2}}$ is a certain kind of shift operator. A formal way of obtaining eq.\eqref{eq:basic_shift} is first to observe $\frac{-{\rm i}}{\pi\hbar} p_{s_1} q_{s_2} =  s_2 \, \frac{\partial}{\partial s_1} $ using eq.\eqref{eq:p_and_q_j}, then to expand formally $e^{s_2\,\frac{\partial}{\partial s_1}}f = f + s_2 \frac{\partial f}{\partial s_1} + \frac{s_2^2}{2!} \frac{\partial^2 f}{\partial s_1^2} + \cdots$, and finally recall the Taylor series expansion from calculus. A rigorous proof is not difficult either; see e.g. the arguments given in \cite[\S6.4]{KS}. Likewise, we have $(e^{\frac{{\rm i}}{\pi\hbar} q_{t_1} p_{t_2}}f)(t_1,s_1,t_2,s_2) = f(t_1,s_1,t_2-t_1,s_2)$. Combining, we obtain
\begin{align}
\label{eq:quadratic_part_explicit}
(e^{\frac{-1}{2\pi\hbar}(P_1Q_2-P_1^*Q_2^*)}f)(t_1,s_1,t_2,s_2) = f(t_1,s_1+s_2,t_2-t_1,s_2), ~ \forall f\in L^2(\mathbb{R}^4,\,dt_1ds_1dt_2ds_2).
\end{align}
That is, ${\bf S}_{12}$ is an example of an operator on $L^2(\mathbb{R}^n)$ induced by an invertible linear map $\mathbb{R}^n\to \mathbb{R}^n$. Basic properties of such maps are systematically studied in \cite{Kim_phase}, which we recall also for later use.

\begin{definition}[special linear operators on $L^2(\mathbb{R}^n)$; \cite{Kim_phase}]
\label{def:saso}
Let $\mathrm{SL}_\pm(n,\mathbb{R})$ be the group of all $n\times n$ real matrices of determinant $\pm 1$, acting from the right on the space $\mathbb{R}^n$ of real $n$-dimensional row vectors via the usual multiplication of a row vector and a matrix. For each $\mathbf{c} \in \mathrm{SL}_\pm(n,\mathbb{R})$ define the {\em special linear} operator $\mathbf{S}_\mathbf{c}$ on $L^2(\mathbb{R}^n, \bigwedge_{i=1}^n da_i)$ as
$$
(\mathbf{S}_\mathbf{c}f)(\mathbf{a}) = f( \mathbf{a} \, \mathbf{c} ), \qquad \forall \mathbf{a}\in\mathbb{R}^n,
$$
where $\mathbf{a}$ stands for the row vector $(a_1 \, a_2 \, \cdots \, a_n)$. 
\end{definition}

\begin{proposition}[see \cite{Kim_phase}]
\label{prop:saso}
$\mathbf{S}_\mathbf{c}$ is unitary. The assignment $\mathbf{c} \mapsto \mathbf{S}_\mathbf{c}$ is an injective group homomorphism from $\mathrm{SL}_\pm(n,\mathbb{R})$ to the group of all unitary maps on $L^2(\mathbb{R}^n)$. That is, we have
\begin{align}
\label{eq:saso_multiplicativity}
\mathbf{S}_\mathbf{c} \, \mathbf{S}_{\mathbf{c}'} = \mathbf{S}_{\mathbf{c} \, \mathbf{c}'}.
\end{align}

\vs

One has $\mathbf{S}_{\mathbf{id}}=\mathrm{id}$ and $\mathbf{S}_\mathbf{c}^{-1} = \mathbf{S}_{\mathbf{c}^{-1}}$. \qed
\end{proposition}
We recognize our ${\bf S} = {\bf S}_{12} =e^{\frac{-1}{2\pi\hbar}(P_1Q_2 - P_1^* Q_2^*)}$ as being of this form.
\begin{lemma}
\label{lem:S_12_as_saso}
One has
$$
e^{\frac{-1}{2\pi\hbar}(P_1Q_2 - P_1^* Q_2^*)} = \mathbf{S}_{\left(\begin{smallmatrix} 1 & 0 & -1 & 0 \\ 0 & 1 & 0 & 0 \\ 0 & 0 & 1 & 0 \\ 0 & 1 & 0 & 1 \end{smallmatrix} \right)} \quad\mbox{on}\quad L^2(\mathbb{R}^4, \, dt_1 \, ds_1\, dt_2 \, ds_2). \qed
$$
\end{lemma}
What will be convenient is the conjugation action of ${\bf S}_{\bf c}$ on the position and the momentum operators.
\begin{lemma}[the conjugation action of ${\bf S}_{\bf c}$; {\cite[Lem.3.18]{Kim_phase}}]
\label{lem:conjugation_action_of_S_c}
For ${\bf c} = (c_{ij})_{i,j} \in {\rm SL}_\pm (n,\mathbb{R})$, write ${\bf c}^{-1} = (c^{ij})_{i,j}$. Then one has
\begin{align}
\label{eq:saso_conjugation_on_p}
\textstyle {\bf S}_{\bf c} \, {\rm i} \frac{\partial}{\partial t_i} \, {\bf S}_{\bf c}^{-1} & \textstyle = \sum_{j=1}^n c^{ij} ({\rm i} \frac{\partial}{\partial t_j}), \\ 
\label{eq:saso_conjugation_on_q}
\textstyle {\bf S}_{\bf c} \, t_i \, {\bf S}_{\bf c}^{-1} & \textstyle = \sum_{j=1}^n c_{ji} t_j,
\end{align}
which are understood as equalities of self-adjoint operators, or as those of operators $\mathscr{D} \to \mathscr{D}$.
\end{lemma}
We shall also use the basic fact from functional analysis that conjugation by a unitary operator $U$ commutes with the functional calculus, i.e. 
$$
U\, f(A) \, U^{-1} = f(UAU^{-1})
$$
either when $A$ is a self-adjoint operator and $f : \mathbb{R} \to \mathbb{C}$ is a measurable function, or when $A$ is a normal operator and $f : \mathbb{C} \to \mathbb{C}$ is a measurable function (see \cite{RS}).  

\vs

We are now ready to prove the sought-for property of ${\bf F}$:
\begin{proposition}[a decomposition map of the tensor square]
\label{prop:decomposition_of_tensor_square}
If both the Hilbert spaces $\mathscr{H}{\otimes}\mathscr{H}$ and $M{\otimes}\mathscr{H}$ are realized as
\begin{align*}
L^2(\mathbb{R}^2,dt_1\,ds_1) {\otimes} L^2(\mathbb{R}^2, dt_2\,ds_2) \cong L^2(\mathbb{R}^4, dt_1\,ds_1\,dt_2\,ds_2).
\end{align*}
Then the unitary map
\begin{align}
\nonumber
\mathbf{F} ~ : ~ \mathscr{H} {\otimes}\mathscr{H} ~ \to ~ M {\otimes}\mathscr{H}
\end{align}
defined in eq.\eqref{eq:bf_F2} intertwines the action of the algebra on the two representations $\mathscr{H} {\otimes}\mathscr{H} $ and $M {\otimes}\mathscr{H}$; that is, the intertwining equation
\begin{align}
\label{eq:intertwining_equation}
  \mathbf{F} \, \pi^{(2)}(u) = (1\otimes \pi(u)) \, \mathbf{F}
\end{align}
holds for all elements $u$ of the algebra $\mathcal{C}_{\mathbf{q},\mathbf{q}^\vee}$.
\end{proposition}

{\it Proof.} As the first step we claim the following conjugation actions of ${\bf S}={\bf S}_{12} = e^{\frac{-1}{2\pi\hbar}(P_1Q_2-P_1^*Q_2^*)}$:
\begin{align}
\label{eq:S_conjugation_on_P2Q2}
{\bf S}^{-1} \, P_2 \, {\bf S} = P_1 + P_2, \qquad
{\bf S}^{-1} \, Q_2 \, {\bf S} = Q_2.
\end{align}
A formal and heuristic way of seeing this is to use
\begin{align}
\label{eq:technical_conjugation}
e^A \, B \, e^{-A} = B + [A,B], \qquad \mbox{if $[A,B]$ commutes with both $A$ and $B$,}
\end{align}
and the observations $[P_2, \frac{-1}{2\pi \hbar}(P_1 Q_2 - P_1^* Q_2^*)] = P_1$ and $[Q_2, \frac{-1}{2\pi \hbar} (P_1 Q_2 - P_1^* Q_2^*)]=0$. Later, one may also formally use the functional calculus version of this conjugation statement
\begin{align}
\label{eq:technical_conjugation3}
e^A \, f(B)  = f(B+[A,B]) \, e^A, \qquad f(B) \, e^A  = e^A \, f(B+[B,A]).
\end{align}
For our setting, one can prove the above statements rigorously by using Lemmas \ref{lem:S_12_as_saso} and \ref{lem:conjugation_action_of_S_c}, which we leave to the readers. 

\vs

For the intertwining equation for $u = Z_1$, observe
\begin{align*}
\mathbf{F}^{-1} \, (1\otimes \pi(Z_1)) \, \mathbf{F} & = \Psi^\hbar(Q_1 + P_2 - Q_2) \, \ul{ {\bf S}^{-1} \, e^{P_2} \, {\bf S} } \, (\Psi^\hbar(Q_1 + P_2 - Q_2))^{-1} \qquad \mbox{($\because$ eq.\eqref{eq:bf_F2}, \eqref{eq:S_12}, \eqref{eq:representation_q})} \\
& = \Psi^\hbar(Q_1+P_2-Q_2) \, e^{P_1+P_2} \, (\Psi^\hbar(Q_1+P_2-Q_2))^{-1} \qquad \mbox{($\because$ eq.\eqref{eq:S_conjugation_on_P2Q2})} \\
& \stackrel{\vee}{=} e^{P_1+P_2}  \stackrel{{\rm eq.}\eqref{eq:pi_12}}{=} \pi_{12}(Z_1),
\end{align*}
as desired, where the checked equality holds because $P_1+P_2$ strongly commutes with $Q_1+P_2-Q_2$. For the intertwining equation for $u=Z_2$, note
\begin{align*}
{\bf F}^{-1} \, (1\otimes \pi(Z_2)) \, {\bf F} & = \Psi^\hbar(Q_1+P_2-Q_2) \, \ul{ {\bf S}^{-1} \, e^{Q_2} \, {\bf S} } \, (\Psi^\hbar(Q_1+P_2-Q_2))^{-1} \qquad \mbox{($\because$ eq.\eqref{eq:bf_F2}, \eqref{eq:S_12}, \eqref{eq:representation_q})} \\
& = \Psi^\hbar(Q_1+P_2-Q_2) \, e^{Q_2} \, (\Psi^\hbar(Q_1+P_2-Q_2))^{-1} \qquad \mbox{($\because$ eq.\eqref{eq:S_conjugation_on_P2Q2})} \\
& \stackrel{\vee}{=} e^{Q_2} + e^{Q_2+(Q_1+P_2-Q_2)} = e^{Q_1+P_2} + e^{Q_2} \stackrel{\eqref{eq:pi_12}}{=} \pi_{12}(Z_2),
\end{align*}
as desired, where the checked equality holds due to Prop.\ref{prop:operator_version_of_difference_equations_for_Psi}(1), which is the operator version of a difference equation for the modular double compact quantum dilogarithm $\Psi^\hbar$. It is a straightforward exercise to check that $Q_1+P_2-Q_2 = (q_{s_1} + q_{t_2} - q_{s_2}) + {\rm i} (-p_{t_1} + p_{s_2} + p_{t_2})$ and $Q_2 = q_{s_2} - {\rm i} p_{t_2}$ indeed satisfy the conditions of Prop.\ref{prop:operator_version_of_difference_equations_for_Psi}. Similarly, in particular by using the remaining version of the operator version of the difference equation as in Prop.\ref{prop:operator_version_of_difference_equations_for_Psi}(2), one can show the intertwining equations for $u = Z_1^\vee$ and $u=Z_2^\vee$ too.  \qed

\begin{remark}
In fact, the intertwining equations cannot be asked to hold for all vectors in the entire Hilbert space, and should only be asked for the vectors in some suitable dense subspaces, in the style of Fock-Goncharov's Schwartz spaces. We refer the readers to \cite{FG09} for more details. For the purpose of the present paper, we will not worry so much about the precise domains.
\end{remark}

\section{Decomposing the tensor cube: operator $\mathbf{T}$}
\label{sec:T}

On the tensor cube $\mathscr{H}_1 {\otimes} \mathscr{H}_2 {\otimes} \mathscr{H}_3$, where each $\mathscr{H}_j$ is isomorphic to the representation $\mathscr{H}$, one defines the algebra action $\pi_{123}$ as
$$
\pi_{123}(u) := (\pi_1\otimes \pi_2\otimes \pi_3)( (\Delta \otimes 1)(\Delta(u))) = (\pi_1\otimes \pi_2\otimes \pi_3)( (1\otimes \Delta)(\Delta(u))),
$$
for all algebra elements $u\in \mathcal{C}_{\mathbf{q},\mathbf{q}^\vee}$, where the latter equality is by the coassociativity $(\Delta \otimes 1) \circ \Delta = (1\otimes\Delta)\circ \Delta$ of the coproduct. Choosing to use $(\Delta \otimes 1)\circ \Delta$ can be encoded as the parenthesizing $(\mathscr{H}_1 {\otimes}\mathscr{H}_2){\otimes}\mathscr{H}_3$, while $(1\otimes \Delta)\circ \Delta$ as $\mathscr{H}_1 {\otimes} (\mathscr{H}_2{\otimes} \mathscr{H}_3)$. Let us decompose this representation into a direct integral of $\mathscr{H}$, or more precisely into the tensor product of a trivial multiplicity module and $\mathscr{H}$. We can do so by iterating the decomposition of the tensor square of $\mathscr{H}$ obtained in the previous section. To keep track of the factors we write eq.\eqref{eq:bf_F_as_arrow} as
$$
\mathbf{F} ~ : ~ \mathscr{H}_j{\otimes}\mathscr{H}_k ~ \to ~ M_{jk}^l {\otimes}\mathscr{H}_l.
$$

\begin{figure}[htbp!]
\centering
\begin{pspicture}[showgrid=false](0,0)(3.3,2.8)
\rput[bl](0,-0.5){
\PstTriangle[unit=1.5,PolyName=P]
\pcline(P1)(P2)\ncput*{1}
\pcline(P2)(P3)\ncput*{2}
\pcline(P3)(P1)\ncput*{3}
\rput[l]{30}(P2){\hspace{1,7mm}$\bullet$}
}
\end{pspicture}
\caption{A triangle representing $M_{12}^3$, where $\mathscr{H}_1 {\otimes} \mathscr{H}_2 \cong M_{12}^3 \otimes \mathscr{H}_3$}
\label{fig:one_triangle}
\end{figure}
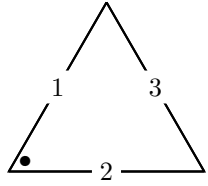

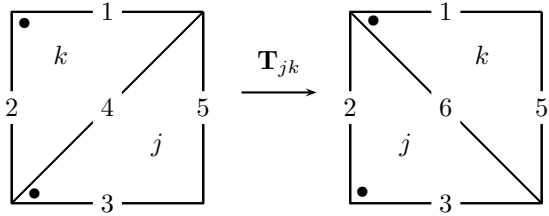
\begin{figure}[htbp!]
\centering
\begin{pspicture}[showgrid=false](0,0)(8.4,3.5)
\rput[bl](0;0){
\PstSquare[unit=1.8,PolyName=P]
\pcline(P1)(P2)\ncput*{1}
\pcline(P2)(P3)\ncput*{2}
\pcline(P3)(P4)\ncput*{3}
\pcline(P4)(P1)\ncput*{5}
\pcline(P1)(P3)\ncput*{4}
\rput[l]{-45}(P2){\hspace{1,4mm}$\bullet$}
\rput[l]{22}(P3){\hspace{2,4mm}$\bullet$}
}
\rput[l](2.5,1.3){$j$}
\rput[l](1.2,2.5){$k$}
\rput[bl](4.5;0){
\PstSquare[unit=1.8,PolyName=P]
\pcline(P1)(P2)\ncput*{1}
\pcline(P2)(P3)\ncput*{2}
\pcline(P3)(P4)\ncput*{3}
\pcline(P4)(P1)\ncput*{5}
\pcline(P2)(P4)\ncput*{6}
\rput[l]{-23}(P2){\hspace{2,4mm}$\bullet$}
\rput[l]{41}(P3){\hspace{1,4mm}$\bullet$}
}
\rput[l](5.8,1.3){$j$}
\rput[l](6.8,2.5){$k$}
\rput[l](3.7,2){\pcline{->}(0,0)(1;0)\Aput{${\bf T}_{jk}$}}
\end{pspicture}
\caption{The move representing ${\bf T} : M_{43}^5 {\otimes} M_{12}^4 \to M_{23}^6 {\otimes} M_{16}^5$}
\label{fig:T}
\end{figure}

The two decompositions of $\mathscr{H}_1 {\otimes} \mathscr{H}_2 {\otimes} \mathscr{H}_3$ are presented as follows:
\begin{align}
\label{eq:intertwining_diagram}
\begin{array}{l}
\xymatrix{
(\mathscr{H}_1 {\otimes} \mathscr{H}_2) {\otimes} \mathscr{H}_3 \ar[r]^-{\mathbf{F}_{12}} \ar[d]^{\mathrm{id}} &
M_{12}^4 {\otimes} \mathscr{H}_4 {\otimes} \mathscr{H}_3 \ar[r]^-{\mathbf{F}_{23}} &
M_{12}^4 {\otimes} M_{43}^5 {\otimes} \mathscr{H}_5 \ar[r]^-{\mathbf{P}_{(12)}} &
M_{43}^5 {\otimes} M_{12}^4 {\otimes} \mathscr{H}_5 \ar@{.>}[d] \\
\mathscr{H}_1 {\otimes} (\mathscr{H}_2 {\otimes} \mathscr{H}_3) \ar[r]^-{\mathbf{F}_{23}} &
\mathscr{H}_1 {\otimes} M_{23}^6 {\otimes} \mathscr{H}_6 \ar[r]^-{\mathbf{P}_{(12)}} &
M_{23}^6 {\otimes} \mathscr{H}_1 {\otimes} \mathscr{H}_6 \ar[r]^-{\mathbf{F}_{23}} &
M_{23}^6 {\otimes} M_{16}^5 {\otimes} \mathscr{H}_5,
}
\end{array}
\end{align}
where, $\mathbf{F}_{ab}$ means that $\mathbf{F}$ being applied to the $a$-th and the $b$-th tensor factors, which do not necessarily coincide with the subscript indices of the corresponding representation $\mathscr{H}_j$, while $\mathbf{P}_{(a\,b)} = \mathbf{P}_{(a\,b)}^{-1}$ means the transposition permutation of the $a$-th and the $b$-th tensor factors. From Prop.\ref{prop:decomposition_of_tensor_square} and Lem.\ref{lem:trivial_representation} we see that all solid arrows in the above diagram are unitary isomorphisms that intertwine the algebra actions. The reason why the very last $\mathscr{H}$ is labeled by $\mathscr{H}_5$ for both decompositions is because we will encode the decomposition $\mathscr{H}_1 {\otimes} \mathscr{H}_2 \to M_{12}^3 {\otimes} \mathscr{H}_3$ as the dotted triangle as in Fig.\ref{fig:one_triangle}, and thus the two decompositions $(\mathscr{H}_1 {\otimes} \mathscr{H}_2) {\otimes} \mathscr{H}_3 \to M_{43}^5 {\otimes} M_{12}^4 {\otimes} \mathscr{H}_5$ and $\mathscr{H}_1 {\otimes} (\mathscr{H}_2 {\otimes} \mathscr{H}_3) \to M_{23}^6 {\otimes} M_{16}^5 {\otimes} \mathscr{H}_5$ respectively as the dotted quadrilaterals as in Fig.\ref{fig:T} \footnote{The Figures \ref{fig:one_triangle} and \ref{fig:T} are taken from \cite{FKi}.}.

\vs

Now we will construct the map $M_{43}^5 {\otimes} M_{12}^4 {\otimes} \mathscr{H}_5 \to M_{23}^6 {\otimes} M_{16}^5 {\otimes} \mathscr{H}_5$, appearing as the rightmost vertical dotted arrow in the diagram in eq.\eqref{eq:intertwining_diagram}, as the unique map that makes this diagram to commute; it is $\mathbf{F}_{23} \, \mathbf{P}_{(12)} \, \mathbf{F}_{23} \, \mathbf{F}_{12}^{-1} \, \mathbf{F}_{23}^{-1} \, \mathbf{P}_{(12)}^{-1}$, which can be read directly from the arrows of the diagram in eq.\eqref{eq:intertwining_diagram}. To simplify this expression a bit, notice that $\mathbf{P}_{(12)} \, (\sim) \, \mathbf{P}_{(12)}^{-1}$ is same as $(\sim)$ with all subscripts $1$ and $2$ exchanged. So $\mathbf{F}_{23} \, \mathbf{P}_{(12)} \, \mathbf{F}_{23} \, \mathbf{F}_{12}^{-1} \, \mathbf{F}_{23}^{-1} \, \mathbf{P}_{(12)}^{-1} = \mathbf{F}_{23} \, \mathbf{F}_{13} \, \mathbf{F}_{21}^{-1} \, \mathbf{F}_{13}^{-1}$.
\begin{proposition}
\label{prop:construction_of_T}
The map for the rightmost vertical dotted arrow in the diagram in eq.\eqref{eq:intertwining_diagram} is
$$
\mathbf{F}_{23} \, \mathbf{F}_{13} \, \mathbf{F}_{21}^{-1} \, \mathbf{F}_{13}^{-1} = \mathbf{T}_{12} = \mathbf{T} \otimes 1 ~ : ~ M_{43}^5 {\otimes} M_{12}^4 {\otimes} \mathscr{H}_5 \to M_{23}^6 {\otimes} M_{16}^5 {\otimes} \mathscr{H}_5
$$
for some uniquely-determined unitary map
$$
\mathbf{T} : M {\otimes} M \to M {\otimes} M.
$$
\end{proposition}

In particular, this operator $\mathbf{F}_{23}\mathbf{F}_{13} \mathbf{F}_{21}^{-1}\mathbf{F}_{13}^{-1}$, which a priori acts on a triple tensor product, only acts on the first two tensor factors $M$, but not on the third factor $\mathscr{H}$. A heuristic proof for an existence of such $\mathbf{T}$ could be obtained by using the fact that the above map $M \otimes M \otimes \mathscr{H} \to M \otimes M \otimes \mathscr{H}$ is an intertwiner, together with the type of irreducibility of the representation on $\mathscr{H}$ in Conjecture \ref{conj:strong_irreducibility}, which tells us that this map $M \otimes M \otimes \mathscr{H} \to M \otimes M \otimes \mathscr{H}$ acts as scalar on the third tensor factor. We shall prove the above proposition by an explicit computation, but let us postpone it a little bit.

\vs

To summarize, this unitary map
$$
\mathbf{T} ~ : ~ M_{43}^5 {\otimes} M_{12}^4 \to M_{23}^6 {\otimes} M_{16}^5
$$
encodes the parenthesizing change $(\mathscr{H}_1{\otimes}\mathscr{H}_2){\otimes}\mathscr{H}_3 \leadsto \mathscr{H}_1 {\otimes} (\mathscr{H}_2 {\otimes} \mathscr{H}_3)$, in the sense of the commutative diagram in eq.\eqref{eq:intertwining_diagram}.

\vs

Now, consider the tensor fourth power representation $\mathscr{H} {\otimes} \mathscr{H} {\otimes} \mathscr{H} {\otimes} \mathscr{H}$. The action is defined via the coproduct, which is well-defined because of the coassociativity. We now decompose this into a direct integral of $\mathscr{H}$, or again, into the tensor product of a trivial multiplicity module and $\mathscr{H}$. There are five ways of decomposing step by step using $\mathbf{F} : \mathscr{H} {\otimes} \mathscr{H} \to M{\otimes}\mathscr{H}$, which are encoded by the five choices of parenthezing, as presented in the following diagram
\begin{align*}
\begin{array}{c}
((\mathscr{H}_1 {\otimes}  \mathscr{H}_2) {\otimes} \mathscr{H}_3) {\otimes} \mathscr{H}_4 ~ \to ~
(\mathscr{H}_1 {\otimes} (\mathscr{H}_2 {\otimes} \mathscr{H}_3)) {\otimes} \mathscr{H}_4 ~ \to ~
\mathscr{H}_1 {\otimes} ( (\mathscr{H}_2 {\otimes} \mathscr{H}_3 ) {\otimes} \mathscr{H}_4) \\
\searrow \hspace{70mm} \swarrow \\
(\mathscr{H}_1 {\otimes} \mathscr{H}_2) {\otimes} (\mathscr{H}_3 {\otimes} \mathscr{H}_4) ~ \to ~
\mathscr{H}_1 {\otimes} (\mathscr{H}_2 {\otimes} (\mathscr{H}_3 {\otimes} \mathscr{H}_4))
  \end{array}
\end{align*}
where all arrows are identity maps; hence all maps are unitary $\mathcal{C}_{\mathbf{q},\mathbf{q}^\vee}$-intertwining maps, and the diagram commutes. Descending this $\mathscr{H}$-level diagram to the $M$-level diagram by decomposing the representations according to the parentheses, we obtain the corresponding commutative diagram of unitary maps
\begin{align*}
  \begin{array}{c}
M_{64}^n {\otimes} M_{53}^6 {\otimes} M_{12}^5 ~ \stackrel{\mathbf{T}_{23}}{\longrightarrow} ~ M_{64}^n {\otimes} M_{23}^9 {\otimes} M_{19}^6 ~ \stackrel{\mathbf{T}_{13}}{\longrightarrow} ~ M_{94}^8 {\otimes} M_{23}^9 {\otimes} M_{18}^n \\
\mathbf{T}_{12} \, \searrow \hspace{50mm} \swarrow \, \mathbf{T}_{12} \\
M_{34}^7 {\otimes} M_{57}^n {\otimes} M_{12}^5 ~ \stackrel{\mathbf{T}_{23}}{\longrightarrow} ~ M_{34}^7 {\otimes} M_{27}^8 {\otimes} M_{18}^n,
\end{array}
\end{align*}
yielding:
\begin{proposition}[the pentagon relation for $\mathbf{T}$]
\label{prop:TT_TTT}
One has
$$
\mathbf{T}_{23} \mathbf{T}_{12} = \mathbf{T}_{12} \mathbf{T}_{13} \mathbf{T}_{23}. \quad \qed
$$
\end{proposition}

Now let us show Prop.\ref{prop:construction_of_T}, which we postponed:

\vs

{\it Proof of Prop.\ref{prop:construction_of_T}.} We shall show
\begin{align}
\label{eq:F_pentagon}
\mathbf{F}_{13} \mathbf{F}_{21} = \mathbf{F}_{21} \mathbf{F}_{23} \mathbf{F}_{13}.
\end{align}
Then we would have $\mathbf{F}_{23} \, \mathbf{F}_{13} \, \mathbf{F}_{21}^{-1} \, \mathbf{F}_{13}^{-1} = \mathbf{F}_{21}^{-1}$, so that Prop.\ref{prop:construction_of_T} indeed holds with
\begin{align}
\label{eq:T_as_F}
\mathbf{T}_{12} = \mathbf{F}_{21}^{-1}.
\end{align}
Recall 
$$
{\bf F}_{ij} = {\bf S}_{ij} (\Psi^\hbar(Q_i+P_j-Q_j))^{-1},
$$
where $\mathbf{S}_{ij} = e^{\frac{-1}{2\pi\hbar}(P_i Q_j - P_i^* Q_j^*)}$. In addition to eq.\eqref{eq:S_conjugation_on_P2Q2}, we also need the conjugation action of ${\bf S}_{12}$ on the remaining $P_1$ and $Q_1$. We write them as
\begin{align}
\label{eq:S_conjugation_on_PQs}
{\bf S}_{ij}^{-1} \, P_j \, {\bf S}_{ij} = P_i + P_j, \quad
{\bf S}_{ij}^{-1} \, Q_j \, {\bf S}_{ij} = Q_j, \quad
{\bf S}_{ij}^{-1} \, P_i \, {\bf S}_{ij} = P_i, \quad
{\bf S}_{ij}^{-1} \, Q_i \, {\bf S}_{ij} = Q_i-Q_j,
\end{align}
The first two are from eq.\eqref{eq:S_conjugation_on_P2Q2}, and the latter two can be seen heuristically from $[P_i, \frac{-1}{2\pi \hbar}(P_i Q_j - P_i^* Q_j^*)] = 0$ and $[Q_i, \frac{-1}{2\pi \hbar} (P_i Q_j - P_i^* Q_j^*)]=-Q_j$, and can be proved rigorously by using Lemmas \ref{lem:S_12_as_saso} and \ref{lem:conjugation_action_of_S_c}. It is obvious that ${\bf S}_{ij}^{-1} P_k {\bf S}_{ij} = P_k$ and ${\bf S}_{ij}^{-1} Q_k {\bf S}_{ij} = Q_k$ if $i,j,k$ are all distinct. Note first that
\begin{align*}
{\bf F}_{12} \, {\bf F}_{21} & = {\bf S}_{13} \ul{ (\Psi^\hbar(Q_1+P_3-Q_3))^{-1} \, {\bf S}_{21} } (\Psi^\hbar(Q_2+P_1-Q_1))^{-1} \\
& = {\bf S}_{13} {\bf S}_{21} (\Psi^\hbar(Q_1 + P_3-Q_3))^{-1} (\Psi^\hbar(Q_2 + P_1 - Q_1))^{-1} \quad \mbox{($\because$ eq.\eqref{eq:S_conjugation_on_PQs})}.
\end{align*}
On the other hand, observe
\begin{align*}
& {\bf F}_{21} {\bf F}_{23} {\bf F}_{13} \\
& = {\bf S}_{21} \ul{ (\Psi^\hbar(Q_2 + P_1 - Q_1))^{-1} {\bf S}_{23}  }(\Psi^\hbar(Q_2+P_3-Q_3))^{-1} {\bf S}_{13} (\Psi^\hbar(Q_1 + P_3 - Q_3))^{-1} \\
& \stackrel{{\rm eq.}\eqref{eq:S_conjugation_on_PQs}}{=} {\bf S}_{21} {\bf S}_{23} \ul{ (\Psi^\hbar( Q_2 - Q_3 + P_1 - Q_1))^{-1} (\Psi^\hbar(Q_2+P_3-Q_3))^{-1} {\bf S}_{13} } (\Psi^\hbar(Q_1 + P_3 - Q_3))^{-1} \\
& \stackrel{{\rm eq.}\eqref{eq:S_conjugation_on_PQs}}{=} {\bf S}_{21} {\bf S}_{23} {\bf S}_{13} (\Psi^\hbar(Q_2-\cancel{Q_3}+P_1-Q_1+\cancel{Q_3}))^{-1} (\Psi^\hbar(Q_2 + P_1+P_3-Q_3))^{-1} (\Psi^\hbar(Q_1+P_3-Q_3))^{-1}.
\end{align*}
First, we claim that
\begin{align}
\label{eq:S_pentagon}
\mathbf{S}_{13} \mathbf{S}_{21} = \mathbf{S}_{21} \mathbf{S}_{23} \mathbf{S}_{13}
\end{align}
holds. This could boil down via Prop.\ref{prop:saso} and Lem.\ref{lem:S_12_as_saso} to showing some equality of $6\times 6$ matrices. Equivalently, one can also show this directly as an equality of operators on
$$
L^2(\mathbb{R}^2, dt_1\,ds_1) {\otimes} L^2(\mathbb{R}^2, dt_2\,ds_2) {\otimes} L^2(\mathbb{R}^2, dt_3\,ds_3) \cong L^2(\mathbb{R}^6, \, dt_1\,ds_1\,dt_2\,ds_2\,dt_3\,ds_3),
$$
using the formula in eq.\eqref{eq:quadratic_part_explicit} for $\mathbf{S}_{ij}$: for any $f\in L^2(\mathbb{R}^6, \, dt_1\,ds_1\,dt_2\,ds_2\,dt_3\,ds_3)$ note
\begin{align*}
( \mathbf{S}_{13} (\mathbf{S}_{21} f) )(t_1,s_1,t_2,s_2,t_3,s_3) & = (\mathbf{S}_{21} f)(t_1,s_1+s_3,t_2,s_2,t_3-t_1,s_3) \\
& = f(t_1-t_2, s_1+s_3, t_2, s_2+(s_1+s_3), t_3-t_1,s_3), \\
(\mathbf{S}_{21}(\mathbf{S}_{23}(\mathbf{S}_{13}f)))(t_1,s_1,t_2,s_2,t_3,s_3) & = (\mathbf{S}_{23}(\mathbf{S}_{13} f))(t_1-t_2,s_1,t_2,s_2+s_1,t_3,s_3) \\
& = (\mathbf{S}_{13} f)(t_1-t_2,s_1,t_2,s_2+s_1+s_3,t_3-t_2,s_3) \\
& = f(t_1-t_2,s_1+s_3, t_2,s_2+s_1+s_3,t_3-t_2-(t_1-t_2),s_3),
\end{align*}
hence indeed $\mathbf{S}_{13} \mathbf{S}_{21}f  = \mathbf{S}_{21} \mathbf{S}_{23} \mathbf{S}_{13} f$.

\vs

It then remains to show the equality
\begin{align*}
& (\Psi^\hbar(Q_1 + P_3-Q_3))^{-1} (\Psi^\hbar(Q_2 + P_1 - Q_1))^{-1} \\
& = (\Psi^\hbar(Q_2+P_1-Q_1))^{-1} (\Psi^\hbar(Q_2 + P_1+P_3-Q_3))^{-1} (\Psi^\hbar(Q_1+P_3-Q_3))^{-1},
\end{align*}
which follows from Prop.\ref{prop:pentagon_identity_for_Psi}, the pentagon identity for $\Psi^\hbar$. Indeed, it is straightforward to check whether $Q_2 + P_1 - Q_1 = (q_{s_2} + q_{t_1} - q_{s_1}) + {\rm i} (-p_{t_2} + p_{s_1} + p_{t_1}) = {\bf x} + {\rm i} \hbar {\bf y}$ and $Q_1+P_3-Q_3 = (q_{s_1}+q_{t_3}-q_{s_3}) + {\rm i} (-p_{t_1} + p_{s_3} + p_{t_3}) = {\bf x}' + {\rm i} \hbar {\bf y}'$ satisfy the conditions of Prop.\ref{prop:pentagon_identity_for_Psi}. \qed

\vs

Now, some words on the uniqueness of $\mathbf{T}$. Notice that the construction of $\mathbf{T}$ in Prop.\ref{prop:construction_of_T} depends on the choice of the unitary intertwiner $\mathbf{F} : \mathscr{H} {\otimes}\mathscr{H} \to M {\otimes} \mathscr{H}$; we used our operator $\mathbf{F}$ constructed in \S\ref{sec:F}. However, $\mathbf{F}$ is not a unique unitary intertwiner from $\mathscr{H}{\otimes}\mathscr{H}$ to $M{\otimes} \mathscr{H}$. Suppose $\mathbf{F}' : \mathscr{H}{\otimes}\mathscr{H} \to M {\otimes} \mathscr{H}$ is \emph{any} unitary intertwining map. Then $\mathbf{F}' \circ \mathbf{F}^{-1} : M {\otimes} \mathscr{H} \to M {\otimes} \mathscr{H}$ is an unitary intertwining map. That is, $\mathbf{F}' \, \mathbf{F}^{-1}$ acts trivially by conjugation on $1\otimes \pi(u)$ for each $u\in \mathcal{C}_{\mathbf{q},\mathbf{q}^\vee}$, which are operators on $L^2(\mathbb{R}^2, \,dt_1\,ds_1) {\otimes} L^2(\mathbb{R}^2,\,dt_2\,ds_2)$. Using the irreducibility of the representation $(\mathscr{H},\pi)$ of $\mathcal{C}_{\mathbf{q},\mathbf{q}^\vee}$ of Conjecture \ref{conj:strong_irreducibility}, one can deduce that $\mathbf{F}' \, \mathbf{F}^{-1}$ must act as a scalar times identity on the second tensor factor $L^2(\mathbb{R}^2,\,dt_2\,ds_2) \equiv \mathscr{H}$ of the tensor product $M{\otimes} \mathscr{H}$. Thus there must be some unitary operator $U$ on $L^2(\mathbb{R}^2) \equiv M$ such that $\mathbf{F}' \, \mathbf{F}^{-1} = U\otimes 1$, hence
\begin{align}
\label{eq:bf_F_prime}
\mathbf{F}' = (U\otimes 1) \circ \mathbf{F},
\end{align}
compactly written as $\mathbf{F}'_{12} = U_1 \, \mathbf{F}_{12}$. Conversely, for \emph{any} unitary operator $U$ on $L^2(\mathbb{R}^2)$ the operator $(U\otimes 1) \circ \mathbf{F}$ is a unitary intertwiner from $\mathscr{H}{\otimes}\mathscr{H}$ to $M{\otimes}\mathscr{H}$. Note now that
\begin{align*}
& \mathbf{F}_{13}' \, \mathbf{F}_{21}' \, (\mathbf{F}_{13}')^{-1} \, (\mathbf{F}_{23}')^{-1}
= U_1 \, \mathbf{F}_{13} \, U_2 \, \mathbf{F}_{21} \, \mathbf{F}_{13}^{-1} \, U_1^{-1} \, \mathbf{F}_{23}^{-1} \, U_2^{-1} \\
& \stackrel{\vee}{=}
(U_1U_2) \, \mathbf{F}_{13} \, \mathbf{F}_{21} \, \mathbf{F}_{13}^{-1} \, \mathbf{F}_{23}^{-1} \, (U_1U_2)^{-1}
= (U_1U_2) \, \mathbf{T}_{12} \, (U_1U_2)^{-1}
= ( (U\otimes U) \mathbf{T} (U\otimes U)^{-1} ) \otimes 1,
\end{align*}
where, for the checked equality we used the fact that two expressions with the subscripts not intersecting commute with each other (because they are acting on mutually distinct tensor factors). In particular, the dotted arrow map $\mathbf{F}_{13}' \, \mathbf{F}_{21}' \, (\mathbf{F}_{13}')^{-1} \, (\mathbf{F}_{23}')^{-1} : M {\otimes} M {\otimes} \mathscr{H} \to M{\otimes} M {\otimes} \mathscr{H}$ constructed from the intertwiner $\mathbf{F}'$ is also of the form $\mathbf{T}' \otimes 1$, i.e. acting only on the first two tensor factors, where $\mathbf{T}' = (U\otimes U) \, \mathbf{T} \, (U\otimes U)^{-1}$. This establishes Prop.\ref{prop:construction_of_T} not just for our $\mathbf{F}$, but when $\mathbf{F}$ is any unitary intertwiner from $\mathscr{H}{\otimes}\mathscr{H}$ to $M{\otimes}\mathscr{H}$, and moreover also the uniqueness of the resulting $\mathbf{T}$ up to conjugation by a unitary operator on $L^2(\mathbb{R}^2){\otimes}L^2(\mathbb{R}^2)$ of the form $U\otimes U$. Among the non-unique choices for the intertwining map $\mathbf{F} : \mathscr{H} {\otimes} \mathscr{H} \to M {\otimes} \mathscr{H}$, our choice is special with respect to the simple relationship of it as in eq.\eqref{eq:T_as_F} with the corresponding $\mathbf{T}$; there might be some other choices that are more natural with respect to other aspects, analogously as discussed in \cite[Remark 4.7]{FKi}.

\section{Dual representation}
\label{sec:dual_representation}

\begin{definition}[the Hom space]
\label{def:Hom_space}
Let $V$ and $W$ be complex Hilbert spaces. Denote by $|| \cdot ||_V$ and $||\cdot ||_W$ their Hilbert space norms. The space of bounded (or continuous) linear operators $V\to W$ is denoted by
\begin{align*}
\mathcal{B}(V,W) := \left\{ \, \varphi : V\to W \, \left| \, \mbox{$\varphi$ is $\mathbb{C}$-linear and $||\varphi|| = \sup_{v\in V\setminus\{0\}} \frac{||\varphi (v)||_W}{||v||_V}< \infty$} \, \right. \right\}.
\end{align*}

\vs

Suppose that $V$ and $W$ are separable. Let $\{v_i\}_{i\in \mathbb{N}}$ be an orthonormal basis of $V$. A bounded linear operator $\varphi : V \to W$ (i.e. $\varphi \in \mathcal{B}(V,W)$) is said to be \ul{\em Hilbert-Schmidt} if
$$
\sum_{i \in \mathbb{N}} || \varphi(v_i) ||_W^2 < \infty.
$$
Denote the space of all Hilbert-Schmidt operators $V\to W$ as
$$
{\rm Hom}_{\rm HS} (V,W) := \{ \, \varphi \in \mathcal{B}(V,W) \, | \, \mbox{$\varphi$ is Hilbert-Schmidt} \, \}.
$$
\end{definition}
It is known that the definition of Hilbert-Schmidt does not depend on the choice of an orthonormal basis of the domain Hilbert space $V$, and also that the following inner product on ${\rm Hom}_{\rm HS}(V,W)$ makes it a separable Hilbert space:
$$
\langle \varphi_1, \varphi_2 \rangle_{\rm HS} := \sum_{i\in \mathbb{N}} \langle \varphi_1(v_i), \varphi_2(v_i)\rangle_W, \qquad \varphi_1, \varphi_2 \in {\rm Hom}_{\rm HS}(V,W),
$$
where $\langle \cdot,\cdot\rangle_W$ is the inner product on $W$. See \cite[Appendix B]{PR} (also \cite[\S VI.6]{RS}) for more details and proofs about Hilbert-Schmidt operators. Although we find ${\rm Hom}_{\rm HS}(\cdot,\cdot)$ to be providing the best basic model for the space of maps between representations, the specific choice of a model for the ${\rm Hom}$ space will not matter too seriously for the main purpose of the present paper, as shall be seen later.

\vs

For the dual of a Hilbert space $V$, we take the continuous linear dual $\mathcal{B}(V,\mathbb{C})$ as usual, where $\mathbb{C}$ here is given the standard Hilbert space structure. The following mimics the standard Hopf algebra action on the dual spaces of finite dimensional representations; see e.g. \cite{CP}.

\begin{definition}[the left and the right dual representations]
Let $V$ be a densely-defined representation of a Hopf algebra $\mathcal{C}$ (see Def.\ref{def:densely-defined_representation}), with the antipode $S: \mathcal{C} \to \mathcal{C}$. Define the \ul{\em left dual representation} $V'$ of $V$  as $\mathcal{B}(V,\mathbb{C})$ as a Hilbert space, on which the densely-defined action of $\mathcal{C}$ is given by
$$
(u.\varphi)(v) := \varphi(S(u).v) \qquad\mbox{for $u\in\mathcal{C}$, $\varphi \in \mathcal{B}(V,\mathbb{C})$, $v\in V$}.
$$
Define the \ul{\em right dual representation} ${}'V$ of $V$ as $\mathcal{B}(V,\mathbb{C})$ as a Hilbert space, on which the densely-defined action of $\mathcal{C}$ is given by
$$
(u.\varphi)(v) := \varphi(S^{-1}(u).v) \qquad\mbox{for $u\in\mathcal{C}$, $\varphi \in \mathcal{B}(V,\mathbb{C})$, $v\in V$}.
$$
\end{definition}
One can easily show that the above two indeed define algebra actions, using the fact that $S$ and $S^{-1}$ are algebra anti-homomorphisms. Moreover, one finds that the natural pairing maps
$$
V' \otimes V \to \mathbb{C} \qquad\mbox{and}\qquad
V\otimes {}' V \to \mathbb{C}
$$
are intertwining, where the actions of $\mathcal{C}$ on the left hand side is defined using the coproduct while the $\mathbb{C}$ in the right hand side stands for the one-dimensional trivial representation.

\vs

For a Hilbert space $V$, notice that the vector space $\mathcal{B}(V,\mathbb{C})$ can naturally be realized as $V$ itself, thanks to the Riesz Representation Theorem; more explicitly, if we denote by $\langle \cdot,\cdot \rangle_V$ the inner product of $V$, for each $\varphi \in \mathcal{B}(V,\mathbb{C})$ there exists a unique $v_\varphi \in V$ such that $\langle w, v_\varphi \rangle_V = \varphi(w)$ for all $w\in V$. This map $\mathcal{B}(V,\mathbb{C}) \to V : \varphi \mapsto v_\varphi$ is complex-conjugate-linear and bijective, and moreover is an isometry with respect to the operator norm on $\mathcal{B}(V,\mathbb{C})$ and the usual norm on $V$. 

\vs

In the case when $V$ is the representation $\mathscr{H} = L^2(\mathbb{R}^2, \, dt \, ds)$ of $\mathcal{C} = \mathcal{C}_{\mathbf{q},\mathbf{q}^\vee}$ studied in \S\ref{sec:integrable_representation_of_C_q}, one can tweak this map $\mathcal{B}(V,\mathbb{C}) \to V$ a little bit to turn it into a Hilbert space isomorphism, with the help of the complex-conjugation map ${\,}^- : \mathscr{H} \to \mathscr{H} : f\mapsto \ol{f}$, which is an involutive complex-conjugate-linear bijection. Namely, we realize the vector space $\mathcal{B}(\mathscr{H},\mathbb{C})$ as
$$
\mathcal{B}(\mathscr{H},\mathbb{C}) \equiv L^2(\mathbb{R}^2, \, dt \,ds).
$$
Consider the following pairing map
$$
\langle\langle \, , \, \rangle\rangle ~ : ~
\mathcal{B}(\mathscr{H},\mathbb{C}) \times \mathscr{H} \to \mathbb{C}
$$
defined as
$$
\langle\langle f,g\rangle\rangle := \int_{\mathbb{R}^2} f(t,s) \, g(t,s) \, dt \, ds, \qquad \forall f,g\in L^2(\mathbb{R}^2,\, dt\,ds).
$$
We shall now write the action $\pi'$ on the left dual $\mathscr{H}' \equiv L^2(\mathbb{R}^2, \, dt\,ds)$ and the action ${}'\pi$ on the right dual ${}'\mathscr{H} \equiv L^2(\mathbb{R}^2, \, dt\,ds)$ more explicitly as follows. Observe for $f,g\in L^2(\mathbb{R}^2)$ and $u\in \mathcal{C}_{\mathbf{q},\mathbf{q}^\vee}$ that
\begin{align*}
  \langle\langle \pi'(u).f, g\rangle\rangle = \langle\langle f, \pi(S(u)).g\rangle\rangle = \langle\langle (\pi(S(u)))^\mathrm{tr}.f, \, g\rangle\rangle,
\end{align*}
where $\mathrm{tr}$ stands for the `transpose'; likewise for ${}'\pi$ with $S$ replaced by $S^{-1}$, so
$$
\pi'(u) = (\pi(S(u)))^\mathrm{tr}, \qquad
{}'\pi(u) = (\pi(S^{-1}(u)))^\mathrm{tr}.
$$
As $\langle\langle f,g\rangle\rangle = \langle f,\ol{g}\rangle_{L^2(\mathbb{R}^2)}$, one can see that taking a transpose of an operator $A$ is same as taking the adjoint of the complex-conjugate of the operator, i.e.
$$
A^\mathrm{tr} = ({}^- \circ A\circ {}^-)^*.
$$
At the moment, the symbol $*$ is not to stand for the full-fledged operator adjoint, but only the restriction to some dense subspace; let us consider only the case when $f,g$ belong to the space $\mathscr{D}$ of eq.\eqref{eq:D}. As the adjoint $*$ changes the product order of operators, i.e. $(AB)^* = B^* A^*$, so does $\mathrm{tr}$, i.e. $(AB)^\mathrm{tr} = B^\mathrm{tr} \, A^\mathrm{tr}$; it is easy to see that $A \mapsto A^\mathrm{tr}$ is complex linear, as well as involutive. One finds
$$
q_t^\mathrm{tr} = ({}^-\circ q_t\circ {}^-)^* = q_t^* = q_t, \qquad p_t^\mathrm{tr} = ({}^- \circ p_t \circ {}^-)^* = - p_t^* = - p_t,
$$
and likewise $q_s^\mathrm{tr}=q_s$, $p_s^\mathrm{tr}=-p_s$. Thus
\begin{align*}
& P^\mathrm{tr} = (q_t + {\rm i}\,p_s)^\mathrm{tr} = q_t - {\rm i}\,p_s = P^*, \quad Q^\mathrm{tr}=(q_s - {\rm i}\,p_t)^\mathrm{tr} = q_s + {\rm i} \, p_t = Q^*, \\
& (P^*)^\mathrm{tr} = (q_t - {\rm i}\,p_s)^\mathrm{tr} = q_t + {\rm i}\,p_s = P, \quad (Q^*)^\mathrm{tr}=(q_s + {\rm i}\,p_t)^\mathrm{tr} = q_s - {\rm i} \, p_t = Q, \\
\end{align*}
so
\begin{align*}
\pi'(Z_1) & = (\pi(S(Z_1)))^\mathrm{tr} = (\pi(Z_1^{-1}))^\mathrm{tr} = (\mathbf{Z}_1^{-1})^\mathrm{tr} = (e^{-P})^\mathrm{tr} = e^{-P^*}, \\
\pi'(Z_2) & = (\pi(S(Z_2))^\mathrm{tr} = (\pi(-Z_2 Z_1^{-1}))^\mathrm{tr} = (-\mathbf{Z}_2 \, \mathbf{Z}_1^{-1})^\mathrm{tr} = (-e^Q \, e^{-P})^\mathrm{tr} = - e^{-P^*} \, e^{Q^*}, \\
\pi'(Z_1^*) & = (\pi(S(Z_1^*)))^\mathrm{tr} = (\pi((Z_1^*)^{-1}))^\mathrm{tr} = ((\mathbf{Z}_1^*)^{-1})^\mathrm{tr} = (e^{-P^*})^\mathrm{tr} = e^{-P}, \\
\pi'(Z_2^*) & = (\pi(S(Z_2^*))^\mathrm{tr} = (\pi(-Z_2^* (Z_1^*)^{-1}))^\mathrm{tr} = (-\mathbf{Z}_2^* \, (\mathbf{Z}_1^*)^{-1})^\mathrm{tr} = (-e^{Q^*} \, e^{-P^*})^\mathrm{tr} = - e^{-P} \, e^{Q}, \\
{}'\pi(Z_1) & = (\pi(S^{-1}(Z_1)))^\mathrm{tr} = (\pi(Z_1^{-1}))^\mathrm{tr} = e^{-P^*}, \\
{}'\pi(Z_2) & = (\pi(S^{-1}(Z_2)))^\mathrm{tr} = (\pi(-Z_1^{-1}Z_2))^\mathrm{tr} = (-\mathbf{Z}_1^{-1} \mathbf{Z}_2)^\mathrm{tr} = (- e^{-P} e^Q)^\mathrm{tr} = - e^{Q^*} e^{-P^*}, \\
{}'\pi(Z_1^*) & = (\pi(S^{-1}(Z_1^*)))^\mathrm{tr} = (\pi((Z_1^*)^{-1}))^\mathrm{tr} = e^{-P}, \\
{}'\pi(Z_2^*) & = (\pi(S^{-1}(Z_2^*)))^\mathrm{tr} = (\pi(-(Z_1^*)^{-1}Z_2^*))^\mathrm{tr} = (-(\mathbf{Z}_1^*)^{-1} \mathbf{Z}_2^*)^\mathrm{tr} = (- e^{-P^*} e^{Q^*})^\mathrm{tr} = - e^{Q} e^{-P}.
\end{align*}
Likewise, one obtains
\begin{align*}
& \pi'(Z_1^\vee) = e^{-P^*/({\rm i} \hbar)}, \qquad
\pi'(Z_2^\vee) = - e^{-P^*/({\rm i} \hbar)} \, e^{Q^*/({\rm i} \hbar)} \\
& \pi'((Z_1^\vee)^*) = e^{P/({\rm i} \hbar)}, \qquad
\pi'((Z_2^\vee)^*) = - e^{P/({\rm i} \hbar)} \, e^{-Q/({\rm i} \hbar)} \\
& {}'\pi(Z_1^\vee) = e^{-P^*/({\rm i} \hbar)}, \qquad
{}'\pi(Z_2^\vee) = -  e^{Q^*/({\rm i} \hbar)}\, e^{-P^*/({\rm i} \hbar)}, \\
& {}'\pi((Z_1^\vee)^*) = e^{P/({\rm i} \hbar)}, \qquad
{}'\pi((Z_2^\vee)^*) = -  e^{-Q/({\rm i} \hbar)}\, e^{P/({\rm i} \hbar)}.
\end{align*}

\vs

Let us now build non-unitary $\mathcal{C}_{\mathbf{q},\mathbf{q}^\vee}$-intertwiners from $\mathscr{H}$ to its left and right duals. We make use of the special linear operators ${\bf S}_{\bf c}$ developed in \S\ref{sec:F}. 
\begin{proposition}[weak isomorphisms between $\mathscr{H}$ and its duals]
\label{prop:weak_isomorphisms_with_duals}
The densely-defined operators
$$
C := \mathbf{S}_{\smallmattwo{-1}{-1}{0}{1}} \, \circ \, e^{q_t + \hbar^{-1} p_s} ~ : ~ L^2(\mathbb{R}^2, \, dt\,ds) \to L^2(\mathbb{R}^2,\, dt\,ds)
$$
and
$$
D := \mathbf{S}_{\smallmattwo{-1}{-1}{0}{1}} \, \circ \, e^{-q_t + \hbar^{-1} p_s} ~ : ~ L^2(\mathbb{R}^2, \, dt\,ds) \to L^2(\mathbb{R}^2,\, dt\,ds)
$$
provide non-unitary densely-defined bijective (on $\mathscr{D}$) maps $C : \mathscr{H} \to \mathscr{H}'$ and $D: {}' \mathscr{H} \to \mathscr{H}$ that intertwine the densely-defined actions $\pi,\pi',{}'\pi$ of the algebra $\mathcal{C}_{\mathbf{q},\mathbf{q}^\vee}$.
\end{proposition}

{\it Proof.} From eq.\eqref{eq:technical_conjugation} and eq.\eqref{eq:p_j_and_q_k_Heisenberg_relations} one observes
\begin{align*}
& e^{q_t + \hbar^{-1} p_s} \, p_t \, e^{- (q_t + \hbar^{-1} p_s)} = p_t - \pi {\rm i} \, \hbar, \qquad
e^{q_t + \hbar^{-1} p_s} \, p_s \, e^{- (q_t + \hbar^{-1} p_s)} = p_s, \\
& e^{q_t + \hbar^{-1} p_s} \, q_t \, e^{- (q_t + \hbar^{-1} p_s)} = q_t, \qquad\qquad\qquad
e^{q_t + \, \hbar^{-1} p_s} \, q_s \, e^{- (q_t + \, \hbar^{-1} p_s)} = q_s + \pi {\rm i}.
\end{align*}
These equalities make good senses as operators $\mathscr{D} \to \mathscr{D}$, and can be directly checked (not just heuristically using eq.\eqref{eq:technical_conjugation}). Let us write $\mathbf{c} = \smallmattwo{c_{tt}}{c_{ts}}{c_{st}}{c_{ss}} = \smallmattwo{-1}{-1}{0}{1}$. Then $\mathbf{c}^{-1} = \smallmattwo{c^{tt}}{c^{ts}}{c^{st}}{c^{ss}} = \smallmattwo{-1}{-1}{0}{1}$. From eq.\eqref{eq:saso_conjugation_on_p}--\eqref{eq:saso_conjugation_on_q} and eq.\eqref{eq:p_and_q_j} we have
\begin{align*}
\mathbf{S}_\mathbf{c} \, p_t \, \mathbf{S}_\mathbf{c}^{-1} = -p_t -p_s, \qquad
\mathbf{S}_\mathbf{c} \, p_s \, \mathbf{S}_\mathbf{c}^{-1} = p_s, \qquad
\mathbf{S}_\mathbf{c} \, q_t \, \mathbf{S}_\mathbf{c}^{-1} = -q_t, \qquad
\mathbf{S}_\mathbf{c} \, q_s \, \mathbf{S}_\mathbf{c}^{-1} = -q_t +q_s.
\end{align*}
Hence $C := \mathbf{S}_\mathbf{c} \circ e^{q_t + \hbar^{-1} p_s}$ satisfies the conjugation equations
\begin{align*}
C \, P \, C^{-1} = 
C \, (q_t+{\rm i}\,p_s) \, C^{-1}
& = - q_t + {\rm i}\,p_s = - P^*, \\
C \, Q \, C^{-1} =
C \, (q_s - {\rm i}\,p_t) \, C^{-1}
& = \pi({\rm i}-\hbar) - q_t + q_s + {\rm i}\, (p_t+p_s) \\
& = \pi({\rm i}-\hbar) - P^* + Q^*, \\
C \, P^* \, C^{-1} = 
C \, (q_t-{\rm i}\,p_s) \, C^{-1}
& = - q_t - {\rm i}\,p_s = - P, \\
C \, Q^* \, C^{-1} =
C \, (q_s + {\rm i}\,p_t) \, C^{-1}
& = \pi({\rm i}+\hbar) - q_t + q_s + {\rm i}\, (-p_t-p_s) \\
& = \pi({\rm i}+\hbar) - P + Q,
\end{align*}
as equalities of operators $\mathscr{D} \to \mathscr{D}$. Exponentiating, one obtains
\begin{align*}
& C \, \pi(Z_1) \, C^{-1}
= C \, e^P \, C^{-1} = e^{-P^*} = \pi'(Z_1), \\
& C \, \pi(Z_2) \, C^{-1}
= C \, e^Q \, C^{-1}
= - e^{-\pi\hbar} e^{-P^*+Q^*}
\stackrel{\vee}{=} - e^{-P^*} e^{Q^*} = \pi'(Z_2), \\
& C \, \pi(Z_1^*) \, C^{-1}
= C \, e^{P^*} \, C^{-1} = e^{-P} = \pi'(Z_1^*), \\
& C \, \pi(Z_2^*) \, C^{-1}
= C \, e^{Q^*} \, C^{-1}
= - e^{\pi\hbar} e^{-P+Q}
\stackrel{\vee}{=} - e^{-P} e^{Q} = \pi'(Z_2^*),
\end{align*}
where the checked equality is from the BCH formula in eq.\eqref{eq:BCH} together with $[-P,Q] = 2\pi\hbar\cdot\mathrm{id}$, $[-P^*,Q^*] = -2\pi\hbar\cdot\mathrm{id}$ (see eq.\eqref{eq:PQ_commutator}, \eqref{eq:other_commutators}). Indeed, one can also directly check that these equalities hold as operators $\mathscr{D} \to \mathscr{D}$, not just heuristically by the BCH formula. Likewise,
\begin{align*}
& C \, \pi(Z_1^\vee) \, C^{-1}
= C \, e^{P/({\rm i} \hbar)} \, C^{-1} = e^{-P^*/({\rm i} \hbar)} = \pi'(Z_1^\vee), \\
& C \, \pi(Z_2^\vee) \, C^{-1}
= C \, e^{Q/({\rm i} \hbar)} \, C^{-1} = - e^{\pi/\hbar} e^{(-P^*+Q^*)/({\rm i} \hbar)} \\
& C \, \pi((Z_1^\vee)^*) \, C^{-1}
= C \, e^{- P^*/({\rm i} \hbar)} \, C^{-1} = e^{P/({\rm i} \hbar)} = \pi'( (Z_1^\vee)^* ), \\
&\qquad\qquad \qquad\qquad\qquad\qquad\qquad\qquad \,\,\,
= - e^{-P^*/({\rm i} \hbar)} e^{Q^*/({\rm i} \hbar)} = \pi'(Z_2^\vee), \\
& C \, \pi( (Z_2^\vee)^* ) \, C^{-1}
= C \, e^{-Q^*/({\rm i} \hbar)} \, C^{-1} = - e^{-\pi/\hbar} e^{(P-Q)/({\rm i} \hbar)} \\
&\qquad\qquad \qquad\qquad\qquad\qquad\qquad\qquad \,\,\,
= - e^{P/({\rm i} \hbar)} e^{-Q/({\rm i} \hbar)} = \pi'((Z_2^\vee)^*).
\end{align*}
Thus indeed $C = \mathbf{S}_\mathbf{c} \circ e^{q_t + \hbar^{-1} p_s} : \mathscr{H} \to \mathscr{H}'$ is an intertwiner.

\vs

For the other assertion about $D = \mathbf{S}_\mathbf{c} \circ e^{-q_t + \hbar^{-1} p_s}$, we proceed similarly. From eq.\eqref{eq:technical_conjugation} and eq.\eqref{eq:p_j_and_q_k_Heisenberg_relations},
\begin{align*}
& e^{-q_t + \hbar^{-1} p_s} \, p_t \, e^{-(-q_t + \hbar^{-1} p_s)} = p_t + \pi {\rm i} \, \hbar, \qquad
e^{-q_t + \hbar^{-1} p_s} \, p_s \, e^{-(-q_t + \hbar^{-1} p_s)} = p_s, \\
& e^{-q_t + \hbar^{-1} p_s} \, q_t \, e^{-(-q_t + \hbar^{-1} p_s)} = q_t, \qquad\qquad\qquad
e^{-q_t + \, \hbar^{-1} p_s} \, q_s \, e^{-(-q_t + \, \hbar^{-1} p_s)} = q_s + \pi {\rm i}.
\end{align*}
Hence $D$ satisfies the same conjugation equations as $C$ with $\hbar$ in the right hand side replaced by $-\hbar$; from those one can deduce the following equations
\begin{align*}
& D \, (-P^*) \, D^{-1} = P, \qquad
D \, (\pi(-{\rm i}+\hbar) +  Q^*- P^*) \, D^{-1} = Q, \\
& D \, (-P) \, D^{-1} = P^*, \qquad
D \, (\pi(-{\rm i}-\hbar) + Q - P ) \, D^{-1} = Q^*,
\end{align*}
thus
\begin{align*}
& D \, ({}'\pi(Z_1)) \, D^{-1}
= D \, e^{-P^*} \, D^{-1}
= e^P = \pi(Z_1) \\
& D \, ({}'\pi(Z_2)) \, D^{-1}
= D \, (- e^{Q^*}e^{-P^*}) \, D^{-1}
\stackrel{\vee}{=} D \, e^{\pi(-{\rm i}+\hbar) + Q^*-P^*} \, D^{-1} = e^Q = \pi(Z_2), \\
& D \, ({}'\pi(Z_1^*)) \, D^{-1}
= D \, e^{-P} \, D^{-1}
= e^{P^*} = \pi(Z_1^*) \\
& D \, ({}'\pi(Z_2^*)) \, D^{-1}
= D \, (- e^Qe^{-P}) \, D^{-1}
\stackrel{\vee}{=} D\, e^{\pi(-{\rm i}-\hbar) +Q-P} \, D^{-1} = e^{Q^*} = \pi(Z_2^*),
\end{align*}
and
\begin{align*}
& D \, ({}'\pi(Z_1^\vee)) \, D^{-1}
= D \, e^{-P^*/({\rm i} \hbar)} \, D^{-1}
= e^{P/({\rm i} \hbar)} = \pi(Z_1^\vee) \\
& D \, ({}'\pi(Z_2^\vee)) \, D^{-1}
= D \, (- e^{Q^*/({\rm i} \hbar)}e^{-P^*/({\rm i} \hbar)}) \, D^{-1} \\
& \qquad\qquad\qquad\quad 
\stackrel{\vee}{=} D \, e^{( \pi(-{\rm i}+\hbar) + Q^*-P^*)/({\rm i} \hbar)} \, D^{-1} = e^{Q/({\rm i} \hbar)} = \pi(Z_2^\vee), \\
& D \, ({}'\pi((Z_1^\vee)^*)) \, D^{-1}
= D \, e^{P/({\rm i} \hbar)} \, D^{-1}
= e^{-P^*/({\rm i} \hbar)} = \pi((Z_1^\vee)^*) \\
& D \, ({}'\pi((Z_2^\vee)^*)) \, D^{-1}
= D \, (- e^{-Q/({\rm i} \hbar)}e^{P/({\rm i} \hbar)}) \, D^{-1} \\
& \qquad\qquad\qquad\qquad 
\stackrel{\vee}{=} D\, e^{( \pi({\rm i}+\hbar) -Q+P )/({\rm i} \hbar)} \, D^{-1} = e^{-Q^*/({\rm i} \hbar)} = \pi((Z_2^\vee)^*),
\end{align*}
where the checked equalities are again suppored by the BCH formula and can be directly checked as equalities of operators $\mathscr{D} \to \mathscr{D}$. Thus indeed $D : {}'\mathscr{H} \to \mathscr{H}$ is an intertwiner. \qed

\vs

The irreducibility of the representation $(\mathscr{H},\pi)$ of $\mathcal{C}_{\mathbf{q},\mathbf{q}^\vee}$ in Conjecture \ref{conj:strong_irreducibility} would imply that these isomorphisms $C$ and $D$ are unique intertwining densely-defined bijective linear maps, up to multiplicative scalars.

\section{The order three operator $\mathbf{A}$}
\label{sec:A}

Consider
$$
\mathscr{H}_1 \otimes \mathscr{H}_2 \cong M_{12}^3 \otimes \mathscr{H}_3 \qquad\mbox{and}\qquad
\mathscr{H}_2 \otimes \mathscr{H}_3 \cong M_{23}^1 \otimes \mathscr{H}_1
$$
where each $\mathscr{H}_j$ stands for the representation $(\mathscr{H},\pi)$. We shall construct an operator
$$
\mathbf{A} : M_{12}^3 \to M_{23}^1
$$
using the representation theory of Hopf algebras, using analogous methods to \cite{FKi} \cite{Kim_JPAA}.
\begin{definition}[the invariant subspace]
Let $V$ be a densely-defined representation of a Hopf algebra $\mathcal{C}$ in the sense of Def.\ref{def:densely-defined_representation}. The \ul{\em invariant subspace} of $V$, denoted by $\mathrm{Inv}(V)$, is the set of all elements of $V$ on which $\mathcal{C}$ acts trivially:
$$
\mathrm{Inv}(V) := \{ v\in V \, | \, u.v = \epsilon(u)v, ~ \forall u \in \mathcal{C} \},
$$
where $\epsilon : \mathcal{C} \to \mathbb{C}$ is the counit.
\end{definition}
\begin{definition}[the space of intertwiners]
Let $V,W$ be densely-defined representations of a Hopf algebra $\mathcal{C}$ (see Def.\ref{def:densely-defined_representation}). The \ul{\em space of intertwiners} from $V$ to $W$ is defined as the following space
$$
\mathrm{Hom}_\mathcal{C}(V,W) := \{ \varphi\in \mathrm{Hom}_{\rm HS}(V,W) \, | \, \varphi(u.v) = u.(\varphi(v)), ~ \forall u\in \mathcal{C}, ~ \forall v\in V\},
$$
where ${\rm Hom}_{\rm HS}(\cdot,\cdot)$ is as defined in Def.\ref{def:Hom_space}.
\end{definition}
Of course, the quantifier `$\forall v\in V$' must be understood appropriately, because the action of $\mathcal{C}$ is defined only on a dense subspace.

\vs

For separable complex Hilbert spaces $V$ and $W$, one has a natural Hilbert space isomorphism
\begin{align}
\label{eq:Hom_using_dual}
\mathrm{Hom}_{\rm HS} (V,W) \cong \mathcal{B} (V,\mathbb{C}) \otimes W,
\end{align}
where the right hand side means the Hilbert space tensor product; this statement is well known when $V=W$. An element $\sum f \otimes w$ of the right hand side (which is a finite sum or suitable convergent infinite series), for $f \in \mathcal{B}(V,\mathbb{C})$ and $w \in W$, corresponds to an element of $\mathrm{Hom}_{\rm HS}(V,W)$ that sends each $v\in V$ to $\sum f(v)\,w \in W$. The right hand side of eq.\eqref{eq:Hom_using_dual} can of course be written as $W \otimes \mathcal{B}(V,\mathbb{C})$, because the tensor product of Hilbert spaces is commutative.

\vs

When $V$ and $W$ are densely-defined representations of a Hopf algebra $\mathcal{C}$, we know the actions of $\mathcal{C}$ on $W$, on $\mathcal{B}(V,\mathbb{C})$, and on the tensor products of representations, the last being given via the coproduct. Hence one may naturally construct a $\mathcal{C}$-action on the right hand side of eq.\eqref{eq:Hom_using_dual}, which one can carry over to the left hand side. However, as the coproduct is not co-commutative in general, it does matter whether to write the right hand side of eq.\eqref{eq:Hom_using_dual} as is now or as $W\otimes \mathcal{B}(V,\mathbb{C})$. Also, there are two actions on $\mathcal{B}(V,\mathbb{C})$ as we saw, namely the left action and the right action. Among these four possible ways of defining actions on the right hand side of eq.\eqref{eq:Hom_using_dual}, we find the following two to be natural:
$$
{}'V \otimes W  \qquad\mbox{and}\qquad W\otimes V'.
$$
The first one carries over to the following action on the left hand side, namely on $\mathrm{Hom}_{\rm HS}(V,W)$:
$$
(u.\varphi)(v) := \sum u_2.(\varphi(S^{-1}(u_1).v)) \quad\mbox{for $u\in \mathcal{C}$, $\varphi \in \mathrm{Hom}_{\rm HS}(V,W)$, $v\in V$, where $\Delta u = \sum u_1 \otimes u_2$},
$$
which we call the \ul{\em right action}, while the second one to
$$
(u.\varphi)(v) := \sum u_1.(\varphi(S(u_2).v)) \quad\mbox{for $u\in \mathcal{C}$, $\varphi \in \mathrm{Hom}_{\rm HS}(V,W)$, $v\in V$, when $\Delta u = \sum u_1 \otimes u_2$},
$$
which we call the \ul{\em left action}; one may recognize the left action on $\mathrm{Hom}_{\rm HS}(V,W)$ to be an analog of the standard way to make $\mathrm{Hom}_\mathbb{C}(V,W)$ a representation of $\mathcal{C}$ in the case when $V$ and $W$ are finite dimensional representations (see e.g. \cite{CP}).

\begin{definition}[the left and the right representations on the Hom spaces]
Let $V,W$ be densely-defined representations of a Hopf algebra $\mathcal{C}$, in the sense of Def.\ref{def:densely-defined_representation}. Define $\mathrm{Hom}_{\rm HS}^\mathrm{L}(V,W)$ to be the densely-defined representation of $\mathcal{C}$ on the space $\mathrm{Hom}_{\rm HS}(V,W)$ with the left action, and $\mathrm{Hom}_{\rm HS}^\mathrm{R}(V,W)$ to be that on the space $\mathrm{Hom}_{\rm HS}(V,W)$ with the right action.
\end{definition}

\begin{lemma}
One has
$$
\mathrm{Hom}_\mathcal{C}(V,W) ~ \subseteq ~ \mathrm{Inv}(\mathrm{Hom}^\mathrm{L}_{\rm HS}(V,W)) \cap \mathrm{Inv}(\mathrm{Hom}^\mathrm{R}_{\rm HS}(V,W))
$$
\end{lemma}

{\it Proof.} Suppose $\varphi \in \mathrm{Hom}_\mathcal{C}(V,W) \subset \mathrm{Hom}_{\rm HS}(V,W)$. The left action of $u\in \mathcal{C}$ on $\varphi$ is given by $(u.\varphi)(v) = \sum u_1.(\varphi(S(u_2).v)) = \sum u_1.(S(u_2).(\varphi(v))) = (\sum u_1 \, S(u_2)).(\varphi(v))$; one of the Hopf algebra axioms tells us $\sum u_1 \, S(u_2) = \epsilon(u)$, so the left action of $u$ on $\varphi$ is a trivial action, hence $\varphi\in\mathrm{Inv}(\mathrm{Hom}^\mathrm{L}_{\rm HS}(V,W))$. Meanwhile, the right action of $u$ on $\varphi$ is $(u.\varphi)(v) = \sum u_2.(\varphi(S^{-1}(u_1).v)) = (\sum u_2\, S^{-1}(u_1)).(\varphi(v)) = (S^{-1}(\sum u_1\, S(u_2))).(\varphi(v)) = (S^{-1}(\epsilon(u))).(\varphi(v)) = \epsilon(u).\varphi(v) = \epsilon(u)\cdot\varphi(v)$, hence the right action of $u$ on $\varphi$ is a trivial action, so $\varphi \in \mathrm{Inv}(\mathrm{Hom}^\mathrm{R}_{\rm HS}(V,W))$. \qed

\vs

The philosophy that we shall heuristically take from now on is
$$
\mathrm{Hom}_\mathcal{C}(V,W) = \mathrm{Inv}(\mathrm{Hom}^\mathrm{L}_{\rm HS}(V,W)) = \mathrm{Inv}(\mathrm{Hom}^\mathrm{R}_{\rm HS}(V,W)),
$$
which we do not claim or justify; but keep this in mind in the upcoming constructions and computations.

\vs

In order to build the sought-for operator $\mathbf{A}:M_{12}^3 \to M_{23}^1$, we first establish the following isomorphism 
$$
I_{12}^3 : M_{12}^3 \to \mathrm{Hom}_{\mathcal{C}_{\mathbf{q},\mathbf{q}^\vee}}(\mathscr{H}_3,\mathscr{H}_1\otimes \mathscr{H}_2)
$$
and similarly $I_{23}^1$; it is given by
\begin{align}
\label{eq:I_23_1_formula}
(I_{12}^3 f)(v) = (\mathbf{F}_{12}^3)^{-1}( f\otimes v), \qquad \forall f\in M_{12}^3, \quad \forall v\in \mathscr{H}_3,
\end{align}
where $\mathbf{F}_{12}^3 : \mathscr{H}_1 \otimes \mathscr{H}_2 \to M_{12}^3 \otimes \mathscr{H}_3$ is the intertwiner that we built. Denote the canonical vector space map in eq.\eqref{eq:Hom_using_dual} by
$$
J_{12}^3 : \mathrm{Hom}_{\rm HS}(\mathscr{H}_3,\mathscr{H}_1\otimes\mathscr{H}_2) \to \mathscr{H}_1\otimes\mathscr{H}_2\otimes\mathscr{H}_3'.
$$
Restricting to $\mathrm{Hom}_{\mathcal{C}_{\mathbf{q},\mathbf{q}^\vee}}(\mathscr{H}_3,\mathscr{H}_1\otimes\mathscr{H}_2)$ we obtain a map
$$
(J_{12}^3)^{\mathrm{Inv}} : \mathrm{Hom}_{\mathcal{C}_{\mathbf{q},\mathbf{q}^\vee}}(\mathscr{H}_3,\mathscr{H}_1\otimes\mathscr{H}_2) \to \mathrm{Inv}(\mathscr{H}_1\otimes\mathscr{H}_2\otimes\mathscr{H}_3').
$$
Consider the maps
\begin{align*}
\xymatrix{
\mathscr{H}_1\otimes\mathscr{H}_2\otimes\mathscr{H}_3' \ar[r] &
{}'\mathscr{H}_1 \otimes \mathscr{H}_2 \otimes \mathscr{H}_3 \ar[r] &
\mathrm{Hom}_{\rm HS}(\mathscr{H}_1,\mathscr{H}_2\otimes\mathscr{H}_3) \ar[r] &
\mathscr{H}_2\otimes\mathscr{H}_3 \otimes \mathscr{H}_1',
}
\end{align*}
the first arrow being $D^{-1} \otimes 1 \otimes C^{-1}$, while the latter two arrows are canonical vector space maps in eq.\eqref{eq:Hom_using_dual}; in particular, as both ${}' \mathscr{H}_1$ and $\mathscr{H}_1'$ are $\mathcal{B}(\mathscr{H}_1,\mathbb{C})$ as vector spaces, we see that the composition ${}'\mathscr{H}_1 \otimes \mathscr{H}_2 \otimes \mathscr{H}_3 \to \mathscr{H}_2\otimes\mathscr{H}_3 \otimes \mathscr{H}_1'$ of the latter two arrows is just the permutation map $\mathbf{P}_{(132)}$, i.e. sending $\sum \varphi \otimes \psi \otimes \eta \mapsto \sum \psi \otimes \eta \otimes \varphi$. Restricting the composition of these three arrows to the invariant subspace we obtain a map
$$
(A_{12}^3)^\mathrm{Inv}: \mathrm{Inv}(\mathscr{H}_1\otimes\mathscr{H}_2\otimes\mathscr{H}_3') \to \mathrm{Inv}(\mathscr{H}_2\otimes\mathscr{H}_3\otimes\mathscr{H}_1')
$$
which via the maps $J_{12}^3$ and $J_{23}^1$ is carried over to the map
$$
(A_{12}^3)^\mathrm{Hom} : \mathrm{Hom}_{\mathcal{C}_{\mathbf{q},\mathbf{q}^\vee}}(\mathscr{H}_3,\mathscr{H}_1\otimes\mathscr{H}_2) \to \mathrm{Hom}_{\mathcal{C}_{\mathbf{q},\mathbf{q}^\vee}}(\mathscr{H}_1,\mathscr{H}_2\otimes\mathscr{H}_3)
$$
which in turn via the maps $I_{12}^3$ and $I_{23}^1$ is carried over to the sought-for map
$$
\mathbf{A} = \mathbf{A}_{12}^3: M_{12}^3 \to M_{23}^1.
$$
As $M = M_{12}^3$ is realized as $L^2(\mathbb{R}^2, dt\,ds)$ as a Hilbert space, one can explicitly compute what this operator $\mathbf{A}$ on $L^2(\mathbb{R}^2, \, dt\,ds)$ is.

\section{A formal computation of the operator $\mathbf{A}$}
\label{sec:computation_of_A}

Let $f\in M_{12}^3$. Then $I_{12}^3 f \in \mathrm{Hom}_{\mathcal{C}_{\mathbf{q},\mathbf{q}^\vee}}(\mathscr{H}_3,\mathscr{H}_1\otimes\mathscr{H}_2)$. Suppose $\{v_i\}_{i=1}^\infty$ is an orthonormal basis of the Hilbert space $\mathscr{H}$; let $\{v^i\}_{i=1}^\infty$ be the corresponding dual basis of $\mathscr{H}'$. Then from eq.\eqref{eq:I_23_1_formula} and the canonical map in eq.\eqref{eq:Hom_using_dual}, we see that
$$
(J_{12}^3)(I_{12}^3 f) = \sum_{i=1}^\infty \mathbf{F}^{-1}(f\otimes v_i) \otimes v^i \in \mathscr{H}_1\otimes\mathscr{H}_2\otimes \mathscr{H}_3'.
$$
We then apply $D^{-1}_1 \, C^{-1}_3$ to land in ${}'\mathscr{H}_1\otimes \mathscr{H}_2\otimes \mathscr{H}_3$, then $\mathbf{P}_{(132)}$ to land in $\mathscr{H}_2\otimes \mathscr{H}_3\otimes \mathscr{H}_1'$. This last element must coincide with
$$
J_{23}^1(I_{23}^1 \mathbf{A} f) = \sum_{i=1}^\infty \mathbf{F}^{-1}(\mathbf{A}f\otimes v_i)\otimes v^i \in \mathscr{H}_2\otimes\mathscr{H}_3\otimes \mathscr{H}_1'.
$$
So
$$
\textstyle \mathbf{P}_{(132)} D_1^{-1} C_3^{-1} \, \mathbf{F}_{12}^{-1} (\sum_{i=1}^\infty f\otimes v_i \otimes v^i) = \mathbf{F}_{12}^{-1} \mathbf{A}_1 (\sum_{i=1}^\infty f\otimes v_i \otimes v^i),
$$
where $\mathbf{F}_{12}$ means that we are applying $\mathbf{F}$ to the first and the second tensor factors (so $1,2$ are not directly about the indexing $j$ for $\mathscr{H}_j$'s), and $\mathbf{A}_1$ means that we are applying $\mathbf{A}$ to the first tensor factor. Thus, we must compute
\begin{align}
\label{eq:to_compute_for_A}
\textstyle \mathbf{F}_{12} \, \mathbf{P}_{(132)} \, D_1^{-1}\,C_3^{-1}\, \mathbf{F}_{12}^{-1} (f\otimes (\sum_{i=1}^\infty v_i \otimes v^i)), \quad\mbox{which belongs to}\quad M\otimes\mathscr{H}\otimes\mathscr{H}',
\end{align}
and prove that it is in the form $g \otimes (\sum_{i=1}^\infty v_i \otimes v^i)$ for some $g\in \mathscr{H}$; then this $g$ is our sought-for $\mathbf{A} f$. When computing this, we notice the fact that $\sum_i v_i \otimes v^i$ is the `canonical formal element' of $\mathscr{H} \otimes \mathscr{H}' = L^2(\mathbb{R}^2, dt_2\,ds_2) {\otimes} L^2(\mathbb{R}^2,\,dt_3\,ds_3)$ corresponding to $\mathrm{id}\in \mathrm{Hom}_\mathbb{C}(\mathscr{H},\mathscr{H})$, hence can be thought of formally as $\delta(t_2-t_3) \, \delta(s_2-s_3)$, where $\delta$ here means the Direc delta function(-al), viewed as a generalized eigenvector. As the title of the present section suggests, the readers can take this just formally, at the moment. Therefore, when applied to $g\otimes(\sum_i v_i \otimes v^i) \in M\otimes \mathscr{H}\otimes \mathscr{H}'$, we find
$$
q_{t_2} = q_{t_3}, \quad q_{s_2} = q_{s_3}, \quad
p_{t_2} = - p_{t_3}, \quad p_{s_2} = -p_{s_3}
$$
and hence
\begin{align}
\label{eq:operator_equations_applied_to_delta}
P_2 = P_3^*, \quad Q_2 = Q_3^*, \quad P_2^* = P_3, \quad Q_2^* = Q_3.
\end{align}
Such a proof and the computational result for $g = \mathbf{A} f$ could be rigorously obtained by expressing $\mathbf{F}$ as an explicit integral operator when $f$ is a nice test function, and computing the multiple integral expressing the element in eq.\eqref{eq:to_compute_for_A} with the help of some integral identities involving $\Psi^\hbar$, in the style of \cite{FKi}. However, for that one first needs to establish necessary integral identities for this modular double compact quantum dilogarithm function $\Psi^\hbar$ like Ruijsenaars did for the usual non-compact quantum dilogarithm function $\Phi^\hbar$ in \cite{Rui}, whose results were used in \cite{FKi}. We do not provide such a rigorous treatment in the present paper, and only perform formal computations to figure out what our $\mathbf{A}$ is.

\vs

For our formal computation, we assume without proof that eq.\eqref{eq:to_compute_for_A} is indeed of the form $g\otimes (\sum_{i=1}^\infty v_i \otimes v^i)$, which one might be able to show using representation theory of Hopf algebras by establishing that each map that we used preserves the invariant subspaces; then
\begin{align}
\nonumber
\textstyle \mathbf{F}_{12} \, \mathbf{P}_{(132)} \, D_1^{-1}\,C_3^{-1}\, \mathbf{F}_{12}^{-1} \, (f(t_1,s_1) \, \delta(t_2-t_3)\,\delta(s_2-s_3))
= (\mathbf{A} f)(t_1,s_1) \, \delta(t_2-t_3)\,\delta(s_2-s_3),
\end{align}
for nice test functions $f$, say $f\in \mathscr{D}$. We shall then compute how this operator $\mathbf{A} : M = L^2(\mathbb{R}^2, \, dt_1\,ds_1) \to M = L^2(\mathbb{R}^2,\,dt_1\,ds_1)$ intertwines the normal operators $P_1,Q_1,P_1^*,Q_1^*$, i.e. the conjugation action of $\mathbf{A}$ on these operators, which completely determines $\mathbf{A}$ up to a multiplicative constant because of the Stone-von Neumann type irreducibility of these four operators on $L^2(\mathbb{R}^2)$; that is, any bounded operator (strongly) commuting with all these four must be a scalar operator. To do this, we investigate how $\mathbf{A}$ intertwines the exponential operators $e^{P_1}, e^{Q_1}, e^{P_1^*}, e^{Q_1^*}$ and their Langlands modular dual operators (i.e. the ones marked with $\vee$, which are the $1/({\rm i}\hbar)$-th powers); if we assume Conjecture \ref{conj:strong_irreducibility}, knowing conjugation action on these eight exponential operators is as good as knowing that on $P_1,Q_1,P_1^*,Q_1^*$. We show one such computation of the conjugation action on the exponential operators. 
\begin{align*}
& (\mathbf{F}_{12} \, \mathbf{P}_{(132)}\, D_1^{-1}  \, C_3^{-1} \, \mathbf{F}_{12}^{-1}) \, e^{-P_1} \\
& = \mathbf{F}_{12} \, \mathbf{P}_{(132)}\, D_1^{-1} \, C_3^{-1} \, \Psi^\hbar(Q_1+P_2-Q_2)\,\ul{ {\bf S}_{12}^{-1} \, e^{-P_1} } \qquad \mbox{($\because$ eq.\eqref{eq:bf_F2})} \\
& = \mathbf{F}_{12} \, \mathbf{P}_{(132)}\, D_1^{-1} \, C_3^{-1} \, \ul{ \Psi^\hbar(Q_1+P_2-Q_2)\, (e^{-P_1} } \, {\bf S}_{12}^{-1} ) \qquad (\because {\rm eq.}\eqref{eq:S_conjugation_on_PQs}) \\
& = \mathbf{F}_{12} \, \mathbf{P}_{(132)}\, D_1^{-1} \, \ul{ C_3^{-1} \, (e^{-P_1} + e^{-P_1+Q_1+P_2-Q_2}) } \, \ul{ \Psi^\hbar(Q_1+P_2-Q_2) \, {\bf S}_{12}^{-1} } \qquad (\because \mbox{Prop.\ref{prop:operator_version_of_difference_equations_for_Psi}(1)}) \\
& = \mathbf{F}_{12} \, \mathbf{P}_{(132)}\, \ul{ D_1^{-1}  \,(e^{-P_1} + e^{-P_1+Q_1+P_2-Q_2})}\,  C_3^{-1}  ) \, (\mathbf{F}_{12}^{-1}) \qquad \mbox{($\because$ eq.\eqref{eq:bf_F2})} \\
& = \mathbf{F}_{12} \, \ul{ \mathbf{P}_{(132)}  \,(e^{P_1^*} + e^{\pi(-{\rm i} + \hbar) + Q_1^* +P_2-Q_2}) } \,  D_1^{-1} ) \, C_3^{-1} \, \mathbf{F}_{12}^{-1} \\
& \textstyle  \qquad (\because \mbox{the conjugation action of $D$ from the proof of Prop.\ref{prop:weak_isomorphisms_with_duals}}) \\
& = \mathbf{F}_{12} \,(e^{P_3^*} + e^{\pi(-{\rm i}+\hbar) + Q_3^* +P_1-Q_1}) \,  \mathbf{P}_{(132)}   \, D_1^{-1} \, C_3^{-1} \, \mathbf{F}_{12}^{-1} \\
& = {\bf S}_{12} \, \ul{ \Psi^\hbar(Q_1+P_2-Q_2)^{-1}  \,(e^{P_3^*} + e^{\pi(-{\rm i}+\hbar) + Q_3^* +P_1-Q_1})} \, \mathbf{P}_{(132)} \, D_1^{-1} \, C_3^{-1} \, \mathbf{F}_{12}^{-1} \\
&  = \ul{ {\bf S}_{12} \, ( e^{P_3^*} + e^{\pi(-{\rm i}+\hbar) + Q_3^* +P_1-Q_1} + e^{\pi(-{\rm i}+\hbar) + Q_3^*+P_1+P_2-Q_2} )} \, \Psi^\hbar(Q_1+P_2-Q_2)^{-1} ) \\
& \textstyle \quad \cdot \, \mathbf{P}_{(132)} \, D_1^{-1} \, C_3^{-1} \, \mathbf{F}_{12}^{-1} \qquad \mbox{($\because$ Prop.\ref{prop:operator_version_of_difference_equations_for_Psi}(1))} \\
& =
( e^{P_3^*} + e^{\pi(-{\rm i}+\hbar) + Q_3^* +P_1-Q_1-Q_2} + e^{\pi(-{\rm i}+\hbar) + Q_3^*+P_2-Q_2} ) \, \ul{ {\bf S}_{12} \, \Psi^\hbar(Q_1+P_2-Q_2)^{-1} }  \\
& \textstyle \quad \cdot \mathbf{P}_{(132)} \, D_1^{-1} \, C_3^{-1} \, \mathbf{F}_{12}^{-1} \qquad \mbox{($\because$ eq.\eqref{eq:S_conjugation_on_PQs})} \\
& =
( e^{P_3^*} + e^{\pi(-{\rm i}+\hbar) + Q_3^* +P_1-Q_1-Q_2} + e^{\pi(-{\rm i}+\hbar) + Q_3^*+P_2-Q_2} ) \, \mathbf{F}_{12}\,\mathbf{P}_{(132)}\,D_1^{-1}\,C_3^{-1}\,\mathbf{F}_{12}^{-1}.  \quad \mbox{($\because$ eq.\eqref{eq:bf_F2})}
\end{align*}
Applying to a formal element $f(t_1,s_1)\,\delta(t_2-t_3)\,\delta(s_2-s_3)$, we can translate the above result as
\begin{align*}
& (\mathbf{A}(e^{-P_1}f))(t_1,s_1)\, \delta(t_2-t_3)\,\delta(s_2-s_3) \\
& = ( \underbrace{ e^{P_3^*} + e^{\pi(-{\rm i}+\hbar) + Q_3^* +P_1-Q_1-Q_2} + e^{\pi(-{\rm i}+\hbar) + Q_3^*+P_2-Q_2} } )( (\mathbf{A}f)(t_1,s_1) \, \delta(t_2-t_3) \, \delta(s_2-s_3)).
\end{align*}
As $(\mathbf{A}f)(t_1,s_1) \, \delta(t_2-t_3) \, \delta(s_2-s_3)$ is of the form $g\otimes (\sum_i v_i\otimes v^i) \in M\otimes \mathscr{H} \otimes \mathscr{H}'$, one can use eq.\eqref{eq:operator_equations_applied_to_delta} to transform the above underbraced expression; the first term becomes $e^{P_2}$, the second term $e^{\pi(-{\rm i}+\hbar)+P_1-Q_1} \, e^{-Q_2} \, e^{Q_3^*}$ becomes $e^{\pi(-{\rm i}+\hbar) + P_1-Q_1} \, e^{-Q_2} \, e^{Q_2} = e^{\pi(-{\rm i}+\hbar)+P_1-Q_1}$, and the last term \\$e^{\pi(-{\rm i}+\hbar) + P_2-Q_2} \, e^{Q_3^*}$ becomes $e^{\pi(-{\rm i}+\hbar) + P_2-Q_2} \, e^{Q_2} \stackrel{\vee}{=} e^{-\pi\hbar + \pi(-{\rm i}+\hbar) + P_2-Q_2+Q_2} = - e^{P_2}$ (the checked equality follows e.g. from the BCH formula in eq.\eqref{eq:BCH}). So the underbraced expression becomes just $e^{\pi(-{\rm i}+\hbar) + P_1-Q_1}$. We thus formally proved
$$
\mathbf{A} \, e^{-P_1} = e^{\pi(-{\rm i}+\hbar) + P_1 - Q_1} \, \mathbf{A},
$$
which holds as equality when applied to nice test functions in $L^2(\mathbb{R}^2)$; proving this equality rigorously is a straightforward but tedious job. One can also check that similar computations yield
$$
\mathbf{A} \, e^{-P_1^*} = e^{\pi(-{\rm i}-\hbar)+P_1^*-Q_1^*} \, \mathbf{A}.
$$

\vs

For completeness, we also present here the computation for $e^{-Q_1}$. Let us compute the conjugation $(\mathbf{F}_{12} \, \mathbf{P}_{(132)} \, D_1^{-1}\,C_3^{-1}\, \mathbf{F}_{12}^{-1}) \, e^{-Q_1} (\mathbf{F}_{12} \, C_3 \, D_1 \, \mathbf{P}_{(132)}^{-1} \, \mathbf{F}_{12}^{-1})$ step by step; namely, first conjugation by $\mathbf{F}_{12}^{-1}$, then by $C_3^{-1}$, etc. Of course we break $\mathbf{F}_{12}$ into the two factors ${\bf S}_{12} \, (\Psi^\hbar(Q_1+P_2-Q_2))^{-1}$ as in eq.\eqref{eq:bf_F2}.
\begin{align*}
\mbox{initial} & :  e^{-Q_1} \\
\hspace{-3mm} \mbox{conjugation by ${\bf S}_{12}^{-1}$} & : e^{-Q_1+Q_2} \qquad (\because {\rm eq}.\eqref{eq:S_conjugation_on_PQs} )\\
\mbox{by $\Psi^\hbar(Q_1+P_2-Q_2)$} & : e^{-Q_1+Q_2} + e^{P_2} \quad (\because \mbox{Prop.\ref{prop:operator_version_of_difference_equations_for_Psi}(1)}) \\
\mbox{by $C_3^{-1}$} & : e^{-Q_1+Q_2} +e^{P_2} \\
\mbox{by $D_1^{-1}$} & : e^{\pi({\rm i}-\hbar) + P_1^* - Q_1^* + Q_2} + e^{P_2} \qquad (\because \mbox{proof of Prop.\ref{prop:weak_isomorphisms_with_duals}}) \\
\mbox{by $\mathbf{P}_{(132)}$} & : e^{\pi({\rm i}-\hbar) + P_3^* - Q_3^* + Q_1} + e^{P_1} \\
\mbox{by $\Psi^\hbar(Q_1+P_2-Q_2)^{-1}$} & : e^{\pi({\rm i}-\hbar) + P_3^* - Q_3^* + Q_1} + e^{P_1} + e^{P_1+Q_1+P_2-Q_2} \qquad (\because \mbox{Prop.\ref{prop:operator_version_of_difference_equations_for_Psi}(1)}) \\
\mbox{by ${\bf S}_{12}$} & : e^{\pi({\rm i}-\hbar) + P_3^* - Q_3^* + Q_1+Q_2} + e^{P_1} + e^{P_1+Q_1+P_2-Q_2+Q_2-P_1} \quad(\because {\rm eq}.\eqref{eq:S_conjugation_on_PQs}) \\
& = - e^{Q_1} \, e^{-\pi\hbar + P_3^*-Q_3^*} \, \ul{ e^{Q_2} } + e^{P_1} + e^{Q_1} \, \ul{ e^{P_2} } \\
& = - e^{Q_1} \, e^{-\pi\hbar+P_3^*-Q_3^*} \, e^{Q_3^*} + e^{P_1} + e^{Q_1} \, e^{P_3^*} \\
& \quad \textstyle \quad\mbox{on}\quad g\otimes(\sum_{i=1}^\infty v_i \otimes v^i) \quad(\because{\rm eq}.\eqref{eq:operator_equations_applied_to_delta}) \\
& = - e^{Q_1} \, e^{\pi\hbar - \pi\hbar + P_3^* - Q_3^* + Q_3^*} + e^{P_1} + e^{Q_1} \, e^{P_3^*} \qquad (\because \mbox{e.g. by BCH} ~ {\rm eq}.\eqref{eq:BCH}) \\
& = - e^{Q_1} \, e^{P_3^*} + e^{P_1} + e^{Q_1} \, e^{P_3^*} = e^{P_1},
\end{align*}
hence
$$
\mathbf{A} \, e^{-Q_1} = e^{P_1}\, \mathbf{A},
$$
which one can prove in a rigorous manner too. Similarly one obtains
$$
\mathbf{A} \, e^{-Q_1^*} = e^{P_1^*}\, \mathbf{A}.
$$
Readers can carefully check that similar computations go through for the Langlands modular dual operators, in particular using Prop.\ref{prop:operator_version_of_difference_equations_for_Psi}(2):
\begin{align*}
& \mathbf{A} \, e^{-P_1/({\rm i} \hbar)} = e^{( \pi(-{\rm i}+\hbar) + P_1-Q_1)/({\rm i} \hbar)} \, \mathbf{A}, \qquad \mathbf{A} \, e^{P_1^*/({\rm i} \hbar)} = e^{- ( \pi(-{\rm i} - \hbar) + P_1^*-Q_1^*)/({\rm i} \hbar)} \, \mathbf{A} \\
& \mathbf{A} \, e^{-Q_1/({\rm i} \hbar)} = e^{P_1/({\rm i} \hbar)} \, \mathbf{A}, \qquad\qquad\qquad\qquad \mathbf{A} \, e^{Q_1^*/({\rm i} \hbar)} = e^{-P_1^*/({\rm i} \hbar)} \, \mathbf{A}.
\end{align*}
Then, in the spirit of Conjecture \ref{conj:strong_irreducibility}, we claim that we can obtain the following conjugation actions of $\mathbf{A}$ on $P_1,Q_1,P_1^*,Q_1^*$:
\begin{align}
\nonumber
\begin{array}{ll}
\mathbf{A}\, P_1 \, \mathbf{A}^{-1} = \pi({\rm i}-\hbar) - P_1 + Q_1, & 
\mathbf{A}\, P_1^* \, \mathbf{A}^{-1} = \pi({\rm i}+\hbar) - P_1^* + Q_1^*, \\
\mathbf{A} \, Q_1 \, \mathbf{A}^{-1} = - P_1, & 
\mathbf{A}\, Q_1^* \, \mathbf{A}^{-1} = -P_1^*.
\end{array} 
\end{align}
Let us now drop the subscript $1$ and think of $\mathbf{A}$ as being an operator on $L^2(\mathbb{R}^2, \, dt\,ds)$. As $P_1 = P = q_t + {\rm i} \, p_s$, $P_1^*=P^* = q_t - {\rm i}\,p_s$, $Q_1 = Q = q_s - {\rm i}\,p_t$ and $Q_1^*=Q^* = q_s + {\rm i}\,p_t$, one deduces the following conjugation actions on $q_t,p_t,q_s,p_s$:
\begin{align*}
& \mathbf{A} \, q_t \, \mathbf{A}^{-1} = \pi {\rm i} - q_t + q_s, \qquad\qquad
\mathbf{A} \, q_s \, \mathbf{A}^{-1} = -q_t \\
& \mathbf{A} \, p_t \, \mathbf{A}^{-1} = p_s, \qquad\qquad\qquad\qquad\qquad
\mathbf{A} \, p_s \, \mathbf{A}^{-1} = \pi{\rm i}\,\hbar - p_t - p_s,
\end{align*}
which completely determines $\mathbf{A}$ up to multiplicative constant, thanks to the strong irreducibility of $q_t,p_t,q_s,p_s$ on $L^2(\mathbb{R}^2,\,dt\,ds)$. From
\begin{align*}
& e^{-q_s + \hbar^{-1} p_t} \, q_t \, e^{- (-q_s + \hbar^{-1} p_t)} = q_t + \pi{\rm i}, \qquad\qquad\qquad
e^{-q_s + \, \hbar^{-1} p_t} \, q_s \, e^{- (-q_s + \, \hbar^{-1} p_t)} = q_s, \\
& e^{-q_s + \hbar^{-1} p_t} \, p_t \, e^{- (-q_s + \hbar^{-1} p_t)} = p_t, \qquad\qquad\qquad
e^{-q_s + \hbar^{-1} p_t} \, p_s \, e^{- (-q_s + \hbar^{-1} p_t)} = p_s + \pi{\rm i} \hbar,
\end{align*}
and the conjugation actions of $\mathbf{S}_{\smallmattwo{-1}{-1}{1}{0}}$ as in eq.\eqref{eq:saso_conjugation_on_p}--\eqref{eq:saso_conjugation_on_q} , one can easily check that the following is a solution to the conjugation equations for ${\bf A}$:
\begin{proposition}
\label{prop:bf_A_as_saso}
One has
$$
\mathbf{A} = c \cdot \mathbf{S}_{\smallmattwo{-1}{-1}{1}{0}} \, \circ \, e^{-q_s + \hbar^{-1} p_t}
$$
as operators on $L^2(\mathbb{R}^2,\,dt\,ds)$, where $c$ is some complex constant, and $\mathbf{S}_\mathbf{c}$ is as in Def.\ref{def:saso}. \qed
\end{proposition}
We derived the above formula for $\mathbf{A}$ indirectly by computing its conjugation actions, which involved formal computations. As mentioned, one could directly check if the above final computational result for ${\bf A}$ holds true, in the style of \cite{FKi}; then one could also find out what the constant $c$ is. 

\vs

It is a straightforward exercise to check that $\mathbf{S}_{\smallmattwo{-1}{-1}{1}{0}}$ is of order three; just observe that $\smallmattwo{-1}{-1}{1}{0}$ is an element of $\mathrm{SL}_\pm(2,\mathbb{R})$ of order three, and use Prop.\ref{prop:saso}. Observe now that
\begin{align*}
\mathbf{A}^3 & = c^3 \, \mathbf{S}_{\smallmattwo{-1}{-1}{1}{0}} \, e^{-q_s + \hbar^{-1} p_t} \, \mathbf{S}_{\smallmattwo{-1}{-1}{1}{0}} \,  \, e^{-q_s + \hbar^{-1} p_t} \, \mathbf{S}_{\smallmattwo{-1}{-1}{1}{0}} \, e^{-q_s + \hbar^{-1} p_t} \\
& = c^3 \, \left( \mathbf{S}_{\smallmattwo{-1}{-1}{1}{0}} \right)^2 e^{-q_t+q_s + \hbar^{-1} (-p_t-p_s)}  \,  \, e^{-q_s + \hbar^{-1} p_t} \, \mathbf{S}_{\smallmattwo{-1}{-1}{1}{0}} \, e^{-q_s + \hbar^{-1} p_t} \\
& = c^3 \, \left( \mathbf{S}_{\smallmattwo{-1}{-1}{1}{0}} \right)^3 e^{q_t+ \hbar^{-1} p_s}  \,  \, e^{-q_t+q_s + \hbar^{-1} (-p_t-p_s)} \, \, e^{-q_s + \hbar^{-1} p_t}.
\end{align*}
By a repeated use of the BCH formula in eq.\eqref{eq:BCH} one verifies $e^{q_t+ \hbar^{-1} p_s}  \,  \, e^{-q_t+q_s + \hbar^{-1} (-p_t-p_s)} \, \, e^{-q_s + \hbar^{-1} p_t} = e^{\pi{\rm i}\cdot\mathrm{id}} = - \mathrm{id}$. So ${\bf A}^3$ is a scalar operator:
\begin{proposition}[the order three relation for $\mathbf{A}$]
\label{prop:AAA}
One has
$$
\mathbf{A}^3 = - c^3 \cdot \mathrm{id}.
$$
\end{proposition}
Again, this can also be proved directly, not resorting to the BCH formula. Furthermore, we expect that $\mathbf{A}^3 = \mathrm{id}$ would hold, and that almost the same argument as in the representation theoretic proof of Prop.4.16 of \cite{FKi} should work. As a consequence, in particular, one would deduce that $c$ is a third root of $-1$. Finally, we can ponder about the uniqueness of $\mathbf{A}$, with respect to the choice of the $\mathcal{C}_{\mathbf{q},\mathbf{q}^\vee}$-intertwiner $\mathbf{F} : \mathscr{H}{\otimes}\mathscr{H} \to M{\otimes} \mathscr{H}$. As discussed at the end of \S\ref{sec:T}, modulo Conjecture \ref{conj:strong_irreducibility} about the irreducibility, any other intertwiner $\mathbf{F}' : \mathscr{H} {\otimes}\mathscr{H} \to M \otimes\mathscr{H}$ can be written as $\mathbf{F}'_{12} = U_1 \, \mathbf{F}_{12}$ for some unitary operator $U : M \to M$; see eq.\eqref{eq:bf_F_prime}. Note then
$$
\mathbf{F}'_{12} \, \mathbf{P}_{(132)} \, D_1^{-1} \, C_3^{-1} \, (\mathbf{F}_{12}')^{-1} = U_1\, (\mathbf{F}_{12} \, \mathbf{P}_{(132)} \, D_1^{-1} \, C_3^{-1} \, \mathbf{F}_{12}^{-1}) \, U_1^{-1}.
$$
Applied to elements $f\otimes(\sum_{i=1}^\infty v_i \otimes v^i) \in M \otimes \mathscr{H} \otimes \mathscr{H}'$, one finds that this is $U_1 \, \mathbf{A}_1 \, U_1^{-1}$.

\section{New families of representations of the Kashaev group}

Let us summarize what we obtained so far. Let $M = L^2(\mathbb{R}^2, \,dt\,ds)$. From a certain category of representations of the modular double quantum pseudo-K\"ahler plane Hopf algebra $\mathcal{C}_{\mathbf{q},\mathbf{q}^\vee}$, we realized $M$ as the space of intertwiners $\mathrm{Hom}_{\mathcal{C}_{\mathbf{q},\mathbf{q}^\vee}}(\mathscr{H},\mathscr{H}{\otimes}\mathscr{H})$, and constructed a unitary operator (Prop.\ref{prop:construction_of_T})
\begin{align}
\label{eq:T_explicit}
\mathbf{T}_{12} \stackrel{{\rm eq}.\eqref{eq:T_as_F}}{=} \mathbf{F}_{21}^{-1} \stackrel{{\rm eq}.\eqref{eq:bf_F}}{=} \Psi^\hbar(Q_2+P_1-Q_1) \, {\bf S}_{21}^{-1} ~ : ~ M {\otimes} M \to M {\otimes} M,
\end{align}
where $\Psi^\hbar$ is as in eq.\eqref{eq:Psi_hbar_first_time}, $P_j,Q_j$ are as in \S\ref{sec:integrable_representation_of_C_q}, and ${\bf S}_{ij} = e^{\frac{-1}{2\pi\hbar}(P_i Q_j - P_i^* Q_j^*)}$ is as in \S\ref{sec:F}. Also, from this category we constructed a non-unitary operator (Prop.\ref{prop:bf_A_as_saso})
$$
\mathbf{A} = c \cdot \mathbf{S}_{\smallmattwo{-1}{-1}{1}{0}} \, \circ \, e^{-q_s + \hbar^{-1} p_t} ~ : ~ M \to M,
$$
where $\mathbf{S}_\mathbf{c}$ is as in Def.\ref{def:saso}. So far we proved
$$
\mathbf{T}_{23} \, \mathbf{T}_{12} = \mathbf{T}_{12} \mathbf{T}_{13} \mathbf{T}_{23} \quad\mbox{and}\quad
\mathbf{A}^3= -c^3 \cdot \mathrm{id}.
$$

\vs

Mimicking \cite{FKi}, let us define the following one-real-parameter family of operators on $M$:
\begin{align}
\label{eq:bf_A_m}
\mathbf{A}^{(m)} := e^{\pi{\rm i}(m-1)^2/3} \cdot \mathbf{S}_{\smallmattwo{-1}{-1}{1}{0}} \, \circ \, e^{(m-1)(q_s - \hbar^{-1} p_t)} ~ : ~ M \to M, \qquad \forall m\in \mathbb{R}.
\end{align}
In particular, one has
$$
\mathbf{A} = c\, e^{-\pi{\rm i}/3} \cdot \mathbf{A}^{(0)},
$$
while
$$
\mathbf{A}^{(1)} = \mathbf{S}_{\smallmattwo{-1}{-1}{1}{0}}
$$
is a special linear operator (Def.\ref{def:saso}) and hence is unitary (Prop.\ref{prop:saso}), with
$$
\mathbf{A}^{(m)} = e^{\pi{\rm i}(m-1)^2/3} \cdot \mathbf{A}^{(1)} \circ e^{(m-1)(q_s - \hbar^{-1} p_t)}.
$$
By a similar argument as we used when showing ${\bf A}^3 = - c^3\cdot\mathrm{id}$, one can show
$$
(\mathbf{A}^{(m)})^3 = \mathrm{id}, \qquad \forall m\in \mathbb{R}.
$$

\begin{proposition}[two relations involving both ${\bf A}$ and ${\bf T}$]
\label{prop:ATA_ATA_and_TAT_AAP}
For each $m\in \mathbb{R}$ one has
\begin{align}
\label{eq:ATA_ATA}
\mathbf{A}_1^{(m)} \, \mathbf{T}_{12} \, \mathbf{A}_2^{(m)} = \mathbf{A}_2^{(m)} \, \mathbf{T}_{21} \, \mathbf{A}_1^{(m)}
\end{align}
and
\begin{align}
\label{eq:TAT_AAP}
\mathbf{T}_{12} \, \mathbf{A}_1^{(m)} \, \mathbf{T}_{21} = e^{\pi{\rm i}(m-1)^2/3} \, \mathbf{A}_1^{(m)} \, \mathbf{A}_2^{(m)} \, \mathbf{P}_{(12)}.
\end{align}
\end{proposition}

{\it Proof}. When $m=0$, similar arguments like the representation theoretic proofs of Propositions 4.20 and 4.21 of \cite{FKi}, which are mostly just diagram chasing, should work. But here we present a direct proof using the operator expressions of $\mathbf{T}$ and $\mathbf{A}^{(m)}$. From equations \eqref{eq:bf_A_m}, \eqref{eq:saso_conjugation_on_p} and \eqref{eq:saso_conjugation_on_q}, one has
\begin{align}
\label{eq:bf_A_m_conjugation0}
\left\{ \begin{array}{ll}
\mathbf{A}^{(m)} \, q_t \, (\mathbf{A}^{(m)})^{-1} = -(m-1)\pi{\rm i} - q_t + q_s, &
\mathbf{A}^{(m)} \, q_s \, (\mathbf{A}^{(m)})^{-1} = -q_t \\
\mathbf{A}^{(m)} \, p_t \, (\mathbf{A}^{(m)})^{-1} = p_s, &
\mathbf{A}^{(m)} \, p_s \, (\mathbf{A}^{(m)})^{-1} = -(m-1)\pi{\rm i}\,\hbar- p_t - p_s,
\end{array} \right.
\end{align}
which yield
\begin{align}
\label{eq:bf_A_m_conjugation1}
\left\{ \begin{array}{ll}
\mathbf{A}^{(m)} \, P \, (\mathbf{A}^{(m)})^{-1} = (m-1)\pi(-{\rm i}+\hbar)
 - P + Q, & \\
\mathbf{A}^{(m)} \, P^* \, (\mathbf{A}^{(m)})^{-1} = (m-1)\pi(-{\rm i}-\hbar) 
 - P^* + Q^*, & \\
\mathbf{A}^{(m)} \, Q \, (\mathbf{A}^{(m)})^{-1} = - P, &
\mathbf{A}^{(m)} \, Q^* \, (\mathbf{A}^{(m)})^{-1} = -P^*.
\end{array} \right.
\end{align}
Thus
\begin{align*}
& \mathbf{A}_1^{(m)} \, \mathbf{T}_{12} \, (\mathbf{A}_1^{(m)})^{-1}
\stackrel{{\rm eq}.\eqref{eq:T_explicit}}{=} \mathbf{A}_1^{(m)} \, \Psi^\hbar(Q_2+P_1-Q_1)\, e^{\frac{1}{2\pi\hbar}(P_2 Q_1 - P_2^* Q_1^*)} \, (\mathbf{A}_1^{(m)})^{-1} \\
& \stackrel{{\rm eq}.\eqref{eq:bf_A_m_conjugation1}}{=}  \Psi^\hbar(Q_2+(m-1)\pi(-{\rm i}+\hbar)+Q_1) \, e^{\frac{1}{2\pi\hbar}( -P_2P_1 + P_2^*P_1^* )},
\end{align*}
so that
\begin{align*}
& \mathbf{A}_2^{(m)} \, \mathbf{T}_{21} \, (\mathbf{A}_2^{(m)})^{-1}
= \mathbf{P}_{(12)} \, (\mathbf{A}_1^{(m)} \, \mathbf{T}_{12} \, (\mathbf{A}_1^{(m)})^{-1}) \, \mathbf{P}_{(12)} \\
& = \mathbf{P}_{(12)} \, \Psi^\hbar(Q_2+(m-1)\pi(-{\rm i}+\hbar)+Q_1)  e^{\frac{1}{2\pi\hbar}( -P_2P_1 + P_2^*P_1^* )}) \, \mathbf{P}_{(12)} \\
& = \Psi^\hbar(Q_1+(m-1)\pi(-{\rm i}+\hbar)+Q_2) \, e^{\frac{1}{2\pi\hbar}( -P_1P_2 + P_1^*P_2^* )} \\
& = \mathbf{A}_1^{(m)} \, \mathbf{T}_{12} \, (\mathbf{A}_1^{(m)})^{-1},
\end{align*}
as desired in the first equation \eqref{eq:ATA_ATA}.

\vs

For the second equation \eqref{eq:TAT_AAP}, note first that
\begin{align}
\nonumber
\mathbf{T}_{12} & \stackrel{{\rm eq}.\eqref{eq:T_explicit}}{=} \Psi^\hbar(Q_2+P_1-Q_1)^{-1}  \, e^{\frac{1}{2\pi\hbar}(P_2Q_1 - P_2^*Q_1^*)}\\
\label{eq:T_21_explicit}
& \stackrel{\vee}{=} e^{\frac{1}{2\pi\hbar}(P_2Q_1 - P_2^*Q_1^*)}\,   \Psi^\hbar(Q_2+P_1-\cancel{Q_1}+\cancel{Q_1}-P_2),
\end{align}
where the checked equality is from eq.\eqref{eq:S_conjugation_on_PQs}. Note now
\begin{align*}
& \mathbf{T}_{12}\, \mathbf{A}_1^{(m)} \, \mathbf{T}_{21}\, (\mathbf{A}_1^{(m)})^{-1} \\
& \stackrel{{\rm eq}.\eqref{eq:T_explicit},{\rm eq}.\eqref{eq:T_21_explicit}}{=} e^{\frac{1}{2\pi\hbar}(P_2Q_1 - P_2^*Q_1^*)}\, \Psi^\hbar(P_1-P_2+Q_2)\, \mathbf{A}_1^{(m)} \,(  \Psi^\hbar(Q_1+P_2-Q_2) \, e^{\frac{1}{2\pi\hbar}(P_1Q_2 - P_1^*Q_2^*)} ) \, ( \mathbf{A}_1^{(m)})^{-1} \\
& \stackrel{{\rm eq}.\eqref{eq:bf_A_m_conjugation1}}{=}
e^{\frac{1}{2\pi\hbar}(P_2Q_1 - P_2^*Q_1^*)} \, \ul{ \, \Psi^\hbar(P_1-P_2+Q_2) \, \Psi^\hbar(-P_1+P_2-Q_2) } \, ( \mathbf{A}_1^{(m)} \, e^{\frac{1}{2\pi\hbar}(P_1Q_2 - P_1^*Q_2^*)} \, (\mathbf{A}_1^{(m)})^{-1} ) \\
& \stackrel{\vee}{=} 
e^{\frac{1}{2\pi\hbar}(P_2Q_1 - P_2^*Q_1^*)}
\, \exp\left( \frac{ (q_{t_1} - q_{t_2} + q_{s_2})(p_{s_1} - p_{s_2} - p_{t_2}) }{\pi{\rm i} \hbar} \right)
\, ( \mathbf{A}_1^{(m)} \, e^{\frac{1}{2\pi\hbar}(P_1Q_2 - P_1^*Q_2^*)} \, (\mathbf{A}_1^{(m)})^{-1} ),
\end{align*}
where the last checked equality holds from the involutivity property in eq.\eqref{eq:involutivity_Psi} and $P_1-P_2+Q_2 = (q_{t_1} - q_{t_2} + q_{s_2}) + {\rm i} (p_{s_1} - p_{s_2} - p_{t_2})$. Now, let us show that this unitary operator $\exp\left( \frac{ (q_{t_1} - q_{t_2} + q_{s_2})(p_{s_1} - p_{s_2} - p_{t_2}) }{\pi{\rm i}\hbar} \right)$ is in fact a special linear operator (Def.\ref{def:saso}) on $L^2(\mathbb{R}^4, \, dt_1\,ds_1\,dt_2\,ds_2)$ associated to an element of $\mathrm{SL}_\pm(4,\mathbb{R})$. To see this we first make use of the unitary operator $\mathbf{A}^{(1)}_2$, which is a special linear operator, to observe
$$
\mathbf{A}_2^{(1)} \, \exp\left( \frac{ (q_{t_1} - q_{t_2} + q_{s_2})(p_{s_1} - p_{s_2} - p_{t_2}) }{\pi{\rm i}\hbar} \right) \, (\mathbf{A}_2^{(1)})^{-1} \stackrel{{\rm eq}.\eqref{eq:bf_A_m_conjugation0}}{=} \exp\left( \frac{ (q_{t_1} - q_{s_2})(p_{s_1} + p_{t_2}) }{\pi{\rm i}\hbar} \right),
$$
which one can see is a special linear operator too. Indeed, as $q_{t_1},q_{s_2},p_{s_1},p_{t_2}$ mutually strongly commute, one first writes $\exp\left( \frac{ (q_{t_1} - q_{s_2})(p_{s_1} + p_{t_2}) }{\pi{\rm i}\hbar} \right) = \exp\left( \frac{ (q_{t_1} - q_{s_2})p_{s_1} }{\pi{\rm i}\hbar} \right) \exp\left( \frac{ (q_{t_1} - q_{s_2})p_{t_2} }{\pi{\rm i}\hbar} \right) $ and then observes $( \exp\left( \frac{ (q_{t_1} - q_{s_2})p_{s_1} }{\pi{\rm i}\hbar} \right) f)(t_1,s_1,t_2,s_2) = f(t_1,s_1+t_1-s_2,t_2,s_2)$ for all $f\in L^2(\mathbb{R}^4, \, dt_1\,ds_1\,dt_2\,ds_2)$ similarly as done in \S\ref{sec:F}, and likewise $(\exp\left( \frac{ (q_{t_1} - q_{s_2})p_{t_2} }{\pi{\rm i} \hbar} \right)f)(t_1,s_1,t_2,s_2) = f(t_1,s_1,t_2+t_1-s_2,s_2)$. Combining, we get $(\exp\left( \frac{ (q_{t_1} - q_{s_2})(p_{s_1} + p_{t_2}) }{\pi{\rm i} \hbar} \right) f)(t_1,s_1,t_2,s_2) = f(t_1,s_1+t_1-s_2,t_2+t_1-s_2,s_2)$, hence
$$
\exp\left( \frac{ (q_{t_1} - q_{s_2})(p_{s_1} + p_{t_2}) }{\pi{\rm i}\hbar} \right) = \mathbf{S}_{\left(\begin{smallmatrix} 1 & 1 & 1 & 0 \\ 0 & 1 & 0 & 0 \\ 0 & 0 & 1 & 0 \\ 0 & -1 & -1 & 1 \end{smallmatrix} \right)} \quad\mbox{on} \quad L^2(\mathbb{R}^4, \, dt_1\,ds_1\,dt_2\,ds_2).
$$

\vs

Next,
\begin{align*}
\mathbf{A}_1^{(m)} \, e^{\frac{1}{2\pi\hbar}(P_1Q_2 - P_1^*Q_2^*)} \, (\mathbf{A}_1^{(m)})^{-1} 
& \stackrel{{\rm eq}.\eqref{eq:bf_A_m_conjugation1}}{=}
e^{\frac{1}{2\pi\hbar}( ((m-1)\pi(-{\rm i}+\hbar) - P_1 + Q_1) Q_2 - ((m-1)\pi(-{\rm i}-\hbar) - P_1^* + Q_1^*)Q_2^*)} \\
& = e^{\frac{1}{2\pi\hbar}( (-P_1+Q_1)Q_2 - (-P_1^*+Q_1^*)Q_2^*)} \, e^{\frac{m-1}{2\pi\hbar}( \pi(-{\rm i}+\hbar) Q_2 - \pi(-{\rm i}-\hbar) Q_2^* )} \\
& = e^{\frac{1}{2\pi\hbar}( (-P_1+Q_1)Q_2 - (-P_1^*+Q_1^*)Q_2^*)} \,
e^{(m-1)( q_{s_2} - \hbar^{-1}  p_{t_2}  )} \\
& = \mathbf{A}_1^{(1)} \, e^{\frac{1}{2\pi\hbar}(P_1Q_2 - P_1^*Q_2^*)} \, (\mathbf{A}_1^{(1)})^{-1} \, e^{(m-1)( q_{s_2} - \hbar^{-1}  p_{t_2}  )}
\end{align*}

\vs

On the other hand, from Lem.\ref{lem:S_12_as_saso} and eq.\eqref{eq:bf_A_m} we know how to express $e^{\frac{1}{2\pi\hbar}(P_2Q_1-P_2^*Q_1^*)}$, $\mathbf{A}_2^{(1)}$, $\mathbf{A}_1^{(1)}$, $e^{\frac{1}{2\pi\hbar}(P_1Q_2-P_1^*Q_2^*)}$ as special linear operators; a useful formula to use is $(\mathbf{A}^{(1)})^{-1} = \mathbf{S}_{\smallmattwo{0}{1}{-1}{-1}}$ (see Prop.\ref{prop:saso}). We thus have
\begin{align*}
& \mathbf{T}_{12} \, \mathbf{A}_1^{(m)} \, \mathbf{T}_{21} \, (\mathbf{A}_1^{(m)})^{-1} \\
& = \underbrace{ \mathbf{S}_{\left(\begin{smallmatrix} 1 & 0 & 0 & 0 \\ 0 & 1 & 0 & -1 \\ 1 & 0 & 1 & 0 \\ 0 & 0 & 0 & 1 \end{smallmatrix} \right) } }_{ e^{\frac{1}{2\pi\hbar}(P_2Q_1-P_2^*Q_1^*)} }
\,\,\,
\underbrace{ \mathbf{S}_{\left(\begin{smallmatrix} 1 & 0 & 0 & 0 \\ 0 & 1 & 0 & 0 \\ 0 & 0 & 0 & 1 \\ 0 & 0 & -1 & -1 \end{smallmatrix} \right) } }_{ (\mathbf{A}_2^{(1)})^{-1} }
\,\,\,
\underbrace{ \mathbf{S}_{\left(\begin{smallmatrix} 1 & 1 & 1 & 0 \\ 0 & 1 & 0 & 0 \\ 0 & 0 & 1 & 0 \\ 0 & -1 & -1 & 1 \end{smallmatrix} \right) } }_{ \exp\left( \frac{(q_{t_1}-q_{s_2})(p_{s_1}+p_{t_2})}{\pi{\rm i}\hbar} \right) } 
\,\,\,
\underbrace{ \mathbf{S}_{\left(\begin{smallmatrix} 1 & 0 & 0 & 0 \\ 0 & 1 & 0 & 0 \\ 0 & 0 & -1 & -1 \\ 0 & 0 & 1 & 0 \end{smallmatrix} \right) } }_{ \mathbf{A}_2^{(1)} } 
\\
& \quad \cdot 
\underbrace{ \mathbf{S}_{\left(\begin{smallmatrix} -1 & -1 & 0 & 0 \\ 1 & 0 & 0 & 0 \\ 0 & 0 & 1 & 0 \\ 0 & 0 & 0 & 1 \end{smallmatrix} \right) } }_{ \mathbf{A}_1^{(1)} }
\,\,\,
\underbrace{ \mathbf{S}_{\left(\begin{smallmatrix} 1 & 0 & 1 & 0 \\ 0 & 1 & 0 & 0 \\ 0 & 0 & 1 & 0 \\ 0 & -1 & 0 & 1 \end{smallmatrix} \right) } }_{ e^{\frac{1}{2\pi\hbar}(P_1Q_2-P_1^*Q_2^*)} } 
\,\,\,
\underbrace{ \mathbf{S}_{\left(\begin{smallmatrix} 0 & 1 & 0 & 0 \\ -1 & -1 & 0 & 0 \\ 0 & 0 & 1 & 0 \\ 0 & 0 & 0 & 1 \end{smallmatrix} \right) } }_{  (\mathbf{A}_1^{(1)})^{-1} } \, e^{(m-1)(q_{s_2} - \hbar^{-1} p_{t_2})}.
\end{align*}
Then, using the (4-dimensional version of the) multiplicativity property in eq.\eqref{eq:saso_multiplicativity} of special linear operators, it is a straightforward exercise to verify that the above composition of seven special linear operators yields
$$
\mathbf{T}_{12} \, \mathbf{A}_1^{(m)} \, \mathbf{T}_{21} \, (\mathbf{A}_1^{(m)})^{-1} = \mathbf{S}_{\left(\begin{smallmatrix} 0 & 0 & -1 & -1 \\ 0 & 0 & 1 & 0 \\ 1 & 0 & 0 & 0 \\ 0 & 1 & 0 & 0 \end{smallmatrix} \right) } \, e^{(m-1)(q_{s_2} - \hbar^{-1} p_{t_2})};
$$
just multiply the seven elements of $\mathrm{SL}_\pm(4,\mathbb{R})$ in the subscripts. Meanwhile, one observes
\begin{align*}
\mathbf{P}_{(12)} \, \mathbf{A}_2^{(m)} & =  
e^{\pi{\rm i}(m-1)^2/3} \, \underbrace{ \mathbf{S}_{\left(\begin{smallmatrix} 0 & 0 & 1 & 0 \\ 0 & 0 & 0 & 1 \\ 1 & 0 & 0 & 0 \\ 0 & 1 & 0 & 0 \end{smallmatrix} \right) } }_{ \mathbf{P}_{(12)} } \,\,\, \underbrace{ \mathbf{S}_{\left(\begin{smallmatrix} -1 & -1 & 0 & 0 \\ 1 & 0 & 0 & 0 \\ 0 & 0 & 1 & 0 \\ 0 & 0 & 0 & 1 \end{smallmatrix} \right) } }_{ c^{-1}\, \mathbf{A}_2^{(1)} } \, e^{(m-1)(q_{s_2} - \hbar^{-1} p_{t_2})} \\
& = e^{\pi{\rm i}(m-1)^2/3} \, \mathbf{S}_{\left(\begin{smallmatrix} 0 & 0 & -1 & -1 \\ 0 & 0 & 1 & 0 \\ 1 & 0 & 0 & 0 \\ 0 & 1 & 0 & 0 \end{smallmatrix} \right) } \, e^{(m-1)(q_{s_2} - \hbar^{-1} p_{t_2})},
\end{align*}
where the last equality is easily verified by eq.\eqref{eq:saso_multiplicativity}; just multiply the two elements of $\mathrm{SL}_\pm(4,\mathbb{R})$ in the subscripts. So we explicitly showed
$$
\mathbf{T}_{12} \, \mathbf{A}_1^{(m)} \, \mathbf{T}_{21} \, (\mathbf{A}_1^{(m)})^{-1} = e^{\pi{\rm i}(m-1)^2/3}\, \mathbf{P}_{(12)} \, \mathbf{A}_2^{(m)},
$$
which yields $\mathbf{T}_{12} \, \mathbf{A}_1^{(m)} \, \mathbf{T}_{21}  = \mathbf{A}_1^{(m)} \, \mathbf{P}_{(12)} \, \mathbf{A}_1^{(m)} = e^{\pi{\rm i}(m-1)^2/3} \, \mathbf{A}_1^{(m)} \, \mathbf{A}_2^{(m)} \, \mathbf{P}_{(12)}$, as desired in the second equation \eqref{eq:TAT_AAP}. \qed

\vs

\begin{theorem}[the main result of the present paper]
\label{thm:main}
Let $(\mathscr{H},\pi)$ be the integrable representation (as studied in \S\ref{sec:integral_representation_of_the_modular_double}) of the modular double of the quantum pseudo-K\"ahler plane Hopf algebra $\mathcal{C}_{\mathbf{q},\mathbf{q}^\vee}$ (defined in \S\ref{sec:quantum_pK_plane_algebra}). Let $M \equiv L^2(\mathbb{R}^2, \, dt\,ds)$ be the trivial representation of $\mathcal{C}_{\mathbf{q},\mathbf{q}^\vee}$, and let $\mathbf{F} : \mathscr{H}{\otimes} \mathscr{H} \to M {\otimes} \mathscr{H}$ be the unitary $\mathcal{C}_{\mathbf{q},\mathbf{q}^\vee}$-intertwining map constructed in \S\ref{sec:F}.

\vs

Via the construction in \S\ref{sec:T}, the identity $\mathcal{C}_{\mathbf{q},\mathbf{q}^\vee}$-intertwiner $(\mathscr{H}{\otimes}\mathscr{H}){\otimes}\mathscr{H} \to \mathscr{H}{\otimes}(\mathscr{H}{\otimes}\mathscr{H})$ is encoded, with the help of $\mathbf{F}$, as the the unitary operator $\mathbf{T} : M {\otimes} M \to M{\otimes}M$; see Prop.\ref{prop:construction_of_T}, and also eq.\eqref{eq:T_explicit} for an explicit formula. Via the construction of \S\ref{sec:A}--\ref{sec:computation_of_A},  a certain natural map $\mathrm{Hom}_{\mathcal{C}_{\mathbf{q},\mathbf{q}^\vee}}(\mathscr{H}_3,\mathscr{H}_1{\otimes}\mathscr{H}_2) \to \mathrm{Hom}_{\mathcal{C}_{\mathbf{q},\mathbf{q}^\vee}}(\mathscr{H}_1,\mathscr{H}_2{\otimes}\mathscr{H}_3)$ is encoded, with the help of $\mathbf{F}$, as the non-unitary operator $\mathbf{A} : M \to M$; see Prop.\ref{prop:bf_A_as_saso}. If the intertwiner $\mathbf{F}$ is replaced by another intertwiner $\mathbf{F}' = (U\otimes 1) \circ \mathbf{F}$ for a unitary operator $U : M \to M$, then the corresponding $\mathbf{T}'$ and $\mathbf{A}'$ are related to the previous $\mathbf{T}$ and $\mathbf{A}$ by $\mathbf{T}' = (U\otimes U) \, \mathbf{T} \, (U\otimes U)^{-1}$ and $\mathbf{A}' = U \, \mathbf{A} \, U^{-1}$.

\vs

Define the family of densely-defined operators $\mathbf{A}^{(m)} : M \to M$ by the formula in eq.\eqref{eq:bf_A_m}, so that $\mathbf{A} = c\,e^{-\pi{\rm i}/3}\, \mathbf{A}^{(0)}$ (for some fixed complex number denoted by $c$) and that $\mathbf{A}^{(1)}$ is unitary. These operators satisfy
\begin{align}
\label{eq:Kashaev_relations_operators}
\left\{  \begin{array}{ll}
(\mathbf{A}^{(m)})^3 = \mathrm{id}, & \,\,
\mathbf{T}_{23} \, \mathbf{T}_{12} = \mathbf{T}_{12}\, \mathbf{T}_{13} \, \mathbf{T}_{23}, \\
\mathbf{A}_2^{(m)} \, \mathbf{T}_{21} \, \mathbf{A}_1^{(m)} = \mathbf{A}_1^{(m)} \, \mathbf{T}_{12} \, \mathbf{A}_2^{(m)}, & \,\,
\mathbf{T}_{12} \, \mathbf{A}_1^{(m)} \, \mathbf{T}_{21} = e^{\pi{\rm i}\frac{(m-1)^2}{3}} \mathbf{A}_1^{(m)} \, \mathbf{A}_2^{(m)} \, \mathbf{P}_{(12)}. \quad \qed
  \end{array} \right.
\end{align}
\end{theorem}

Having a pair of operators $(\mathbf{T},\mathbf{A}^{(m)})$ for a Hilbert space $M$ satisfying the relations in eq.\eqref{eq:Kashaev_relations_operators} can be interpreted as having a projective representation of the {\it Kashaev group} of transformations of dotted ideal triangulations of punctured surfaces, appearing in the Kashaev quantization of Teichm\"uller spaces (\cite{Ka1} \cite{Ka2} \cite{FKi} \cite{K16}), in view of the following definition (see e.g. \cite{FKi}):
\begin{definition}
For an index set $I$, the \ul{\em Kashaev group} associated to $I$ is the group presented by the generators $a_j,t_{jk},p_\sigma$, for $j,k\in I$ ($j\neq k$) and permutations $\sigma$ of $I$, with the relations
$$
a_j^3 = 1, \quad t_{k\ell} t_{jk} = t_{jk} t_{j\ell} t_{k\ell}, \quad
a_j t_{jk} a_k = a_k t_{kj} a_j, \quad
t_{jk} a_j t_{kj} = a_j a_k p_{(j\,k)},
$$
for all $j,k,\ell\in I$ with $j,k,\ell$ being mutually distinct, together with the trivial relations: $p_{\mathrm{id}}=1$, $p_{\sigma_1} \, p_{\sigma_2} = p_{\sigma_1\circ \sigma_2}$, $p_\sigma \, a_j = a_{\sigma(j)} \, p_\sigma$, $p_\sigma \, t_{jk} = t_{\sigma(j)\,\sigma(k)} \, p_\sigma$, $a_j a_k = a_k a_j$, $a_j t_{k\ell} = t_{k\ell} a_j$ (for mutually distinct $j,k,\ell$), $t_{jk} t_{\ell n} = t_{\ell n} t_{jk}$ (for mutually distinct $j,k,\ell,n$).
\end{definition}
Moreover, we believe that representation theoretic arguments as in \cite{FKi} would yield the {\em genuine} (as opposed to projective) version of the relations in \eqref{eq:Kashaev_relations_operators} (i.e. without multiplicative constants) for the pair of operators $(\mathbf{T}, \mathbf{A})$ built from the representation theory of $\mathcal{C}_{\mathbf{q},\mathbf{q}^\vee}$; that is, $(\mathbf{T},\mathbf{A})$ is a genuine representation of the Kashaev groups. This would in particular imply that $c = e^{2\pi{\rm i}/3}$.

\section{Further discussion}

Notice that the pair of operators $(\mathbf{T},\mathbf{A}^{(1)})$ provides \emph{unitary} projective representations of the Kashaev groups. Consequently, we obtain a family of new unitary projective representations on a Hilbert space of the mapping class group of any finite-type punctured surface. In the `universal' case, this leads to a projective representation of the (Ptolemy-)Thompson group $T$, via a similar construction in \cite{K16}. Meanwhile, the unitary intertwining operators for mutations obtained in \cite{KS} also provide new (genuine) representations of mapping class groups and the Thompson group $T$, via a similar construction in \cite{FS}. It would be interesting to investigate whether the representation of $T$ obtained from the operators $(\mathbf{T},\mathbf{A}^{(1)})$ of the present paper coincides (up to constants) with the one obtained from the operators in \cite{KS} for the positive cosmological constant. This would be the 3d quantum gravity analog of the similar phenomenon arising in the quantum Teichm\"uller theory, which was studied in \cite{K16}.

\vs

In the meantime, as mentioned in \S\ref{sec:introduction}, the new representations of the Kashaev groups obtained in the present paper call for and hint to the existence of a Kashaev-type quantization of 3d gravity for positive cosmological constant, as an analog of Kashaev's quantum Teichm\"uller theory. What should be the counterpart of Kashaev's ratio coordinates for the moduli spaces of 3d spacetimes? For this we may first need some (generalized-)complexified version of lambda lengths for 3d gravity. We leave this to a future research.

\vs

Another remaining question is about similar constructions to the present paper for the two other values of the cosmological constant. For the negative cosmological constant, the 3d quantum gravity theory essentially becomes Fock-Goncharov's quantization of symplectic double of Teichm\"uller spaces, as mentioned in \cite{KS}. On the other hand, the Kashaev group representation obtained in \cite{FKi} from the quantum plane algebra is the one coming from Kashaev's quantization of one copy of a Teichm\"uller space, not the symplectic double. Hence one might try to define a `symplectic double' version of the quantum plane algebra, and deduce unitary operators for a `symplectic double' version of Kashaev's quantum Teichm\"uller theory. In particular, each of the $\mathbf{F}$ and $\mathbf{T}$ operator will involve two factors of non-compact quantum dilogarithm function $\Phi^\hbar$ (or $e_b$), not just one. For the zero cosmological constant case, one must first come up with a suitable Hopf algebra to start with. Here is one suggestion for such a Hopf algebra together with a representation, inspired by the quantization result \cite{KS}, which unifies all three cases for the cosmological constant. Recall $\mathscr{H} = L^2(\mathbb{R}^2, dt ds)$ in eq.\eqref{eq:H} and the self-adjoint operators $\wh{p}_j,\wh{q}_j$ in eq.\eqref{eq:p_j_and_q_j} on it, for $j\in \{t,s\}$, satisfying the Weyl-relations version of the Heisenberg relations as in eq.\eqref{eq:p_j_and_q_k_Heisenberg_relations}. For a cosmological constant $\Lambda \in \{-1,0,1\}$, consider the direct sum Hilbert space $\mathscr{H} \oplus \mathscr{H}$ (considered in \cite{KS}), and the densely-defined operators $\wh{z}_1 := \smallmattwo{\wh{q}_t}{ -\Lambda\wh{p}_s}{\wh{p}_s}{\wh{q}_t}$ and $\wh{z}_2 := \smallmattwo{\wh{q}_s}{\Lambda \wh{p}_t}{- \wh{p}_t}{\wh{q}_s}$ on it, which would play a role of the $\log$-generators for a sought-for algebra. Here, an element of $\mathscr{H}\oplus\mathscr{H}$ is in the form of $\smallvectwo{\varphi}{\psi}$ for $\varphi,\psi\in \mathscr{H}$, where $\smallmattwo{A}{B}{C}{D} \smallvectwo{\varphi}{\psi} = \smallvectwo{A\varphi + B\psi}{C\varphi + D\psi}$ for operators $A,B,C,D$ on $\mathscr{H}$. Then $[\wh{z}_1,\wh{z}_2] = 2\pi{\rm i}\,\hbar\cdot\hat{\ell}$, where $\hat{\ell} := \smallmattwo{0}{-\Lambda}{1}{0}$ is `like' a scalar satisfying $\hat{\ell}^{\, 2} = - \Lambda \cdot \mathrm{id}$. One suggestion for a Hopf algebra is an algebra generated by $Z_1^{\pm 1},Z_2$ over the `complexified ring of generalized complex numbers' $\mathbb{C}_\Lambda := \mathbb{C}[\ell]/(\ell^2=-\Lambda)$ mod out by the relations $Z_1 Z_2 = e^{2\pi{\rm i}\,\hbar\,\ell} Z_2 Z_1$, whose Hopf algebra structure comes from the coproduct $\Delta Z_1 = Z_1 \otimes Z_1$, $\Delta Z_2 = Z_2 \otimes Z_1 + 1\otimes Z_2$ and the antipode $S(Z_1) = Z_1^{-1}$, $S(Z_2) = -Z_2 Z_1$, with the representation on $\mathscr{H} \oplus \mathscr{H}$ given as $Z_1 \mapsto e^{\wh{z}_1}$, $Z_2 \mapsto e^{\wh{z}_2}$, $\ell \mapsto \wh{\ell}$. See \cite{KS}, in order to seek a motivation for this suggestion, and for more details and subtleties that are omitted here.

\vs

We finish the paper by presenting the formulas for the above-mentioned negative cosmological constant case, without proof. Let $q = e^{\pi{\rm i}\hbar}$ and $q^\vee = e^{\pi{\rm i}/\hbar}$. Consider  $\mathcal{D}_{q,q^\vee}$, defined as the algebra over $\mathbb{Z}[q^{\pm1}, (q^\vee)^{\pm 1}]$ generated by
$$
X, ~ Y, ~ \til{X}, ~ \til{Y}, ~ X^\vee, ~ Y^\vee, ~ \til{X}^\vee, ~ \til{Y}^\vee, \quad\mbox{and the inverses of $X,\til{X},X^\vee,\til{X}^\vee$},
$$
mod out by the relations
\begin{align}
\label{eq:D_q_q_vee_relation}
\left\{
\begin{array}{l}
X \, Y = q^2\, Y\, X, \quad \til{X} \, \til{Y} = q^{-2} \, \til{Y} \, \til{X}, \quad
X^\vee \, Y^\vee = (q^\vee)^2 \, Y^\vee \, X^\vee, \quad \til{X}^\vee \, \til{Y}^\vee = (q^\vee)^{-2} \, \til{Y}^\vee\, \til{X}^\vee, \\ 
 \mbox{each of $X,Y$ (strongly) commutes with each of $\til{X},\til{Y}$,} \\
  \mbox{each of $X^\vee,Y^\vee$ (strongly) commutes with each of $\til{X}^\vee,\til{Y}^\vee$,} \\
 \mbox{each of $X,Y,\til{X},\til{Y}$ (weakly) commutes with each of $X^\vee,Y^\vee,\til{X}^\vee,\til{Y}^\vee$;}
  \end{array} \right.
\end{align}
equipped with the $*$-structure given by
\begin{align}
\label{eq:D_q_q_vee_star_structure}
X^* = X, ~ Y^* = Y, ~ \til{X}^* = \til{X}, ~ \til{Y}^* = \til{Y}, ~ (X^\vee)^* = X^\vee, ~ (Y^\vee)^* = Y^\vee, ~ (\til{X}^\vee)^* = \til{X}^\vee, ~ (\til{Y}^\vee)^* = \til{Y}^\vee,
\end{align}
and the Hopf algebra structure
\begin{align*}
\left\{ \begin{array}{rll}
\mbox{coproduct}: & \Delta(X) = X \otimes X, & \quad \Delta(Y) = Y \otimes X + 1 \otimes Y, \\
& \Delta(\til{X}) = \til{X} \otimes \til{X}, & \quad \Delta(\til{Y}) = \til{Y} \otimes \til{X} +  1 \otimes \til{Y}, \\
& \Delta(X^\vee) = X^\vee \otimes X^\vee, & \quad \Delta(Y^\vee) = Y^\vee \otimes X^\vee + 1\otimes Y^\vee, \\
& \Delta(\til{X}^\vee) = \til{X}^\vee \otimes \til{X}^\vee, & \quad \Delta(\til{Y}^\vee) = \til{Y}^\vee \otimes \til{X}^\vee + 1 \otimes \til{Y}^\vee, \\
\mbox{counit} : & \epsilon(X) = 1= \epsilon(\til{X}), & \quad \epsilon(Y) = 0 = \epsilon(\til{Y}), \\
 & \epsilon(X^\vee) = 1= \epsilon(\til{X}^\vee), & \quad \epsilon(Y^\vee) = 0 = \epsilon(\til{Y}^\vee), \\
\mbox{antipode} : & S(X) = X^{-1}, &  \quad S(Y) = -  Y X^{-1}, \\
& S(\til{X}) = \til{X}^{-1}, &  \quad S(\til{Y}) = -  \til{Y} \til{X}^{-1} , \\
& S(X^\vee) = (X^\vee)^{-1}, & \quad S(Y^\vee) = - Y^\vee (X^\vee)^{-1}, \\
& S(\til{X}^\vee) = (\til{X}^\vee)^{-1}, & \quad S(\til{Y}^\vee) = - \til{Y}^\vee (\til{X}^\vee)^{-1} ,
\end{array} \right.
\end{align*}
compatible with the $*$-structure. We also consider the transcendental relations
\begin{align}
\label{eq:D_transcendental_relation}
X^{1/\hbar} = X^\vee, \qquad \til{X}^{1/\hbar} = \til{X}^\vee, \qquad Y^{1/\hbar} = Y^\vee, \qquad \til{Y}^{1/\hbar} = \til{Y}^\vee.
\end{align}
So this Hopf $*$-algebra $\mathcal{D}_{q,q^\vee}$ could be called the modular double of the `symplectic double of the quantum plane algebra'. One then considers the following positive $*$-representation $\pi$ of $\mathcal{D}_{q,q^\vee}$ on the Hilbert space $\mathscr{H} = L^2(\mathbb{R}^2, dt\,ds)$ in eq.\eqref{eq:H}
\begin{align}
\label{eq:representation_D}
\left\{  \begin{array}{ll}
X \mapsto \pi(X):=\mathbf{X} := e^P, &
Y \mapsto \pi(Y):=\mathbf{Y} := e^Q, \\
\til{X} \mapsto \pi(\til{X}):=\til{\mathbf{X}} := e^{\til{P}}, &
\til{Y} \mapsto \pi(\til{Y}):=\til{\mathbf{Y}} := e^{\til{Q}}, \\
X^\vee \mapsto \pi(X^\vee):=\mathbf{X}^\vee := e^{P/\hbar}, &
Y^\vee \mapsto \pi(Y^\vee):=\mathbf{Y}^\vee := e^{Q/\hbar}, \\
\til{X}^\vee \mapsto \pi(\til{X}^\vee):=\til{\mathbf{X}}^\vee := e^{\til{P}/\hbar}, &
\til{Y}^\vee \mapsto \pi(\til{Y}^\vee):=\til{\mathbf{Y}}^\vee := e^{\til{Q}/\hbar},
  \end{array} \right.
\end{align}
where the self-adjoint operators $P$ and $Q$ are defined as
\begin{align}
P := q_t + p_s, \qquad \til{P} := q_t - p_s, \qquad Q := q_s - p_t,  \qquad \til{Q} := q_s + p_t,
\end{align}
where $p_t = \pi{\rm i}\hbar \, \frac{\partial}{\partial t}$, $p_s = \pi{\rm i}\hbar\,\frac{\partial}{\partial s}$, $q_t = t$, and $q_s=s$ are as in eq.\eqref{eq:p_j_and_q_j} of \S\ref{sec:integrable_representation_of_C_q}. Beware that we are now using the symbols $P$ and $Q$ in a different way from before; in particular, now we have $[P,Q] = 2\pi{\rm i}\, \hbar = - [\til{P},\til{Q}]$ and $[P,\til{P}]=[P,\til{Q}]=[Q,\til{Q}]=[Q,\til{P}]=0$. We then consider the tensor product representation $\pi^{(2)}$ on $\mathscr{H} {\otimes} \mathscr{H} \equiv L^2(\mathbb{R}^4, dt_1\, ds_1\, dt_2\, ds_2)$, where the $\mathcal{D}_{q,q^\vee}$-action is defined via the coproduct, i.e. via the formula eq.\eqref{eq:pi_squared_definition}. Define a unitary operator
\begin{align*}
{\bf F} := {\bf S}_{12} \, \Phi^\hbar(\til{Q}_1 + \til{P}_2 - \til{Q}_2) \, (\Phi^\hbar(Q_1 + P_2 - Q_2))^{-1}
\end{align*}
on $L^2(\mathbb{R}^4, dt_1\, ds_1\, dt_2\, ds_2)$, where the first factor ${\bf S}_{12}$ of the right hand side is the same special linear operator as appeared in Lem.\ref{lem:S_12_as_saso}, the subscripts for $P,Q,\til{P},\til{Q}$ are to be understood accordingly, e.g. $P_1 = q_{t_1} + p_{s_1}$, $\til{Q}_2 = q_{s_2} + p_{t_2}$, etc, and $\Phi^\hbar$ is Faddeev-Kashaev's non-compact quantum dilogarithm as in eq.\eqref{eq:Phi_hbar_first_time} which appears in the Fock-Goncharov quantization of the symplectic double of Teichm\"uller spaces. Using our new symbols $P,Q,\til{P},\til{Q}$, we can express ${\bf S}_{12}$ as
$$
{\bf S}_{12} = e^{\frac{1}{2\pi{\rm i}\,\hbar}(P_1 Q_2 - \til{P}_1 \til{Q}_2)}.
$$
When $M {\otimes} \mathscr{H}$ is realized as $L^2(\mathbb{R}^2, dt_1\, ds_1) {\otimes} L^2(\mathbb{R}^2, dt_2\, ds_2) \cong L^2(\mathbb{R}^4, dt_1\, ds_1\,dt_2\,ds_2)$ with $M \cong L^2(\mathbb{R}^2, dt\,ds)$ a trivial representation of $\mathcal{D}_{q,q^\vee}$, we claim that the intertwining equations as in \eqref{eq:intertwining_equation} hold:
\begin{align}
\nonumber
  \mathbf{F} \, \pi^{(2)}(u) = (1\otimes \pi(u)) \, \mathbf{F}
\end{align}
for all elements $u$ of the algebra $\mathcal{D}_{q,q^\vee}$, which is not difficult to check using the operator versions (\cite{FG09} \cite{G}) of the difference equations $\Phi^\hbar(z+2\pi{\rm i}\,\hbar) = (1+e^{\pi{\rm i}\,\hbar} \, e^z) \, \Phi^\hbar(z)$ and $\Phi^\hbar(z+2\pi{\rm i}) = (1+e^{\pi{\rm i}/\hbar} \, e^{z/\hbar}) \, \Phi^\hbar(z)$. We then follow the recipe of \S\ref{sec:T}, and obtain an analog of Prop.\ref{prop:construction_of_T}, namely
$$
{\bf F}_{23} {\bf F}_{13} {\bf F}_{21}^{-1} {\bf F}_{13}^{-1} = {\bf T}_{12} = {\bf T}\otimes 1 : M {\otimes} M {\otimes} \mathscr{H} \to M {\otimes} M {\otimes} \mathscr{H},
$$
with the unitary map ${\bf T} = {\bf T}_{12}$ on $M{\otimes} M \cong L^2(\mathbb{R}^2, dt_1\,ds_1) {\otimes} L^2(\mathbb{R}^2, dt_2\, ds_2) \cong L^2(\mathbb{R}^4, dt_1\,ds_1\,dt_2\,ds_2)$
$$
{\bf T}_{12} = {\bf F}_{21}^{-1} =  \Phi^\hbar(Q_2 + P_1 - Q_1) \, (\Phi^\hbar(\til{Q}_2 + \til{P}_1 - \til{Q}_1) )^{-1} \,  {\bf S}_{21}^{-1}.
$$
In particular, this operator ${\bf T}$ satisfies the pentagon equation
$$
{\bf T}_{23} {\bf T}_{12} = {\bf T}_{12} {\bf T}_{13} {\bf T}_{23}
$$
as in Prop.\ref{prop:TT_TTT}. One then studies the dual representation to $\mathscr{H}$, and mimics the constructions of \S\ref{sec:dual_representation}--\S\ref{sec:computation_of_A}, to construct an operator ${\bf A}$ on $M$. We claim that the weak isomorphism $D : {}' \mathscr{H} \to \mathscr{H}$ satisfies
\begin{align*}
& D \, (-\til{P}) \, D^{-1} = P, \qquad
D \, ( -\pi{\rm i}(1+\hbar) +  \til{Q}- \til{P} ) \, D^{-1} = Q, \\
& D \, (-P) \, D^{-1} = \til{P}, \qquad
D \, ( \pi{\rm i}(1+\hbar) + Q - P ) \, D^{-1} = \til{Q},
\end{align*}
using which one can compute, similarly as done in \S\ref{sec:computation_of_A}, the conjugation action of $\mathbf{A}$ 
\begin{align*}
& \mathbf{A} \, q_t \, \mathbf{A}^{-1} = - q_t + q_s, \qquad\qquad
\mathbf{A} \, q_s \, \mathbf{A}^{-1} = -q_t \\
& \mathbf{A} \, p_t \, \mathbf{A}^{-1} = p_s, \qquad\qquad\qquad\quad
\mathbf{A} \, p_s \, \mathbf{A}^{-1} = \pi{\rm i}(1+\hbar) - p_t - p_s,
\end{align*}
which implies
$$
\mathbf{A} = c \cdot \mathbf{S}_{\smallmattwo{-1}{-1}{1}{0} } \, \circ \, e^{-(1+\hbar^{-1}) q_s }
$$
for some constant $c$. Then, likewise as we did, consider the family of operators on $M$
$$
{\bf A}^{(m)} := \mathbf{S}_{\smallmattwo{-1}{-1}{1}{0} } \, \circ \, e^{(m-1)(1+\hbar^{-1}) q_s } ~ : ~ M \to M, \qquad \forall m\in \mathbb{R}
$$ 
so that ${\bf A} = c \, {\bf A}^{(0)}$, while ${\bf A}^{(1)}$ is a special linear operator and hence is unitary. Then, straightforward computations of the same sort as in the proof of Prop.\ref{prop:ATA_ATA_and_TAT_AAP} would yield an analog of Prop.\ref{prop:ATA_ATA_and_TAT_AAP}, namely
$$
\mathbf{A}_1^{(m)} \, \mathbf{T}_{12} \, \mathbf{A}_2^{(m)} = \mathbf{A}_2^{(m)} \, \mathbf{T}_{21} \, \mathbf{A}_1^{(m)}, \qquad 
\mathbf{T}_{12} \, \mathbf{A}_1^{(m)} \, \mathbf{T}_{21} =  \mathbf{A}_1^{(m)} \, \mathbf{A}_2^{(m)} \, \mathbf{P}_{(12)}.
$$
In particular, we would obtain a new pair $({\bf T}, {\bf A}^{(1)})$ of operators providing (new) unitary representations of the Kashaev groups, which we expect should be coming from a not-yet-established Kashaev-type (as opposed to Fock-Goncharov) quantization of the symplectic double of Teichm\"uller spaces.


\vspace{-1mm}

\begin{thebibliography}{FKV01}

\vspace{-1,5mm}

\bibitem[B07]{Bai} H. Bai, {\it Quantum Teichm\"uller spaces and Kashaev's $6j$-symbols}, Algebr. Geom. Topol. {\bf 7} (2007), 1541--1560.

\bibitem[B01]{B01} E. W. Barnes, {\it Theory of the double gamma function}, Phil. Trans. Roy. Soc. A {\bf 196} (1901), 265--388.

\bibitem[CP94]{CP} V. Chari and A. Pressley, {\it A Guide to Quantum Groups}, Cambridge Univ. Press, Cambridge, 1994.

\bibitem[CF99]{CF} L. Chekhov and V. V. Fock, {\it A quantum Teichm\"uller space}, Theor. Math. Phys. {\bf 120} (1999), 1245--1259.

\bibitem[CW14]{CW} I. Cohen and E. Wagner, {\it Function algebras on a $2$-dimensional quantum complex plane}, J. Phys.: Conf. Ser. {\bf 563} (2014), 012034

\bibitem[F95]{Faddeev} L. D. Faddeev, {\it Discrete Heisenberg-Weyl group and modular group}, Lett. Math. Phys. {\bf 34} (1995), 249--254.

\bibitem[FK94]{FK} L. D. Faddeev and R. M. Kashaev, {\it Quantum Dilogarithm}, Modern Phys. Lett. A {\bf 9} (1994), 427--434.

\bibitem[FKV01]{FKV} L. D. Faddeev, R. M. Kashaev, and A. Y. Volkov, {\it Strongly coupled quantum discrete Liouville theory, I: Algebraic approach and duality}, Commun. Math. Phys. {\bf 219}(1) (2001), 199--219.

\bibitem[F97]{F} V. V. Fock, {\it Dual Teichm\"uller spaces}, arXiv:hep-th/9702018

\bibitem[FG09]{FG09} V. V. Fock and A. B. Goncharov, {\it The quantum dilogarithm and representations of the quantum cluster varieties}, Invent. Math. \textbf{175} (2009), 223--286.

\bibitem[FK12]{FKi} I. B. Frenkel and H. Kim, {\it Quantum Teichm\"uller space from the quantum plane}, Duke Math. J. {\bf 161}(2) (2012), 305--366.

\bibitem[FS10]{FS} L. Funar and V. Sergiescu, {\it Central extensions of the Ptolemy-Thompson group and quantized Teichm\"uller theory}, J. Topol. {\bf 3} (2010), 29--62.

\bibitem[G08]{G} A. B. Goncharov, ``Pentagon relation for the quantum dilogarithm and quantized $\mathcal{M}^{\rm cyc}_{0,5}$'' in {\it Geometry and Dynamics of Groups and Spaces}, Progr. Math. {\bf 265}, Birkh\"auser, Basel, 2008, pp.415--428.

\bibitem[H13]{H13} B. C. Hall, {\it Quantum Theory for Mathematicians}, Graduate Texts in Math. {\bf 267}, Springer, New York, 2013.

\bibitem[I14]{I14} I. C.H. Ip, {\it On tensor products of positive representations of split real quantum Borel subalgebra $\mathcal{U}_{q\widetilde{q}}(\frak{b}_\mathbb{R})$}, Trans. Amer. Math. Soc. {\bf 370} (2018), 4177-4200.

\bibitem[K94]{Ka0} R. M. Kashaev, {\it Quantum dilogarithm as a 6j symbol}, Mod. Phys. Lett. {\bf A9} (1994), 3757--3768.

\bibitem[K98]{Ka1} R. M. Kashaev, {\it Quantization of Teichm\"uller spaces and the quantum dilogarithm}, Lett. Math. Phys. {\bf 43} (1998), 105--115.

\bibitem[K00]{Ka2} R. M. Kashaev, ``On the spectrum of Dehn twists in quantum Teichm\"uller theory'' in {\it Physics and Combinatorics (Nagoya, 2000)}, World Sci., River Edge, N. J., 2001, pp.63--81.


\bibitem[K16a]{K16} H. Kim, {\it The dilogarithmic central extension of the Ptolemy-Thompson group via the Kashaev quantization}, Adv. Math. {\bf 293} (2016), 305--366.

\bibitem[K16b]{K16b} H. Kim, {\it Ratio coordinates for higher Teichm\"uller spaces}, Math. Z. {\bf 283} (2016), 469--513. 

\bibitem[K19]{Kim_JPAA} H. Kim, {\it Finite dimensional quantum Teichm\"uller space from the quantum torus at root of unity}, J. Pure Appl. Algebra {\bf 223}(3) (2019), 1337--1381.

\bibitem[K21]{Kim_phase} H. Kim, {\it Phase constants in the Fock-Goncharov quantum cluster varieties}, Anal. Math. Phys. {\bf 11} (2021), 2.

\bibitem[KS21]{KS} H. Kim and C. Scarinci, {\it A quantization of moduli spaces of 3-dimensional gravity}, arXiv:2112.13329


\bibitem[NT13]{NT} I. Nidaiev and J. Teschner, {\it On the relation between the modular double of $\mathcal{U}_q(\frak{sl}(2,\mathbb{R}))$ and the quantum Teichm\"uller theory}, arXiv:1302.3454

\bibitem[PT01]{PT} B. Ponsot and J. Teschner, {\it Clebsh-Gordan and Racah-Wigner coefficients for a continuous series of representations of $\mathcal{U}_q(\frak{sl}(2,\mathbb{R}))$}, Commun. Math. Phys. {\bf 224} (2001), 613--655.

\bibitem[PR07]{PR} C. Pr\'ev\^ot and Michael R\"ockner, {\it A Concise Course on Stochastic Partial Differential Equations}, Lect. Notes Math. {\bf 1905}, Springer Berlin, Heidelberg, 2007.


\bibitem[RS80]{RS} M. Reed and B. Simon, {\it Methods of Modern Mathematical Physics. I: Functional Analysis}, revised and enlarged ed., 1980, Academic Press, New York, London, 1972.

\bibitem[R05]{Rui} S. N. M. Ruijsenaars, {\it A unitary joint eigenfunction transform for the $\mathrm{A}\Delta\mathrm{O}$s $\exp(i a_{\pm 1} d/dz) + \exp(2\pi z/a_{\mp})$}, J. Nonlinear Math. Phys. {\bf 12} suppl.2 (2005), 253--294. 

\bibitem[S92]{Sch92} K. Schm\"udgen, {\it Operator Representations of $\mathbb{R}^2_q$}, Publ. Res. Inst. Math. Sci. {\bf 28} (1992), 1029--1061.

\bibitem[S12]{Sch} K. Schm\"udgen, {\it Unbounded Self-adjoint Operators on Hilbert Space}, Graduate Texts in Math. {\bf 265}, Springer, Dordrecht, 2012.

\bibitem[S13]{So} S. B. Sontz, {\it A Reproducing Kernel and Toeplitz Operators in the Quantum Plane}, arXiv:1305.6986

\bibitem[SS85]{SS1} J. Stochel and F. H. Szafraniec, {\it On normal extensions of unbounded operators. I}, J. Operator Theory {\bf 14} (1985), 31--55.

\bibitem[T07]{Te} J. Teschner, ``An analog of a modular functor from quantized Teichm\"uller theory'' in {\it Handbook of Teichm\"uller Theory, vol. I}, IRMA Lect. Math. Theor. Phys. {\bf 11}, Eur. Math. Soc., Z\"urich, 2007, pp.685--760.

\bibitem[W00]{W} S. L. Woronowicz, {\it Quantum exponential function}, Rev. Math. Phys. {\bf 12} no.6 (2000), 873--920.



\end{thebibliography}
\end{document}